\documentclass[12pt]{amsart}
\usepackage[T1]{fontenc}
\usepackage[utf8]{inputenc}
\usepackage{amsmath}
\usepackage{amsfonts}
\usepackage{amssymb}
\usepackage{amsthm}
\usepackage{mathrsfs}
\usepackage{mathtools}
\usepackage{esint}
\usepackage{lmodern}
\usepackage[expansion=false]{microtype}
\usepackage{cite}
\usepackage[hidelinks]{hyperref}
\allowdisplaybreaks[4]
\numberwithin{equation}{section}
\newtheorem{theorem}{Theorem}[section]
\newtheorem{lemma}[theorem]{Lemma}
\newtheorem{proposition}[theorem]{Proposition}
\newtheorem{corollary}[theorem]{Corollary}
\theoremstyle{definition}

\theoremstyle{remark}

\hypersetup{
 pdftitle={Monotone Solutions of the Anisotropic Allen-Cahn Equation in Four Dimensions},
 pdfsubject={Unpublished mathematical manuscript},
 pdfkeywords={anisotropic Allen-Cahn equation, monotone solutions, minimizing graphs, varifolds}
}

\begin{document}
\title[The anisotropic De Giorgi conjecture in $\mathbb{R}^4$]{The anisotropic De Giorgi conjecture in $\mathbb{R}^4$}

\author{Lu Chen}
\address[Lu Chen]{Key Laboratory of Algebraic Lie Theory and Analysis of Ministry of Education, School of Mathematics and Statistics, Beijing Institute of Technology, Beijing
100081, PR China;  Tangshan Research Institute, Beijing Institute of Technology, Tangshan 063000, PR China}
\email{chenlu5818804@163.com}

\author{Jiali Lan}
\address[Jiali Lan]{Key Laboratory of Algebraic Lie Theory and Analysis of Ministry of Education, School of Mathematics and Statistics, Beijing Institute of Technology, Beijing
100081, PR China}
\email{17636268505@163.com}

\author{Haolin Liu}
\address[Haolin Liu]{Key Laboratory of Algebraic Lie Theory and Analysis of Ministry of Education, School of Mathematics and Statistics, Beijing Institute of Technology, Beijing
100081, PR China}
\email{haolinliu2023@126.com}

\address{}

\thanks{The first author was partly supported by the National Natural Science Foundation of China (No.
12622108), a grant from Beijing Institute of Technology (No.2022CX01002). }

\date{}

\begin{abstract}
We establish the anisotropic De Giorgi conjecture in \(\mathbb{R}^4\) under the limiting assumption. More precisely,   every entire solution of the anisotropic Allen–Cahn equation satisfying \(u_{x_4}>0\) and \(u(x,x_4)\to\pm1\) as \(x_4\to\pm\infty\) for each \(x\in\mathbb{R}^3\) is one-dimensional. We show that every nontrivial generalized subgraph globally minimizing the anisotropic perimeter in \(\mathbb{R}^4\) is a halfspace.  Combining an anisotropic tilt-excess decay estimate with a Liouville argument, we then deduce one-dimensional symmetry.
\end{abstract}
\maketitle
\markboth{Monotone anisotropic Allen--Cahn solutions}{Monotone anisotropic Allen--Cahn solutions}

\section{Introduction and Main Results}\label{section:introduction}
This paper is devoted to  discussing  the De Giorgi conjecture for the anisotropic Allen-Cahn equation
\begin{equation}\label{equation:mainPDE}
 \operatorname{div}\left(F(Du)DF(Du)\right)=u^3-u\quad\text{in }\mathbb{R}^4,
\end{equation}
where   $F:\mathbb{R}^{4}\to[0,+\infty)$ is non-negative, convex, and positively homogeneous function of degree one (see Section 2).  The density $F$ geometrically represents a Minkowski norm on $\mathbb{R}^{n+1}$ and physically models the direction-dependent surface tension of an interface.
When $F(p)=|p|$, equation \eqref{equation:mainPDE} reduces to the classical Allen-Cahn equation.  In 1978, De Giorgi formulated the following conjecture:\\
\textbf{Conjecture.} Let $u\in C^2(\mathbb{R}^n)$ be a solution of
\begin{equation}\label{e0.73}\Delta u=u^3-u\end{equation}
 such that
 \[|u|\leq 1,\quad \quad \frac{\partial u}{\partial x_n}>0\]
 in the whole $\mathbb{R}^n$.  Is it true that all level sets $\{u=\lambda\}$ of $u$ are hyperplanes, at least if $n\leq 8$?



The conjecture was proved for $n = 2$  by Ghoussoub and Gui in \cite{GG} and
for  $n=3$ by Ambrosio and Cabre in \cite{AC}.    The proofs for $n = 2$ and $3$ use some techniques in the
linear theory developed by Berestycki, Cafarelli and Nirenberg \cite{BCN} in one of their papers on qualitative properties of solutions of semilinear elliptic equations.

The question of De Giorgi is also connected with the theories of minimal hypersurfaces and phase transitions. The conjecture is sometimes referred to as  "the $\varepsilon$-version of Bernstein problem for minimal graphs". This relation with Bernstein problem is probably the reason why De Giorgi states  "at least if $n\leq8$" in the above quotation.
 Most articles dealing with the question of De Giorgi have also considered the conjecture in a slightly simpler version. It consists of assuming that, in addition,
\begin{equation}\label{e0.72}
\lim_{x_n\to\pm\infty}u(x',x_n)=\pm1\qquad\text{for all}\quad x\in\mathbb{R}^{n-1}.
\end{equation}
Under this condition, Savin~\cite{Savin} proved an “improvement of flatness” theorem for $0$ level sets of local minimizers and  prove the conjecture of De Giorgi for $n\leq8$.
Wang~\cite{Wang} established an improvement of the tilt-excess decay estimate for the Allen–Cahn equation, and use this to give a different proof of Savin’s theorem on the uniform $C^{1,\alpha}$ regularity of flat level sets. This generalizes Allard’s $\varepsilon$-regularity theorem for stationary varifolds
to the setting of Allen–Cahn equations. Counterexamples in $\mathbb{R}^9$ were also constructed by del Pino, Kowalczyk and Wei \cite{dPKW}.

  Recently, significant attention has been focused on the anisotropic Allen-Cahn equation \eqref{equation:mainPDE}.   De Rosa and Pigati \cite{DeRosaPigati}  developed an
anisotropic Allen-Cahn min–max theory and a detailed connection between diffuse critical points and anisotropic minimal hypersurfaces for the surface functional
\[\mathcal{A}_F(M)=\int_{M}F(\nu_M)d\mathcal{H}^{n-1}.\]
The corresponding anisotropic minimal surface equation with Dirichlet boundary conditions was studied in \cite{Chen2026AMS}.
The anisotropic Bernstein problem asks whether critical points of $\mathcal{A}_{F}$  which are
graphs of functions defined on all of    $\mathbb{R}^n$  are necessarily affine functions.    It is known that the answer is positive in dimension $n = 2$ by  Jenkins \cite{J} and in dimension $n = 3$ by   Simon \cite{S}.  Recently, Mooney-Yang \cite{MY} constructed nontrivial entire anisotropic minimal graphs for $n=4$, demonstrating  that such non-parametric examples can exhibit genuinely nonlinear behavior in low dimensions.
Moreover, Chen-Lan \cite{Chen2026C2} established a $C^2$-perturbative Bernstein theorem for $1\leq n\leq7$.
More recently, Lu-Wang-Wei-Wu \cite{LWWW} constructed a
smooth solution of the   anisotropic Allen-Cahn equation for $n\geq 5$ that is monotone in one direction but is not one-dimensional.  The anisotropic Bernstein conjecture suggests
critical dimensions that differ from the isotropic case. Specifically, one expects the anisotropic De Giorgi conjecture holds at least for $n\leq 4$.



Throughout this paper, we assume that $F:\mathbb{R}^4\to[0,\infty)$ is even, positive away from the origin, positively homogeneous of degree one, and of class $C^{4,\alpha}(\mathbb{R}^4\setminus\{0\})$, where $0<\alpha<1$.
The potential $W\in C^{3,\alpha}([-1,1])$ satisfies
\begin{equation}\label{equation:W}
 W(\pm1)=W'(\pm1)=0,\qquad W(s)>0\quad(-1<s<1),
 \qquad W''(\pm1)>0.
\end{equation}

The following are our results for the anisotropic De Giorgi conjecture in $n=4$:

\begin{theorem}\label{theorem:main}
Let $u\in C^2(\mathbb{R}^4)$ be a solution of \eqref{equation:mainPDE}
\begin{equation}\label{equation:layer}
\frac{\partial u}{\partial x_4}>0,\qquad
 \lim_{x_4\to\pm\infty}u(x,x_4)=\pm1  \quad\text{for every }x\in\mathbb{R}^3.
\end{equation}
Then $u$ is one-dimensional.
\end{theorem}



We briefly explain the main steps. Using monotonicity and the
limits at infinity assumption \eqref{equation:layer}, we show that \(u\) is a global minimizer.
The energy bound and \(\Gamma\)-convergence imply that every
blow-down sequence has a subsequence converging in
\(L^1_{\mathrm{loc}}\) to
\(\mathbf{1}_{E^c}-\mathbf{1}_E\), where \(E\) is a global
\(F\)-perimeter minimizer. Monotonicity implies that \(E\) is
a generalized subgraph, and the density estimates exclude
the trivial limits. We prove that every nontrivial globally
minimizing generalized subgraph in \(\mathbb{R}^4\) is a
halfspace. Thus the limiting interface is a hyperplane.

We adapt Wang's method~\cite{Wang} to establish an anisotropic
tilt-excess decay estimate with an \(O(r^{-2})\) error at
scale \(r\). Minimality and exact energy convergence provide
the required energy control on smaller cylinders, without
using an anisotropic monotonicity formula. Starting from
the planar blow-down limit, we iterate this estimate and
sum the changes of normal to obtain a fixed
\(3\)-dimensional subspace \(P\) such that
\[
 \int_{B_R}|\Pi_PDu|^2\,dx\le CR
\]
for all sufficiently large \(R\). The positive function
\(u_{x_4}\) and each \(u_\tau\), \(\tau\in P\), solve the
same linearized equation. The quotient argument of
Ambrosio--Cabr\'e~\cite{AC}, together with the above growth
estimate, shows that \(u_\tau/u_{x_4}\) is constant.
The resulting translation invariance and the growth
estimate force these constants to vanish. Hence
\(u_\tau\equiv0\) for every \(\tau\in P\), and \(u\) is
one-dimensional.

We note that Liu--Luo--Wang--Wei--Wei--Wu~\cite{LLWWWW}
and Chan--Fern\'andez-Real--Figalli--Florit-Simon--Serra~\cite{CFFFS}
recently and independently  proved that every bounded stable entire solution
of the Allen--Cahn equation \eqref{e0.73} in \(\mathbb{R}^3\)
is one-dimensional. Consequently, the classical De Giorgi
conjecture in \(\mathbb{R}^4\) holds without the limiting
assumption \eqref{e0.72}. Motivated by this result, we will study the anisotropic
De Giorgi conjecture without the limiting assumption
\eqref{e0.72} in future work.

This paper is organized as follows. In Section~\ref{section:foundations}, we obtain the energy and gradient estimates.   In Section~\ref{section:compactness}, we prove the convergence of minimizing solutions and  establish the associated varifold identity. In Section~\ref{section:sharpgraph}, we classify the limiting minimizing subgraphs. In Sections~\ref{section:normalization} and~\ref{section:flux}, we present the graph approximation and the convergence estimates. Then in Section~\ref{section:onestep}, we prove an  improvement of flatness result. Finally, in Section~\ref{section:iteration}, we complete the proof of Theorem~\ref{theorem:main}.

\section{Preliminary Estimates}\label{section:foundations}

We write $X=(x,x_4)\in\mathbb{R}^3\times\mathbb{R}$. Let $ B'_r(x)$ be an open ball in $\mathbb{R}^3$ and  $B_r(X)$ an open ball in $\mathbb{R}^4$. We denote by $Q_r(x)=B'_r(x)\times(-r,r)\subset\mathbb{R}^4$ the finite cylinder over $B_r'(x)$. Let $P$ denote a hyperplane in $\mathbb{R}^4$, and $\Pi_P$ the orthogonal projection onto it. We also let  $G(3)$  denote  the Grassmann manifold of unoriented $3$-dimensional hyperplanes in $\mathbb{R}^{4}$.

Throughout this paper, let $F:\mathbb{R}^{4}\to[0,+\infty)$ be a nonnegative convex function of class $C^4(\mathbb{R}^{4}\setminus\{0\})$, which is even and positively homogeneous of degree 1, so that
\begin{equation}\label{e0.13}F(t\xi)=|t|F(\xi)\quad\text{for any }t\in\mathbb{R}, \xi\in\mathbb{R}^{4}.\end{equation}
Furthermore, assume that $F$ satisfies
\begin{equation}\label{e0.14}\alpha|\xi|\leq F(\xi)\leq\beta|\xi|,\quad\forall\xi\in\mathbb{R}^4\end{equation}
for some positive constant $\alpha\leq\beta$. Because of \eqref{e0.13} we can assume, without loss of generality, that the convex closed set
\[K=\{x\in\mathbb{R}^{n+1}: F(x)\leq 1\}\]
has  measure $|K|$ equal to the measure $w_{4}$ of the unit ball in $\mathbb{R}^{4}$. Sometimes, we say that $F$ is the  gauge of $K$. The dual metric of $F$ is
\[F^o(x)=\sup_{\xi\in K}\left<\xi,x\right>.\]
It is easy to verify that $F^o(x):\mathbb{R}^{4}\to[0,+\infty)$ is also a convex, positively homogeneous function of degree one. In fact, $F^o$ and $F$ are polar to each other in the sense that:
\[F^o(x)=\sup_{\xi\neq0}\frac{\left<x,\xi\right>}{F(\xi)}\qquad F(\xi)=\sup_{x\neq0}\frac{\left<x,\xi\right>}{F^o(x)}.\]
It is clear that $F^o$ is the gauge function of the set
\[K^o=\{x\in\mathbb{R}^{4}: F^o(x)\leq 1\}.\]
We say that $K^o$ and $K$ are polar to each other, and we denote the measure of $K^o$ by $\kappa_{4}$.  For further details, we refer to the literature \cite{AFT,LA,R}.

\medskip
Let $G\subset\mathbb{R}^{4}$ be an open set. The  total variation of a function $u\in BV(G)$ with respect to a gauge function $F$ is given by
\[\int_{G}|Du|_{F}dx=\sup\left\{\int_{G}u\mathrm{div}\sigma dx: \sigma\in C_0^1(G,\mathbb{R}^{n+1}),\  F^o(\sigma)\leq 1\right\}.\]
This yields the "generalized" definition of perimeter of a set $E$ with respect to $F$:
\[P_{F}(E,G)=\int_{G}|D\chi_{E}|_{F}dx=\sup\left\{\int_{E}\mathrm{div}\sigma dx: \sigma\in C_0^1(G,\mathbb{R}^{4}),\  F^o(\sigma)\leq 1\right\}.\]
The following co-area formula
\[\int_{G}|Du|_{F}dx=\int_{-\infty}^{+\infty}P_{F}(\{u>s\},G)ds,\quad\forall u\in BV(G),\]
and the identity
\[P_{F}(E,G)=\int_{G\cap\partial^* E}F(\nu)d\mathcal{H}^{3}\]
hold, where $\partial^* E$ is the reduced boundary of $E$ and $\nu$ is the outer normal to $E$ (see \cite{AB}). When $G=\mathbb{R}^4$, we  denote $P_{F}(E,G)$ simply by $P_{F}(E)$. If $F(\xi)=|\xi|$, the anisotropic perimeter  reduces to the classical perimeter, denoted by $P(E,G)$, and simply $P_{F}(E)$ when $G=\mathbb{R}^4$.

\medskip

For convenience, we set
\begin{equation}\label{equation:Phi}
 \Phi(p)=\frac{1}{2}F(p)^2,\qquad A(p)=D\Phi(p),\qquad A(0)=0.
\end{equation}
 We assume that there are constants $0<\lambda\le\Lambda$ such that
\begin{equation}\label{equation:ellipticity}
 \lambda|\xi|^2\le D^2\Phi(p)[\xi,\xi]
 \le\Lambda|\xi|^2
 \quad\text{for }p\ne0,\ \xi\in\mathbb{R}^4.
\end{equation}
We recall the following lemma:
\begin{lemma}(\cite{WX}) The following two statements about $F$ are equivalent:
\begin{itemize}
\item[(i)] $\operatorname{Hess}(F^2)$ is positive definite in $\mathbb{R}^{4}\setminus\{0\}$;
\item[(ii)] The restriction of
\[F_{\xi\xi}(\xi)=(F_{\xi_i\xi_j}(\xi))_{i,j=1}^{4}\]
to $\xi^{\perp}$ is a positive definite endomorphism $\xi^{\perp}\to \xi^{\perp}$ for all $\xi\in\mathbb{S}^{3}$.
\end{itemize}
\end{lemma}

It follows from the homogeneity of $F$ that  for $p\in\mathbb{R}^4\setminus\{0\}$,
\begin{equation}\label{equation:Euler}
 A(p)\cdot p=2\Phi(p),\quad D^2\Phi(p)p=A(p),\quad\text{and}\quad A(tp)=tA(p)\ \forall t\in\mathbb{R}.
\end{equation}
Let $\gamma(t)=q+t(p-q)$. Since $A(\cdot)$ is locally Lipschitz away from the origin and  continuous at the origin   with $A(0)=0$, we derive that $A(\gamma(t))$ is absolutely continuous on $[0,1]$. Therefore, we have \[
 A(p)-A(q)=\int_{0}^1\frac{d}{dt}A(\gamma(t))
 =
 \int_{0}^{1}
 D^{2}\Phi(\gamma(t)) (p-q)\,\mathrm{d}t.
\]
Then, the lower bound and the  upper bound  in \eqref{equation:ellipticity} give
\begin{equation}\label{equation:fluxglobal}
 (A(p)-A(q))\cdot(p-q)\ge\lambda|p-q|^2,
 \qquad |A(p)-A(q)|\le\Lambda|p-q|.
\end{equation}

\medskip
 Throughout this paper, we consider the  anisotropic Allen-Cahn equation in the general form as
 \begin{equation}\label{e0.4}
 \operatorname{div}D\Phi(Du)=W'(u)\quad\text{in }\mathbb{R}^4,
\end{equation}
where \(W \in C^3(\mathbb{R})\) is a double well potential  satisfying \(W \ge 0\), \(W(\pm 1) = 0\), and \(W > 0\) in \((-1, 1)\). Furthermore, there exists \(\gamma \in (0, 1)\) such that \(W' < 0\) on \((\gamma, 1)\) and \(W' > 0\) on \((-1, -\gamma)\), and \(W'' \ge \kappa > 0\) for all \(|x| \ge \gamma\). A typical example is \(W(u) = (1 - u^2)^2/4\), which gives \eqref{equation:mainPDE}. Without loss of generality, we will always assume $|u|\leq1$.

\medskip
Through a scaling $u_{\varepsilon}(X):=u(\varepsilon^{-1}X)$, we get the singularly perturbed version of the anisotropic Allen-Cahn equation:
\begin{equation}\label{e0.1}
\varepsilon\operatorname{div}D\Phi(Du_{\varepsilon})=\frac{1}{\varepsilon}W'(u_{\varepsilon})\quad\text{in }\mathbb{R}^4.
\end{equation}
This equation arises as the Euler-Lagrange equation of the energy functional (after adding suitable boundary conditions)
\begin{equation}\label{equation:scaledenergy}
 \mathcal{E}_{\varepsilon,F}(u_{\varepsilon};U)=\int_Ue_{\varepsilon,F}(u_{\varepsilon})\,dx\quad\text{with}\quad
 e_{\varepsilon,F}(u_{\varepsilon})=\varepsilon\Phi(Du_{\varepsilon})+\varepsilon^{-1}W(u_{\varepsilon}).
\end{equation}
We say $u_\varepsilon$  is a minimizer (or a minimizing solution) if for every ball $B$ in the definition domain of $u_\varepsilon$,
\[\int_B\varepsilon\Phi(Du_{\varepsilon})+\varepsilon^{-1}W(u_{\varepsilon})\,dx\leq \int_B\varepsilon\Phi(Dv)+\varepsilon^{-1}W(v)\,dx\]
for any $v\in H^1(B)$ satisfying $v=u_\varepsilon$ on $\partial B$. We call $u_\varepsilon$  is a global minimizer  if the energy inequality holds on the entire space $\mathbb{R}^4$, i.e.
 \[\int_{\mathbb{R}^4}\varepsilon\Phi(Du_{\varepsilon})+\varepsilon^{-1}W(u_{\varepsilon})\,dx\leq \int_{\mathbb{R}^4}\varepsilon\Phi(Dv)+\varepsilon^{-1}W(v)\,dx\]
 for any $v\in H_{loc}^1(\mathbb{R}^4)$ such that $v - u_\varepsilon \in H^1_0(\mathbb{R}^4)$.

Under the assumptions on $W$, there exists a one-dimensional solution $g(t)$ defined on $t\in(-\infty,+\infty)$, satisfying
\begin{equation}\label{e0.2}
g''(t)=W'(g(t)),\quad \lim_{t\to\pm\infty}g(t)=\pm1,
\end{equation}
where the convergence rate is exponential. The first integral for \eqref{e0.2} can be written as $g'(t)=\sqrt{2W(g(t))}>0$. Throughout, $c_W$  denotes the constant defined by
\begin{equation}\label{equation:action}   c_W=\int_{-\infty}^{\infty}g'(t)^2dt=\int_{-\infty}^{\infty}\left(\frac{1}{2}g'(t)^2+W(g(t))\right)dt.\end{equation}


\label{section:calibration}
\begin{lemma}\label{proposition:calibration}
Let $u$ be a solution of \eqref{e0.4} satisfying \eqref{equation:layer}. Then $u$ is a global minimizer of $\mathcal{E}_{1,F}$.
\end{lemma}
\begin{proof}

For each fixed $X=(x,x_4)\in\mathbb{R}^4$, we set $U(X,\tau)=u(x,x_4+\tau)$. Then,  the map $\tau\mapsto U(X,\tau)$ is strictly increasing from $-1$ to $1$ and by the implicit function theorem, for any $s\in (-1,1)$, there is a unique $\tau=\tau(X,s)$ such that $U(X,\tau)=s$. We denote
\[
 p(X,s)=D_XU(X,\tau(X,s))\in\mathbb{R}^4,\qquad v(X,s)=U_\tau(X,\tau(X,s))>0.
\]
A direct computation yields
\[
 \frac{\partial p_j}{\partial X_i}=U_{X_iX_j}-  U_{X_j\tau}U_{X_i}/v,
 \qquad \frac{\partial p_j}{\partial s}=U_{X_j\tau}/v\qquad\text{for }1\leq i,j\leq4.
\]
Since $\operatorname{div}\left(A(D_X U(X,\tau))\right)=W'(U(X,\tau))=W'(s)$, applying   Euler's identity \eqref{equation:Euler}, we deduce that
\[\operatorname{div}_{X}A\left(p(X,s)\right)=\sum_{i,j=1}^4\frac{\partial A_i(p)}{\partial p_j}\left(U_{X_iX_j}-  U_{X_j\tau}U_{X_i}/v\right)=W'(s)-A(p)\cdot D_XU_\tau/v\]
and
\[\partial_s\left(\Phi(p)-W(s)\right)=A(p)\cdot D_XU_\tau/v-W'(s).\]
Therefore, we obtain
\begin{equation}\label{equation:calfield}
\operatorname{div}_{(X,s)}\left(A(p(X,s)),\Phi(p(X,s))-W(s)\right)=0.
\end{equation}

Fix a bounded smooth domain $\Omega_R\subset\mathbb{R}^4$ and let
\[
w\in H^1(\Omega_R),
\qquad
w-u\in H_0^1(\Omega_R).
\]
Set $T(s)=\max\{-1,\min\{s,1\}\}$, we have
\[
T(w)-u\in H_0^1(\Omega_R),
\qquad
D(T(w))=\chi_{\{-1<w<1\}}Dw
\quad\text{a.e.},
\]
where $\chi$ denotes the characteristic function of the set $\{-1<w<1\}$.
Moreover, $W\ge0$ on $\mathbb{R}$ and $W(\pm1)=0$
give
\[
\mathcal{E}_{1,F}(T(w);\Omega_R)
\le
\mathcal{E}_{1,F}(w;\Omega_R).
\]
Thus, replacing $w$ by $T(w)$, we may assume
$|w|\le1$ a.e. in $\Omega_R$.

 Since $-1<u(X)<1$, we set
\[
d_R=1-\|u\|_{L^\infty(\Omega_R)}>0,
\qquad
w_k=\left(1-\frac1k\right)w+\frac1k u,
\qquad k\ge2.
\]
Then
\[
w_k-u\in H_0^1(\Omega_R),
\qquad
\|w_k\|_{L^\infty(\Omega_R)}
\le1-\frac{d_R}{k}<1,
\qquad
w_k\longrightarrow w
\quad\text{in }H^1(\Omega_R).
\]
Integrating \eqref{equation:calfield} from $u(X)$
to $w_k(X)$ and applying the Leibniz rule in the
weak sense, we obtain
\[\begin{split}
&\operatorname{div}_X \left( \int_{u(X)}^{w_k(X)} A(p(X,s)) ds \right) - A(p(X,w_k)) \cdot Dw_k + A(p(X,u)) \cdot Du\\
&  + \Phi(p(X,w_k)) - W(w_k) - \Phi(p(X,u)) + W(u) = 0.
\end{split}\]
Integrating this identity over \(\Omega_R\) and applying the divergence theorem, we get
\[\begin{split}
 &- \int_{\Omega_R} A(p(X,w_k)) \cdot Dw_k dX + \int_{\Omega_R} A(p(X,u)) \cdot Du dX\\
 & + \int_{\Omega_R} \left[ \Phi(p(X,w_k)) - W(w_k) - \Phi(p(X,u)) + W(u) \right] dX = 0.
\end{split}\]

Since $\tau(X,u(X))=0$, we have $p(X,u(X))=Du$.
Integration over $\Omega_R$, together with
$A(Du)\cdot Du=2\Phi(Du)$, therefore yields
\[
\int_{\Omega_R}
\bigl[A(p_k)\cdot Dw_k-\Phi(p_k)+W(w_k)\bigr]\,dX
=
\int_{\Omega_R}
\bigl[\Phi(Du)+W(u)\bigr]\,dX.
\]
By convexity of $\Phi$ and Euler's identity,
\begin{align}
\Phi(Dw_k)+W(w_k)
&\ge
\Phi(p_k)+A(p_k)\cdot(Dw_k-p_k)+W(w_k)
\notag\\
&=
A(p_k)\cdot Dw_k-\Phi(p_k)+W(w_k).
\label{equation:calineq}
\end{align}
Consequently, we get
$
\mathcal{E}_{1,F}(u;\Omega_R)
\le
\mathcal{E}_{1,F}(w_k;\Omega_R).
$

Finally, simce $w_k\to w$ and $Dw_k\to Dw$ a.e.,
the quadratic growth of $\Phi$ and the boundedness
of $W$ on $[-1,1]$ imply
\[
0\le \Phi(Dw_k)+W(w_k)
\le C\bigl(1+|Dw|^2+|Du|^2\bigr)
\in L^1(\Omega_R).
\]
Dominated convergence yields
\begin{equation}\label{e0.5}
\int_{\Omega_R}
\bigl[\Phi(Dw)+W(w)\bigr]\,dX
\ge
\int_{\Omega_R}
\bigl[\Phi(Du)+W(u)\bigr]\,dX.
\end{equation}
Therefore, we conclude that $u$ is a global minimizer.
\end{proof}

\begin{lemma}\label{lemma:growth}
If $v$ is a  global minimizer of $\mathcal{E}_{1,F}$ with $|v|\le1$, then
\begin{equation}\label{equation:growth}
 \mathcal{E}_{1,F}(v;B_R(X))\le CR^3\quad \text{for }R\ge1.
\end{equation}
Moreover, the   global minimizer $v_\varepsilon$ of $\mathcal{E}_{\varepsilon,F}$ satisfies
\begin{equation}\label{equation:growth}
 \mathcal{E}_{\varepsilon,F}(v_\varepsilon;B_r(X))\le Cr^3\quad \text{for }r\ge\varepsilon.
\end{equation}
Here $C$ depends only on $F$ and $W$.
\end{lemma}
 \begin{proof}
Let $h$  be a radial Lipschitz  function in $\mathbb{R}^4$ satisfying
\[
h \equiv -1 \text{ on } B_R(X), \qquad h \equiv 1 \text{ on } \mathbb{R}^4 \setminus B_{R+2}(X), \qquad |Dh| \le C_0.
\]
Let $w=\min\{v,h\}$,  the minimality of $v$ implies that
\[\mathcal{E}_{1,F}(v;B_{R+2})\le \mathcal{E}_{1,F}(w;B_{R+2}).\]
 We decompose \(B_{R+2}\) into
\[
B_{R+2} = B_R \cup \bigl(\{v \le h\} \cap (B_{R+2} \setminus B_R)\bigr) \cup \bigl(\{v > h\} \cap (B_{R+2} \setminus B_R)\bigr),
\]
by the construction of \(h\), we obtain
\[\mathcal{E}_{1,F}(v; B_R) + \mathcal{E}_{1,F}\bigl(v; \{v > h\} \cap (B_{R+2} \setminus B_R)\bigr) \le \mathcal{E}_{1,F}\bigl(h; \{v > h\} \cap (B_{R+2} \setminus B_R)\bigr).\]
Since \(|Dh| \le C_0\) and \(|h| \le 1\), we have
\[
\Phi(Dh) + W(h) \le \frac{\Lambda}{2} C_0^2 + \|W\|_{L^\infty([-1,1])} =: C_1.
\]
Therefore, we conclude that
\[\mathcal{E}_{1,F}(h; \{v > h\} \cap (B_{R+2} \setminus B_R)) \le C_1 |B_{R+2} \setminus B_R| \le C_2 R^3\]
 for $R \ge 1$.

Finally, we define
$
\widetilde v(z)=v_\varepsilon(X+\varepsilon z).
$
By scaling, $\widetilde v$ is a global minimizer of
$\mathcal{E}_{1,F}$ with $|\widetilde v|\le1$. Thus,
for $r\ge\varepsilon$,
\[
\begin{aligned}
\mathcal{E}_{\varepsilon,F}(v_\varepsilon;B_r(X))
=\varepsilon^3
  \mathcal{E}_{1,F}(\widetilde v;B_{r/\varepsilon}(0))
\le C\varepsilon^3(r/\varepsilon)^3
=Cr^3.
\end{aligned}
\]
\end{proof}

\begin{lemma}\label{lemma:micro}
 Let $v\in H^1_{\mathrm{loc}}(\mathbb{R}^4)$ be a weak
solution of \eqref{e0.4} with $|v|\le1$.
There exist $\beta\in(0,1)$ and $C>0$ such that
\begin{equation}\label{e0.6}
\sup_{X\in\mathbb{R}^4}
\|v\|_{C^{1,\beta}(B_1(X))}\le C,
\end{equation}
where $C$ is independent of $X$.
Moreover, for every $X_0\in\mathbb{R}^4$ and $R>0$,
if
$\operatorname*{ess\,inf}_{B_R(X_0)}|Dv|>0$,
then $v\in C^{2,\beta}(B_{R/2}(X_0)).$
\end{lemma}

\begin{proof}
For any $X\in\mathbb{R}^4$, we choose $\eta\in C_c^\infty(B_4(X))$ such that   $0\le\eta\le1$, $\eta=1$  in $B_{3}(X)$ and $|D\eta|\leq C$. Testing the equation \eqref{e0.4} by $\eta^2v$, we obtain
\[\int_{B_4(X)}\eta^2 A(Dv)\cdot Dv +2\int_{B_4(X)}\eta v A(Dv)\cdot \nabla \eta=-\int_{B_4}\eta^2vW'(v).\]
From \eqref{equation:fluxglobal},  we have
\[A(Dv)\cdot Dv\geq\lambda|Dv|^2\quad\text{and}\quad |A(Dv)|\leq\Lambda |Dv|.\]
Using the Young's inequality, we get
\[\begin{split}
\lambda\int_{B_4(X)}\eta^2|Dv|^2\leq& 2\Lambda\int_{B_4(X)}\eta|Dv||D\eta|+\|W'\|_{\infty}\int_{B_4(X)}\eta^2\\
\leq&\frac{\lambda}{2}\int_{B_4(X)}\eta^2|Dv|^2+\frac{2\Lambda^2}{\lambda}\int_{B_4(X)}|D \eta|^2+\|W'\|_{\infty}\int_{B_4(X)}\eta^2.
\end{split}\]
Therefore, we derive that
\[\int_{B_3(X)}|Dv|^2\leq\lambda\int_{B_4(X)}\eta^2|Dv|^2\leq C_1\int_{B_4(X)}(|D\eta|^2+\eta^2)\leq C_2,\]
where $C_2>0$ is a constant depending only on $\lambda,\Lambda,\|W\|_{\infty}$.

Define
\[
a(p)=
\begin{cases}
D^2\Phi(p), & p\ne0,\\
\lambda\mathrm{I}, & p=0.
\end{cases}
\]
Then, $\lambda\mathrm{I}\le a\le\Lambda\mathrm{I}$  and for all $p,q\in\mathbb{R}^4$,
\[
A(q)-A(p)
=
\left(\int_0^1 a((1-s)p+sq)\,ds\right)(q-p).
\]
For $e\in S^3$ and $0<|h|<1/4$, set
\[v_h(x)=\frac{v(x+he)-v(x)}{h}.\]
Then, we have
\begin{equation}\label{equation:diffquotient}
\operatorname{div}(B_hDv_h)=c_hv_h
\quad\text{in }B_2,\end{equation}
where
\[B_h(x)=\int_0^1
a\bigl((1-s)Dv(x)+sDv(x+he)\bigr)\,ds\]
and\[  c_h(x)=\int_0^1W''\bigl((1-s)v(x)+sv(x+he)\bigr)\,ds.\]
Moreover,
\begin{equation}\label{equation:diffquotient}\lambda\mathrm{I}\le B_h\le\Lambda\mathrm{I}
\quad\text{and}\quad
|c_h|\le\|W''\|_{L^\infty([-1,1])}.\end{equation}
The Caccioppoli inequality, local boundedness and interior H\"older estimates imply that for some $\beta\in(0,1)$,
\[
\|v_h\|_{H^1(B_{3/2})}
+\|v_h\|_{C^{0,\beta}(B_1)}
\le C\|v_h\|_{L^2(B_2)}
\le C\|Dv\|_{L^2(B_3)}
\le C,
\]
where we used the definition of $v_h$ in the second inequality. Moreover, since $v_h\to\partial_e v$ in $L^2(B_2)$ as $h\to0$, taking $e=e_i$ for $1\le i\le4$, we conclude
\[
v\in H^2_{\mathrm{loc}}(\mathbb{R}^4),
\quad\text{and}\quad
\sup_{X\in\mathbb{R}^4}
\|v\|_{C^{1,\beta}(B_1(X))}\le C.
\]

Finally, suppose $|Dv|\ge c>0$ a.e. in $B_R(X_0)$, then $\lambda\mathrm{I}\le D^2\Phi(Dv)\le\Lambda\mathrm{I}.$ Differentiating both sides of  \eqref{e0.4}  with respect to $x_i$, we obtain
\[
\operatorname{div}(D^2\Phi(Dv)\cdot Dv_i)=W''(v)v_i
\quad\text{in }B_R(X_0),
\]
where $W''(v)v_i\in C^{0,\beta}(B_R(X_0))$.
By the interior Schauder estimates,  \(v_i \in C^{1,\beta}(B_{R/2}(X_0))\) for every \(1\leq i \leq4\), and therefore
\[
v \in C^{2,\beta}(B_{R/2}(X_0)).
\]

\end{proof}

\begin{corollary}\label{corollary:micro}
Let $\varepsilon>0$ and let
$v_\varepsilon\in H^1_{\mathrm{loc}}(\mathbb{R}^4)$ be a weak solution of \eqref{e0.1} with $|v_\varepsilon|\le1$.
For any $X\in\mathbb{R}^4$, define $\widetilde v(z)=v_\varepsilon(X+\varepsilon z)$. Then  \(\widetilde v\) satisfies
\[
\operatorname{div}_zD\Phi(D_z\widetilde v)
=W'(\widetilde v)
\quad\text{in }\mathbb{R}^4,
\qquad |\widetilde v|\le1,
\]
and applying Lemma~\ref{lemma:micro} to \(\widetilde v\), the rescaling gives
\begin{equation}\label{equation:scaledmicro}
\begin{aligned}
&\|v_\varepsilon\|_{L^\infty(B_\varepsilon(X))}
+\varepsilon
 \|Dv_\varepsilon\|_{L^\infty(B_\varepsilon(X))}
+\varepsilon^{1+\beta}
 [Dv_\varepsilon]_{C^{0,\beta}(B_\varepsilon(X))}\\
&\qquad
=\|\widetilde v\|_{L^\infty(B_1)}
+\|D\widetilde v\|_{L^\infty(B_1)}
+[D\widetilde v]_{C^{0,\beta}(B_1)}
\le C,
\end{aligned}
\end{equation}
where $C$ is independent of $X$ and $\varepsilon$.
\end{corollary}

\begin{lemma}\label{lemma:modicaclear}
Let $\varepsilon>0$ and let
$u_\varepsilon\in H^1_{\mathrm{loc}}(\mathbb{R}^4)$
be a weak solution of \eqref{e0.1} with
$|u_\varepsilon|\le1$.
Let $b_w\in(0,1)$ and $\kappa>0$ be such that $W''(s)\ge\kappa$  for $1-b_w\le|s|\le1$,
and set \[\Sigma_\varepsilon=\{X\in Q_2:|u_\varepsilon(X)|\le1-b_w\}.\]
Then, for any $X\in Q_1$,
\begin{equation}\label{equation:clearing}
 e_{\varepsilon,F}(u_\varepsilon)\le C\varepsilon^{-1}
 \left(e^{-c\operatorname{dist}(X,\Sigma_\varepsilon)/\varepsilon}+e^{-c/\varepsilon}\right),
\end{equation}
where $\operatorname{dist}(X,\Sigma)$ denotes the distance from $X$ to $\Sigma$.
\end{lemma}

\begin{proof}
Applying the anisotropic Modica inequality established in \cite[Theorem 5.1]{DeRosaPigati}, we obtain
\begin{equation}\label{equation:Modica}
\varepsilon\Phi(Du_\varepsilon)
\le \varepsilon^{-1}W(u_\varepsilon).
\end{equation}
Hence
\[
e_{\varepsilon,F}(u_\varepsilon)
\le 2\varepsilon^{-1}W(u_\varepsilon).
\]
For $X\in Q_1\cap\Sigma_\varepsilon$, the conclusion
follows from the boundedness of $W$ on $[-1,1]$.

Fix $X\in Q_1\setminus\Sigma_\varepsilon$ and set
\[
d=\operatorname{dist}(X,\Sigma_\varepsilon),
\qquad
d_0=\operatorname{dist}(Q_1,\partial Q_2)>0,
\qquad
r=\frac12\min\{d,d_0\},
\]
where $\operatorname{dist}(X,\varnothing)=+\infty$. Then $\overline{B_r}(X)\subset Q_2\setminus\Sigma_\varepsilon$ and \(|u_\varepsilon|>1-b_w\) in \(B_r(X)\). Without loss of generality, we assume that \(u_\varepsilon>1-b_w\), and set
\[
w=1-u_\varepsilon,
\qquad
c(X)=\int_0^1 W''(1-tw(X))\,dt.
\]
Then \(0\le w\le b_w\), and
\begin{equation}\label{e0.9}
\operatorname{div}D\Phi(Dw)=-\varepsilon^{-2}W'(1-w).
\end{equation}
Since $W'(1)=0$, we have
\[
-W'(1-w)=c(X)w\quad\text{and}\quad \kappa\le c(X)\le K:=\|W''\|_{L^\infty([-1,1])}.
\]
Therefore, by Euler's identity, \eqref{e0.9} can be rewritten as
\begin{equation}\label{e0.7}
\operatorname{div}(\mathcal B Dw)
=\varepsilon^{-2}aw
\quad\text{in }B_r,
\end{equation}
where $\lambda I\le\mathcal B\le\Lambda I$ is defined by
\[
\mathcal B=
\begin{cases}
D^2\Phi(Dw), & Dw\ne0,\\
\lambda I, & Dw=0.
\end{cases}
\]

Multiplying on both sides of \eqref{e0.7} by \(\eta^2w\) with \(\eta\in C_c^\infty(B_r(X))\) and integrating  by parts, we obtain
\[-\int_{B_r(X)}\eta^2\mathcal{B}Dw\cdot Dw-2\int_{B_r(X)}\eta w \mathcal{B}Dw\cdot D\eta=\int_{B_r(X)}c(X)\varepsilon^{-2}\eta^2w^2.\]
Using Young's inequality, we get
\[\begin{split}\left|2\int_{B_r(X)}\eta w \mathcal{B}Dw\cdot D\eta\right|\leq&2\Lambda \int_{B_r(X)}\eta |w| |Dw|| D\eta|\\
\leq&\frac{\lambda}{2}\int_{B_r(X)}\eta^2 |Dw|^2+\frac{2\Lambda^2}{\lambda}\int_{B_r(X)}w^2 |D\eta|^2.\end{split}\]
Therefore, we obtain
\begin{equation}\label{e0.8}
\frac{\lambda}{2}\int_{B_r(X)} \eta^2|Dw|^2
+\frac{\kappa}{\varepsilon^2}\int_{B_r(X)} \eta^2w^2
\le \frac{2\Lambda^2}{\lambda}\int_{B_r(X)} w^2|D\eta|^2.
\end{equation}

For \(r>0\), choose \(\eta\in C_c^\infty(B_r(X))\) such that
\[
\eta=1\ \text{in }B_{r-L\varepsilon}(X),\quad
\eta=0\ \text{in }\mathbb{R}^4\setminus B_r(X),\quad\text{and}\quad
|D\eta|\le (L\varepsilon)^{-1}.
\]
Substituting into the energy inequality \eqref{e0.8} yields
\[
\int_{B_{r-L\varepsilon}(X)}w^2
\le \frac{2\Lambda^2}{\lambda\kappa L^2}
\int_{B_r(X)}w^2.
\]
Fix \(L\) so large that $\theta:=\frac{2\Lambda^2}{\lambda\kappa L^2}\le\frac12.$
Then, iterating over $B_{r_k}(X)$ with $r_k=r-kL\varepsilon,\  k=0,1,\dots$,
we obtain
\[
\int_{B_{r-kL\varepsilon}} w^2 \le \theta^k \int_{B_r} w^2.
\]
Taking $k = \lfloor r/(2L\varepsilon) \rfloor$, so that $r - kL\varepsilon \ge r/2$, yields
\[
\int_{B_{r/2}} w^2 \le e^{-c_0r/(2\varepsilon)} \int_{B_r} w^2,
\]
where $c_0 := -\log\theta/L > 0$.
Since $-\operatorname{div}(\mathcal{B}Dw) \le 0$ and $w \ge 0$, the local boundedness estimate for nonnegative subsolutions gives
\[
w^2 \le Cr^{-4}\int_{B_{r/2}} w^2 \le Cr^{-4}e^{-c_0r/(2\varepsilon)}\int_{B_r} w^2 \le Ce^{-c_0r/(2\varepsilon)}.
\]
Finally, using \(W(u_\varepsilon)\le Cw^2\) and the Modica inequality \eqref{equation:Modica},
we conclude
\[e_{\varepsilon,F}(u_\varepsilon)
\le C\varepsilon^{-1}
\left(e^{-c_0 d/(4\varepsilon)}+e^{-c_0 d_0/(4\varepsilon)}\right)\le C\varepsilon^{-1}
\left(e^{-c\operatorname{dist}(X,\Sigma_\varepsilon)/\varepsilon}+e^{-c/\varepsilon}\right),\]
where $c=\frac{c_0}{4}\min\{1,d_0\}>0$.
\end{proof}

\section{Limits of Minimizing Solutions}\label{section:compactness}

In this section, we study the asymptotic behavior of minimizing solutions. We begin with the following compactness result.

\begin{proposition}\label{proposition:phasecompact}
Let \(\Omega\subset\mathbb{R}^4\) be open and let
\(\varepsilon_k\to0^+\). Let   \(F_k\to F_\infty\) in \(C^{4}(\mathbb{S}^3)\) and let \(u_{\varepsilon_k}\) be local minimizers of \(\mathcal{E}_{\varepsilon_k,F_k}\) with \(|u_{\varepsilon_k}|\le 1\). Then, after passing to a subsequence,
\begin{equation}\label{equation:phaselimit}
 u_{\varepsilon_k}\to\mathbf{1}_{E^c}-\mathbf{1}_E\quad\text{in }L^1_{\mathrm{loc}},
\end{equation}
where $E$ locally minimizes $P_{F_\infty}$.  Moreover, for any relatively compact Lipschitz domain \(U\subset\subset\Omega\) with \(P_{F_\infty}(E;\partial U)=0\),
\begin{equation}\label{e0.11}
\lim_{k\to\infty}\mathcal{E}_{\varepsilon_k,F_k}(u_{\varepsilon_k};U)= c_WP_{F_\infty}(E;U).
\end{equation}
If   \(u_{\varepsilon_k}\) are global minimizers, then for any fixed \(s\in(-1,1)\) and \(b\in(0,1)\), the level sets \(\{u_{\varepsilon_k}=s\}\) and transition layers \(\{|u_{\varepsilon_k}|\le1-b\}\) converge locally  in Hausdorff distance to \(\partial E\).
\end{proposition}

\begin{proof}
For any compact set $V\subset\subset\Omega$, the local comparison in Lemma~\ref{lemma:growth} gives
\[
\sup_k\mathcal{E}_{u_{\varepsilon_k},F_k}(v_k;V)<\infty.
\]
Moreover, since $F_k\to F_\infty$, for some \(\eta_k\to0\), we have
$
(1-\eta_k)\Phi_\infty\le\Phi_k\le(1+\eta_k)\Phi_\infty,
$
and hence
\begin{equation}\label{equation:energycomparison}
(1-\eta_k)\mathcal{E}_{\varepsilon_k,F_\infty}(z;V)
\le\mathcal{E}_k(z;V)
\le(1+\eta_k)\mathcal{E}_{\varepsilon_k,F_\infty}(z;V).
\end{equation}

By the anisotropic $\Gamma$-convergence result \cite[Theorem~4.1]{DeRosaPigati} and the energy comparison \eqref{equation:energycomparison}, after passing to a subsequence, there exists a locally finite perimeter set \(E\) such that \eqref{equation:phaselimit} holds and
\begin{equation}\label{equation:liminf}
c_WP_{F_\infty}(E;V)
\le\liminf_{k\to\infty}\mathcal{E}_k(u_{\varepsilon_k};V).
\end{equation}
Using the Young's inequality, we also have
\begin{equation}\label{equation:YoungBV}
\varepsilon\Phi(Dv)+\varepsilon^{-1}W(v)
\ge\sqrt{2W(v)}F(Dv)
=F(D\Psi(v)).
\end{equation}

Fix a Lipschitz domain \(U\subset\subset\Omega\) with
\(P_{F_\infty}(E;\partial U)=0\), and let \(G\) be a set
of locally finite perimeter such that
\(G\triangle E\subset\subset U\).
Choose smooth domains \(U_\rho\) and \(V\) satisfying
\[
\overline U\subset U_\rho\subset\subset V\subset\subset\Omega,
\qquad
U_\rho\to\overline U
\quad\text{as }\rho\to0.
\]
For fixed \(\rho>0\),  the anisotropic $\Gamma$-convergence result \cite[Theorem~4.1]{DeRosaPigati} implies that there exists
  a recovery sequence
\(z_k\in H^1(U_\rho)\), \(|z_k|\le1\), such that
\[
z_k\longrightarrow\mathbf{1}_{G^c}-\mathbf{1}_G
\quad\text{in }L^1(U_\rho),
\qquad
\mathcal{E}_k(z_k;U_\rho)
\longrightarrow c_WP_{F_\infty}(G;U_\rho).
\]
Since \(G=E\) on \(U_\rho\setminus\overline U\), we obtain
\[
z_k-v_k\to0
\quad\text{in }L^1(U_\rho\setminus\overline U).
\]

Applying Lemma 15 in \cite{FVdensity}, we deduce that for every \(\delta>0\) and  \(k\) sufficiently large,
there exists \(w_k\in H^1(V)\) such that
\[
w_k=z_k\ \text{in }U,\qquad
w_k=u_{\varepsilon_k}\ \text{near }\partial V,
\]
and
\[
\mathcal{E}_k(w_k;V)
\le
\mathcal{E}_k(u_{\varepsilon_k};V\setminus\overline U)
+\mathcal{E}_k(z_k;U_\rho)+\delta.
\]
Furthermore, the minimality of $u_{\varepsilon_k}$ implies that $\mathcal{E}_k(u_{\varepsilon_k};V)\leq \mathcal{E}_k(w_k;V)$, combining this with above inequality, we have
\[
\mathcal{E}_k(u_{\varepsilon_k};U)
\le\mathcal{E}_k(z_k;U_\rho)+\delta.
\]
Since \(P_{F_\infty}(G;\partial U)=0\), letting first
\(k\to\infty\), then \(\delta\to0\) and
\(\rho\to0\), we obtain
\begin{equation}\label{e0.10}\limsup_{k\to\infty}\mathcal{E}_k(u_{\varepsilon_k};U)\le c_WP_{F_\infty}(G;U).\end{equation}
Taking \(G=E\) in \eqref{e0.10} and combining with \eqref{equation:liminf}, we obtain \eqref{e0.11}. Moreover, since
\[c_W P_{F_\infty}(E;U)\leq \liminf_{k\to\infty}\mathcal{E}_k(u_{\varepsilon_k};U)\leq \limsup_{k\to\infty}\mathcal{E}_k(u_{\varepsilon_k};U)\leq c_W P_{F_\infty}(G;U),\]
we deduce that $P_{F_\infty}(E;U)\leq P_{F_\infty}(G;U)$   for every finite perimeter set \(G\) with \(G\triangle E\subset\subset U\). Therefore, \(E\) locally minimizes \(P_{F_\infty}\) in \(\Omega\).

\medskip
Suppose now that the \(u_{\varepsilon_k}\) are global minimizers. Then, by \eqref{equation:energycomparison}, for all sufficiently large $k$ and every $\phi_k\in H^1(V)$ with $\phi_k=u_{\varepsilon_k}$ on $\partial V$,
\[\mathcal{E}_{\varepsilon_k,F_\infty}(v_k;V)\leq\frac{1+\eta_k}{1-\eta_k} \mathcal{E}_{\varepsilon_k,F_\infty}(\phi_k;V)\leq 2  \mathcal{E}_{\varepsilon_k,F_\infty}(\phi_k;V).\]
From Corollary \ref{corollary:micro}  and the anisotropic Modica inequality
\eqref{equation:Modica}, we get
\[
\varepsilon_k\|Du_{\varepsilon_k}\|_{L^\infty(\mathbb{R}^4)}\le C,
\qquad
\varepsilon_k^2\Phi_\infty(Du_{\varepsilon_k})
\le\frac{W(u_{\varepsilon_k})}{1-\eta_k}\le2W(u_{\varepsilon_k}).
\]
Therefore, it follows from \cite[Corollary~2]{FVdensity} that  for any fixed
\(\theta\in(0,1)\), if $|u_{\varepsilon_k}(X_k)|<\theta$ and $X_k\to X$, then  $X\in\partial E.$
Choosing \(\theta>|s|\), respectively \(\theta>1-b\), yields
\begin{equation}\label{e0.12}\limsup_{k\to\infty}(K\cap\{u_{\varepsilon_k}=s\})\subset\partial E\quad\text{and}\quad \limsup_{k\to\infty}(K\cap\{|u_{\varepsilon_k}|\leq 1-b\})\subset\partial E\end{equation}
for every compact set $K$.

Conversely, for \(X\in\partial E\) and \(r>0\), we have
\[
|E\cap B_r(X)|>0,\qquad |E^c\cap B_r(X)|>0.
\]
By \eqref{equation:phaselimit}, for all sufficiently large \(k\),
\[
\inf_{B_r(X)}u_{\varepsilon_k}<s<\sup_{B_r(X)}u_{\varepsilon_k},
\]
which implies that $B_r(X)\cap\{u_{\varepsilon_k}=s\}\ne\varnothing.$
A finite covering argument then shows that for every compact \(K\subset\partial E\) and every \(\varepsilon>0\),
\[
\sup_{x\in K}\operatorname{dist}\bigl(x,\{u_{\varepsilon_k}=s\}\bigr)\le\varepsilon
\]
for all sufficiently large \(k\). Taking \(s=0\) and using \(\{u_{\varepsilon_k}=0\}\subset\{|u_{\varepsilon_k}|\le1-b\}\), we also obtain
\[
\sup_{x\in K}\operatorname{dist}\bigl(x,\{|u_{\varepsilon_k}|\le1-b\}\bigr)\le\varepsilon.
\]
Combining above two inequalities with \eqref{e0.12}, we conclude that the level sets \(\{u_{\varepsilon_k}=s\}\) and transition layers \(\{|u_{\varepsilon_k}|\le1-b\}\) converge locally in Hausdorff distance to \(\partial E\).

\end{proof}

Let \(u_\varepsilon\in H^1_{\mathrm{loc}}(\mathbb R^4)\) be a local minimizer of \(\mathcal{E}_{\varepsilon,F}\) with \(|u_\varepsilon|\le1\). For \(Y\in C_c^1(\mathbb R^4;\mathbb R^4)\), set
\[
\phi_s(X)=X+sY(X),
\qquad v_s=u_\varepsilon\circ\phi_s^{-1}.
\]
 The first inner variation of \(\mathcal{E}_{\varepsilon,F}\) at \(s=0\) gives
\begin{equation}\label{equation:stress}
 \int_{\mathbb{R}^4} e_{\varepsilon,F}(u_\varepsilon)\operatorname{div} Y
       -\varepsilon A(Du_\varepsilon)\cdot(DY)^TDu_\varepsilon=0,
\end{equation}
where $ e_{\varepsilon,F}$ is defined by \eqref{equation:scaledenergy}.
Define the varifold $V_\varepsilon$ of $u_\varepsilon$ by
\begin{equation}\label{equation:diffuseV}
 V_\varepsilon(\varphi)=\int_{\{Du_\varepsilon\ne0\}}
 \varphi(X,\mathrm{I}-\nu_\varepsilon\otimes\nu_\varepsilon)
 \varepsilon|Du_\varepsilon|F(Du_\varepsilon),\quad\text{for any }\varphi\in C_c^\infty( \mathbb{R}^4\times G(3)),
\end{equation}
where $\nu_\varepsilon=Du_\varepsilon/|Du_\varepsilon|$ and $G(3)$  denotes the Grassmann manifold of unoriented $3$-dimensional hyperplanes in $\mathbb{R}^{4}$.

Using the first variation formula for anisotropic varifolds and the homogeneity of $F$, we obtain
\[\delta_FV_\varepsilon(Y)=\int_{\mathbb{R}^4} \left[\varepsilon F(Du_\varepsilon)^2\operatorname{div} Y-\varepsilon A(Du_\varepsilon)\cdot(DY)^TDu_\varepsilon\right]dX.\]
Substituting  \eqref{equation:stress} into the above formula,  the first variation $\delta_FV_\varepsilon$ can be rewritten as
  \begin{equation}\label{equation:discrepancystress}\delta_FV_\varepsilon(Y)=\int_{\mathbb{R}^4}\left[\varepsilon\Phi(Du_\varepsilon)-\frac{W(u_\varepsilon)}{\varepsilon}\right]\operatorname{div} Y\ dX.\end{equation}

\begin{lemma}\label{lemma:varifold-stationary}
Let \(u_{\varepsilon_k} \) be as in Proposition~\ref{proposition:phasecompact}, and let \(V_{\varepsilon_k}\) be the varifold defined by \eqref{equation:diffuseV} with \(F=F_k\). Then any subsequential varifold limit \(V\) of \(V_{\varepsilon_k}\) is \(F_\infty\)-stationary, i.e.,
\[
\delta_{F_\infty}V(Y)=0
\qquad\text{for all }Y\in C_c^1(\mathbb R^4;\mathbb R^4).
\]
\end{lemma}

\begin{proof}
It follows from the uniform ellipticity and the local energy bounds in Lemma \ref{lemma:growth} that for any bounded open set $U\subset\mathbb{R}^4$,
\[
\|V_{\varepsilon_k}\|(U)=
\int_U
\varepsilon_k|Du_{\varepsilon_k}|
F_k(Du_{\varepsilon_k})dX\le C\int_U \varepsilon_k\Phi_k(Du_{\varepsilon_k})dX\leq C\mathcal{E}_{\varepsilon_k,F_k}(u_{\varepsilon_k};U)\leq C_U.
\]
Thus, after passing to a further subsequence, $V_{\varepsilon_k}\rightharpoonup V$ locally  in the sense of varifolds.
Recall that
\[\Psi(s)=\int_{-1}^s\sqrt{2W(t)}dt,
\quad c_W=\Psi(1).\]
Then, passing to a further subsequence, it follows from  Proposition \ref{proposition:phasecompact}  that there exists a set $E$ such that
\[
\Psi(u_{\varepsilon_k})
\longrightarrow \Psi( \mathbf{1}_{E^c}-\mathbf{1}_{E})=c_W\mathbf{1}_{E^c}
\quad\text{in }L^1_{\mathrm{loc}}(\mathbb{R}^4).
\]
Let $U\subset\subset\mathbb{R}^4$ be a Lipschitz domain with $P_{F_\infty}(E;\partial U)=0$. By the lower semicontinuity of anisotropic variation,
\(F_k\to F_\infty\) uniformly on \(\mathbb S^3\),
and  Young's inequality \eqref{equation:YoungBV}, we have
\[
\begin{aligned}
c_WP_{F_\infty}(E;U)
\le \liminf_{k\to\infty}
\int_U F_k(D\Psi(u_{\varepsilon_k}))\,dX
&\le \limsup_{k\to\infty}
\int_U F_k(D\Psi(u_{\varepsilon_k}))\,dX\\
&\le \lim_{k\to\infty}
\mathcal{E}_{\varepsilon_k,F_k}(u_{\varepsilon_k};U)
=c_WP_{F_\infty}(E;U).
\end{aligned}
\]
Consequently,
\[
\int_U\left(\sqrt{\varepsilon_k\Phi_k(Du_{\varepsilon_k})}-\sqrt{\varepsilon_k^{-1}W(u_{\varepsilon_k})}\right)^2\,dx=
\mathcal{E}_{\varepsilon_k,F_k}(u_{\varepsilon_k};U)
-\int_U F_k(D\Psi(u_{\varepsilon_k}))\,dx\to0
\]
as $k\to\infty$.
By Cauchy-Schwarz inequality,
\begin{equation}\label{equation:discrepancyvanishing}
\begin{aligned}
&\int_U\left|\varepsilon_k\Phi_k(Du_{\varepsilon_k})-\varepsilon_k^{-1}W(u_{\varepsilon_k})\right|\,dX\\
\le&
\left(\int_U\left(\sqrt{\varepsilon_k\Phi_k(Du_{\varepsilon_k})}+\sqrt{\varepsilon_k^{-1}W(u_{\varepsilon_k})}\right)^2\,dX\right)^{1/2}
\left(\int_U\left(\sqrt{\varepsilon_k\Phi_k(Du_{\varepsilon_k})}-\sqrt{\varepsilon_k^{-1}W(u_{\varepsilon_k})}\right)^2\,dX\right)^{1/2}\\
\le&
\left(
2\mathcal{E}_{\varepsilon_k,F_k}(u_{\varepsilon_k};U)
\right)^{1/2}
\left(\int_U\left(\sqrt{\varepsilon_k\Phi_k(Du_{\varepsilon_k})}-\sqrt{\varepsilon_k^{-1}W(u_{\varepsilon_k})}\right)^2\,dX\right)^{1/2}\to0\quad\text{as }k\to\infty.\\
\end{aligned}
\end{equation}

Fix \(Y\in C_c^1(\mathbb R^4;\mathbb R^4)\),  and choose a bounded open set $U\supset \operatorname{spt} Y$.  The first variation
formula and \eqref{equation:discrepancystress} give
\[
\begin{aligned}
|\delta_{F_\infty}V_{\varepsilon_k}(Y)|\leq& |\delta_{F_k}V_{\varepsilon_k}(Y)|+|\delta_{F_\infty}V_{\varepsilon_k}(Y)-\delta_{F_k}V_{\varepsilon_k}(Y)|\\
\le&
|\delta_{F_k}V_{\varepsilon_k}(Y)|
+
C\|F_k-F_\infty\|_{C^1(\mathbb S^3)}
\|DY\|_{L^\infty}\|V_{\varepsilon_k}\|(U)\\
\le&
\|\operatorname{div}Y\|_{L^\infty}
\int_U\left|\varepsilon_k\Phi_k(Du_{\varepsilon_k})-\varepsilon_k^{-1}W(u_{\varepsilon_k})\right|dX\\
&+
C\|F_k-F_\infty\|_{C^1(\mathbb S^3)}
\|DY\|_{L^\infty}\|V_{\varepsilon_k}\|(U).
\end{aligned}
\]
Passing to the  limit as $k\to\infty$ gives
\[
\delta_{F_\infty}V(Y)
=
\lim_{k\to\infty}
\delta_{F_\infty}V_{\varepsilon_k}(Y)
=0.
\]
Hence the limit  varifold $V$ is
\(F_\infty\)-stationary.
\end{proof}

\section{Minimizing Subgraphs}\label{section:sharpgraph}

Let $E\subset\Omega\times\mathbb{R}$ be a locally minimizer of $P_F$ in $\mathbb{R}^4$. for the perimeter-minimizing set \(E\), the interface \(\partial E\) splits into a regular set and a singular set:
\[
\partial E = \operatorname{Reg}(\partial E) \cup \operatorname{Sing}(\partial E).
\]
 The regular set \(\Sigma:=\operatorname{Reg}(\partial E)\) is an open subset of \(\partial E\), consisting of points where \(\partial E\) is locally a \(C^1\) hypersurface.  The singular points \(S:=\operatorname{Sing}(\partial E)\) are those at which the regular set has density zero.

 The regularity and compactness theory established in Schoen-Simon-Almgren  \cite{SSA} and  Schoen-Simon in \cite{SS} claim that for an $n$-dimensional  minimizing hypersurface in $\mathbb{R}^{n+1}$,  its  interior singular set $S$ satisfies $\mathcal{H}^{n-2}(S)=0$. On the regular set $\Sigma$, let $\mathcal{S}$ be the Euclidean shape operator and choose the normal of a subgraph to point upward.






\begin{lemma}\label{lemma:sharpdensity}
Let $\Omega\subset\mathbb{R}^3$ be open and let $E\subset\Omega\times\mathbb{R}$ be a local minimizer of $P_F$ in $\mathbb{R}^4$.  If $X$ is   a density-one point of $E$, i.e. $\lim_{r\to0}\frac{|E\cap B_r(X)|}{|B_r(X)|}=1$, then
\begin{equation}\label{equation:sharpbounds}
 P_F(E;B_r(X))\le Cr^3\qquad\text{and}
 \qquad |E\cap B_r(X)|\ge cr^4.
\end{equation}
Moreover, if $X$ is   a density-one point of $E^c$, then $|E^c\cap B_r(X)|\ge cr^4.$
\end{lemma}


\begin{proof}
Let $M(r)=|E\cap B_r(X)|$ and set $E':=E\setminus B_r(X)$. Then, the minimality of $E$ implies that $P_F(E;B_R(X))\le  P_F(E';B_R(X))$, where $B_R(X)\supset B_r(X)$ is a ball in $R^4$. Since $E'=E$ on $B_R(X) \setminus \overline{B_r(X)}$,
  we obtain
\[
 P_F(E;B_r(X))\le \beta\mathcal{H}^3(E\cap\partial B_r(X))=\beta M'(r).
\]
 Moreover, the Euclidean perimeter of $E\cap B_r(X)$
  \[P(E\cap B_r(X))=\mathcal{H}^3(E\cap\partial B_r(X))+\mathcal{H}^3(\partial^*E\cap B_r(X))\leq M'(r)+\alpha^{-1}P_F(E;B_r(X)).\]
 The classical isoperimetric inequality therefore gives
\[ M(r)^{3/4}\le C_0 M'(r).\]
At a density-one point $M(r)>0$ for every small positive $r$. Then, we derive that
\[\frac{d}{dr}M(r)^{1/4}=\frac{1}{4}M(r)^{-3/4}M'(r)\geq \frac{1}{4C_0}.\]
  Since $M(0)=0$, integrating above inequality from $0$ to $r$, we get $M(r)\geq (4C_0)^{-4}r^4:=cr^4$. Moreover, for almost every $r>0$, by \eqref{e0.14} and  the co-area formula that
  \[P_F(E;B_r(X))\leq \beta M'(r)=\beta \mathcal{H}^3(E\cap\partial B_r(X))\le\beta\mathcal{H}^3(\partial B_r(X))=2\beta w_4r^3:=Cr^3.\]
 Since \(P_F(E; B_r(X))\) is monotone nondecreasing with respect to $r$, for any \(r > 0\) we may choose a decreasing sequence \(r_j \to r\) such that the estimate holds for every \(r_j\). By monotonicity,
\[
P_F(E; B_r(X)) \le P_F(E; B_{r_j}(X)) \le C r_j^3.
\]
Letting \(j \to \infty\) yields
\[
P_F(E; B_r(X)) \le C r^3\qquad\text{for every}\quad r>0. \]

Finally, since \(F\) is even, \(E^c\) is also minimizing, and the same argument yields
\(|E^c\cap B_r(X)|\ge cr^4\) if \(X\) is a density-one point of \(E^c\).
Approximating boundary points by interior density-one points in each phase gives
\[
|E\cap B_r(X)|\ge cr^4,\qquad |E^c\cap B_r(X)|\ge cr^4
\]
for every \(X\in\partial E\). Moreover, by norm equivalence and the perimeter estimate,
\[
\mathcal H^3(\partial E\cap B_r(X))
\le \alpha^{-1}P_F(E;B_r(X))
\le Cr^3.
\]
Since the singular set \(S\) satisfies \(\mathcal H^3(S)=0\) by the regularity theorem \cite{SSA,SS}, the above estimate holds for the full boundary \(\partial E\).

\end{proof}

\begin{lemma}\label{lemma:heightBV}
Let $\Omega\subset\mathbb{R}^3$ be open and let $E\subset\Omega\times\mathbb{R}$ be a local minimizer of $P_F$ in $\mathbb{R}^4$.  Suppose that $ E=\{(y,t):t<f(y)\}$ is a subgraph in $\Omega\times\mathbb{R}$, where $|f(y)|<\infty$ for almost everywhere $y\in\Omega.$
Then, $f\in L^\infty_{\mathrm{loc}}(\Omega)\cap BV_{\mathrm{loc}}(\Omega)$ and for every $K\subset\subset\Omega$, there is $M_K$ such that
\begin{equation}\label{equation:boundaryheight}
 \partial E\cap(K\times\mathbb{R})\subset K\times[-M_K,M_K].
\end{equation}
\end{lemma}
\begin{proof}
Since $f$ is finite almost everywhere, dominated convergence gives
\begin{equation}\label{equation:emptytranslate}
 E-se_4\longrightarrow\varnothing
 \quad\text{in }L^1_{\mathrm{loc}}(\Omega\times\mathbb{R})
 \quad\text{as }s\to+\infty.
\end{equation}
Assume by contradiction  that \(\operatorname{ess\,sup}_{K_1} f = +\infty\) for some compact \(K_1 \Subset \Omega\). Then there exist \(s_j \to \infty\) such that
\[
|K_1 \cap \{f > s_j + 1\}| > 0.
\]
And hence, the set \(E_j := E \cap (K_1 \times (s_j, s_j+1))\) has positive measure for each \(j\). Choose a density-one point \(X_j = (y_j, t_j) \in E_j\), then the translated point \(X_j - s_j e_4 \in K_1 \times (0,1)\).

Fix \(r > 0\) with \(r < \operatorname{dist}(K_1, \partial \Omega)\). Since \(X_j - s_j e_4\) is a density-one point of \(E - s_j e_4\), Lemma \ref{lemma:sharpdensity} yields
\[
|(E - s_j e_4) \cap B_r(X_j - s_j e_4)| \ge c r^4,
\]
which is a contradiction with \eqref{equation:emptytranslate}. Consequently, \(f\) is locally essentially   bounded above on \(\Omega\), and the same argument applied to   \(E^c\) provides the essential lower bound.  Therefore \(f \in L^\infty_{\mathrm{loc}}(\Omega)\).


For almost every $y\in K$, there exists $M>0$ such that \(|f(y)|\le M\). If \(X=(y,t)\in K\times\mathbb R\) with \(t>M_K\geq M+1\), then \(B_r(X)\cap E=\varnothing\) for some \(r>0\), so \(X\notin\partial E\) by  Lemma \ref{lemma:sharpdensity}.   The same argument applies to \(t<-M_K\leq-M-1\), we obtain
\[
\partial E\cap(K\times\mathbb R)\subset K\times[-M_K,M_K],
\]
which proves \eqref{equation:boundaryheight}.

Choose \(K'\subset\subset K\) and \(T>M_K\). For almost every \(y\in K'\),
\[
f(y)+T=\int_{-T}^{T}\mathbf 1_E(y,t)\,dt.
\]
Differentiating in distributions and using the finite perimeter of \(E\) in \(K'\times[-T,T]\), we deduce that \(D_y f\) is a finite measure. Therefore \(f\in BV_{\mathrm{loc}}(\Omega)\).

\end{proof}


\begin{lemma}\label{lemma:verticalpositive}
Let $\Omega\subset\mathbb{R}^3$ be open and let $ E=\{(y,t):t<f(y)\}$ be a subgraph in $\Omega\times\mathbb{R}$ with $|f(y)|<\infty$   a.e. in $\Omega.$  Assume in addition that   $E $ be a local minimizer of $P_F$ in $\mathbb{R}^4$. Let \(\nu\) be the unit outward normal to \(\partial E\)  and set \(q:=\nu\cdot e_4\).
 Then, $q>0$ on the regular set $\Sigma:=\operatorname{Reg}\partial E$, and $f\in W^{1,1}_{\mathrm{loc}}(\Omega)$.
\end{lemma}


\begin{proof}
Since $E$ is a subgraph, it is obvious that $q\ge0$.  We first notice that by \cite{DuYang2024},   $q$ satisfies the following anisotropic Jacobi equation on $\Sigma$:
\begin{equation}\label{equation:JacobiSharp}
 \operatorname{div}_\Sigma(D^2F(\nu)\nabla_\Sigma q)+\operatorname{tr}_{\Sigma} (D^2F(\nu)\mathcal{S}^2 )q=0,
\end{equation}
where $\nu$ is a unit normal vector field along $\Sigma$ and $\mathcal{S}$ is the Euclidean shape operator on $\Sigma$.

  Since \(\lambda I \le D^2F(\nu)|_{T\Sigma} \le \Lambda I\) for any $\xi\in T\Sigma$ , it follows from the strong maximum principle that, if $q$ vanished at a regular point, then $q\equiv0$ on this regular connected component of $\Sigma$.  Consequently,  \(e_4\) is tangent to \(\Sigma\) everywhere, and this component  is locally a vertical product $\Gamma\times\mathbb{R}$ with \(\mathcal{H}^2(\Gamma)>0\).


By the  regularity theorem established in \cite{SSA,SS}, $\mathcal{H}^1(S)=0$. Thus $\mathcal{H}^1(\pi S)=0$, where $\pi(y,t)=y$. Choose $y\in\Gamma\setminus\pi S$.   Since \(q=0\), the vertical line \(\{y\}\times\mathbb{R}\) lies entirely in  \(\Sigma\). Suppose it terminates at a finite height \(t_0\), i.e., \((y,t_0)\in\partial E\). If \((y,t_0)\in\Sigma\), since \(\Sigma\) is locally a \(C^1\) graph and \(e_4\) is tangent, the curve can be extended locally, contradicting   the definition of $t_0$. If  \((y,t_0)\in S\), then \(y\in\pi S\), contradicting   the choice of \(y\).  Therefore,   the complete line $\{y\}\times\mathbb{R}$ belongs to $\partial E$,  which is a contradiction with \eqref{equation:boundaryheight}. This proves $q>0$.

Set $S'=\pi S$. By \eqref{equation:boundaryheight}, \(S'\cap K\) is compact for every \(K\subset\subset\Omega\), and \(\mathcal H^1(S')=0\). Since \(q>0\) on \(\Omega\setminus S'\), \(\pi|_{\partial E}\) is a local diffeomorphism, hence \(\partial E\) is locally a \(C^1\) graph over \(\Omega\setminus S'\). For each \(y\in\Omega\setminus S'\), there exists a unique $t$ such that $(y,t)\in \partial E$. Indeed, if there were \(t_1<t_2\) with \((y,t_1),(y,t_2)\in\partial E\), then by the subgraph structure \(E=\{(y,t):t<f(y)\}\) we would have \(t_1=f(y)=t_2\), a contradiction. Therefore, \(\partial E\) is a \(C^1\) graph over \(\Omega\setminus S'\). Consequently, \(f\in C^1(\Omega\setminus S')\), and the singular part \(D^s f\) is supported on \(S'\).

  Since \(\mathcal{H}^1(S')=0\) and \(\mathcal{H}^3(S)=0\), we have
$
\mathcal{H}^3(\partial E\cap(S'\times\mathbb{R}))=0.
$
Hence \(P_F(E;S'\times\mathbb{R})=0\). By the  area formula,
\[
\int_{S'}F\left(\frac{dD^{s}f}{d|D^{s}f|},0\right)d|D^{s}f|
\le P_F(E;S'\times\mathbb{R})=0.
\]
Since \(F((\zeta,0))\ge \alpha|\zeta|>0\) for $\zeta\in\mathbb{R}^{3}\setminus\{0\}$, it follows that \(D^{s}f=0\). Therefore \(f\in W^{1,1}_{\mathrm{loc}}(\Omega)\).

\end{proof}

\begin{lemma}\label{lemma:capacity}
Let \(\Omega\subset\mathbb{R}^3\) be open and let \(E\subset\Omega\times\mathbb{R}\) be a local minimizer of \(P_F\) in \(\mathbb{R}^4\).
For every compact \(K\subset\subset\partial E\), there exist Lipschitz functions \(\zeta_j:\partial E\to[0,1]\) such that
\begin{equation}\label{equation:capacitysmall}
\zeta_j=0 \text{ near } S,\qquad
\zeta_j\to 1 \text{ on } \Sigma,\qquad
\int_\Sigma |\nabla_\Sigma \zeta_j|^2 \to 0.
\end{equation}
Moreover, if \(v\in W^{1,2}_{\mathrm{loc}}(\Sigma)\), then for every \(\eta\in C_c^\infty(\partial E)\) there exists a sequence \(\{\varphi_k\}\subset C_c^\infty(\Sigma)\) such that
\[
\|\varphi_k-\eta v\|_{W^{1,2}(\Sigma)}\to0.
\]
If \(v\ge0\), the approximants can be chosen nonnegative.
\end{lemma}


\begin{proof}
Since \(\mathcal H^1(S)=0\), by compactness, for each \(j\) we cover \(S\) by finite balls \(B_{r_{ji}}(X_{ji})\) with $\sum_i r_{ji}<\frac1j,$ and thus $ \max_i r_{ji}<\frac1j$.
For each \(i\), choose a  cutoff function \(\chi_{ji}\) such that
\[
\chi_{ji}=1 \text{ on } B_{r_{ji}}(X_{ji}),\qquad
\chi_{ji}=0 \text{ outside } B_{2r_{ji}}(X_{ji}),\qquad
|D\chi_{ji}|\le \frac{C}{r_{ji}}.
\]
Since \(\mathcal H^3(\Sigma\cap B_R)\le CR^3\), we get
\[
\int_\Sigma |\nabla_\Sigma \chi_{ji}|^2
\le \frac{C}{r_{ji}^2}\mathcal H^3(\Sigma\cap B_{2r_{ji}})
\le C r_{ji}.
\]
Set \(\zeta_j:=1-\max_i \chi_{ji}\). Then \(\zeta_j=0\) near \(S\), \(\zeta_j\to1\) on \(\Sigma\), and summing over \(i\) yields
\[
\int_\Sigma |\nabla_\Sigma \zeta_j|^2
\le \sum_i \int_\Sigma |\nabla_\Sigma \chi_{ji}|^2
\le C\sum_i r_{ji}\to0.
\]

Let \(v\in W^{1,2}_{\mathrm{loc}}(\Sigma)\) and fix \(\eta\in C_c^\infty(\partial E)\). For \(L>0\), set \(T_Lv=\max\{-L,\min(v,L)\}\). For fixed \(L\), since \(\zeta_j=0\) near \(S\), \(\zeta_j\to1\) on \(\Sigma\), and \(\int_\Sigma|\nabla_\Sigma\zeta_j|^2\to0\), we have
\[
\nabla_\Sigma(\zeta_j T_Lv)
=\zeta_j\mathbf{1}_{\{|v|<L\}}\nabla_\Sigma v+T_Lv\nabla_\Sigma\zeta_j.
\]
By dominated convergence,  \(\zeta_j T_Lv\to T_Lv\) in \(L^2(\Sigma)\) and \(\zeta_j\mathbf{1}_{\{|v|<L\}}\nabla_\Sigma v\to\mathbf{1}_{\{|v|<L\}}\nabla_\Sigma v\) in \(L^2(\Sigma)\),  while
\[
\|T_Lv\nabla_\Sigma\zeta_j\|_{L^2(\Sigma)}\le L\|\nabla_\Sigma\zeta_j\|_{L^2(\Sigma)}\to0.
\]
  Hence \(\zeta_j T_Lv\to  T_Lv\) in  \( W^{1,2}(\Sigma)\). On the other hand, since \(T_Lv\to v\) in \(W^{1,2}(\Sigma)\) as \(L\to\infty\), a diagonal argument yields a sequence \(j(L)\) such that \(\zeta_{j(L)}T_Lv\to v\) in \(W^{1,2}(\Sigma)\).  Multiplying by \(\eta\in C_c^\infty(\partial E)\) gives
\[
\eta\zeta_{j(L)}T_Lv\to \eta v \quad\text{in }W_0^{1,2}(\Sigma).
\]
Since \(\zeta_{j(L)}T_Lv=0\) near \(S\), by standard mollification on  \(\Sigma\), there exists a sequence
\(\{\varphi_k\}\subset C_c^\infty(\Sigma)\) such that
\[
\|\varphi_k-\eta v\|_{W^{1,2}(\Sigma)}\to0
\quad\text{as }k\to\infty.
\]
If \(v\ge0\), the approximate sequence  can be chosen so that \(\varphi_k\ge0\).

\end{proof}


\begin{proposition}\label{proposition:logq}
Let \(\Omega \subset \mathbb{R}^3\) be open and let \(E = \{(y,t) : t < f(y)\}\) be a local minimizer of \(P_F\) in \(\Omega \times \mathbb{R}\) with \(f \in L_{\mathrm{loc}}^\infty(\Omega)\).  Set   \(q = \nu \cdot e_4\), and \(z = -\log q\). Then for every \(K \subset\subset \Omega\), there exists \(C_K > 0\) such that
\[
z(X) \le C_K, \qquad \text{for all}\quad X \in \Sigma \cap (K \times \mathbb{R}),
\]
equivalently, $q(X) \ge e^{-C_K} > 0 $ for any $X \in \Sigma \cap (K \times \mathbb{R})$.

\end{proposition}

\begin{proof}
By the preceding positivity result, \(0<q\le1\) on \(\Sigma\),
so \(z\ge0\).
The local height and perimeter bounds imply that, for every
compact \(K'\subset\subset\Omega\), there exists \(M_{K'}\)
such that
\[
\partial E\cap(K'\times\mathbb R)
\subset K'\times[-M_{K'},M_{K'}],
\qquad
\mathcal H^3\bigl(\Sigma\cap(K'\times\mathbb R)\bigr)<\infty.
\]
All integrals below are taken on \(\Sigma\) with respect to
\(d\mathcal H^3\). A base cutoff \(\eta\) is regarded as the
function \(\eta(y,t)=\eta(y)\) restricted to \(\Sigma\).
The singular cutoffs are chosen on a fixed compact ambient
neighborhood of the supports under consideration.

We first prove a local gradient estimate for \(z=-\log q\). Since \(\nabla_\Sigma q=-q\nabla_\Sigma z\), the Jacobi equation \eqref{equation:JacobiSharp} becomes
\begin{equation}\label{equation:logJacobi}
\operatorname{div}_\Sigma(D^2F(\nu)\nabla_\Sigma z)
=
\operatorname{tr}_\Sigma(D^2F(\nu)\mathcal{S}^2)
+\langle D^2F(\nu)\nabla_\Sigma z,\nabla_\Sigma z\rangle.
\end{equation}
Fix \(0\le\eta\le1\), with \(\eta\in C_c^\infty(\Omega)\),
and let \(\zeta_j\) be the cutoff functions from
Lemma~\ref{lemma:capacity}.
Multiplying the above equality by \(\eta^2\zeta_j^2\)
and integrating by parts gives
\[
\int \eta^2\zeta_j^2\left(
\operatorname{tr}_\Sigma(D^2F(\nu)\mathcal{S}^2)
+\langle D^2F(\nu)\nabla_\Sigma z,\nabla_\Sigma z\rangle
\right)
=
-2\int \eta\zeta_j
\langle D^2F(\nu)\nabla_\Sigma z,\nabla_\Sigma(\eta\zeta_j)\rangle.
\]
By the uniform ellipticity of \(D^2F(\nu)\) and Young's inequality, we obtain
\[
\int \eta^2\zeta_j^2\left(
\operatorname{tr}_\Sigma(D^2F(\nu)\mathcal{S}^2)
+\frac12\langle D^2F(\nu)\nabla_\Sigma z,\nabla_\Sigma z\rangle
\right)
\le
C\int|\nabla_\Sigma(\eta\zeta_j)|^2.
\]
Letting \(j\to\infty\), using \(\zeta_j\to1\) on \(\Sigma\), \(\int_\Sigma|\nabla_\Sigma\zeta_j|^2\to0\), and Fatou's lemma, we obtain
\[
\int \eta^2\left(
\operatorname{tr}_\Sigma(D^2F(\nu)\mathcal{S}^2)
+\frac12\langle D^2F(\nu)\nabla_\Sigma z,\nabla_\Sigma z\rangle
\right)\leq C\int|\nabla_\Sigma\eta|^2<\infty.
\]
In particular, $\int \eta^2|\nabla_\Sigma z|^2\leq\alpha^{-1}\int  \eta^2 \langle D^2F(\nu)\nabla_\Sigma z,\nabla_\Sigma z\rangle<\infty.$ Therefore,
\begin{equation}\label{equation:gradlog}
 \int\eta^2(\operatorname{tr}_\Sigma(D^2F(\nu)\mathcal{S}^2)+|\nabla_\Sigma z|^2)<\infty.
\end{equation}
Choosing \(\eta\equiv1\) near any given compact
\(K'\subset\subset\Omega\), we obtain
\[
\int_{\Sigma\cap(K'\times\mathbb R)}
|\nabla_\Sigma z|^2<\infty.
\]

Define
\[
Z:=F(\nu)e_4-qDF(\nu).
\]
By Euler's identity, $Z\cdot\nu=F(\nu)(e_4\cdot\nu)-qDF(\nu)\cdot\nu=F(\nu)q-qF(\nu)=0$ and $\operatorname{div}_\Sigma Z=$
Let \(\{e_i\}\) be a local orthonormal frame on \(\Sigma\). Let \(\{e_i\}\) be a local orthonormal frame on \(\Sigma\).
With the convention \(\nabla_{e_i}\nu=-\mathcal S e_i\),
we have
\[
\nabla_\Sigma F(\nu)
=
-\mathcal S DF(\nu)^T,
\qquad
\nabla_\Sigma q=-\mathcal S e_4^T.
\]
It follows that
\[
\operatorname{div}_\Sigma(F(\nu)e_4)
=
\sum_i
\bigl\langle
\nabla_{e_i}(F(\nu)e_4),e_i
\bigr\rangle
=
\bigl\langle\nabla_\Sigma F(\nu),e_4^T\bigr\rangle
=
-\bigl\langle
\mathcal S DF(\nu)^T,e_4^T
\bigr\rangle.
\]
On the other hand,
\[
\operatorname{div}_\Sigma(qDF(\nu))
=\langle\nabla_\Sigma q,DF(\nu)^T\rangle+q\operatorname{div}_\Sigma DF(\nu)
=-\langle\mathcal{S}e_4^T,DF(\nu)^T\rangle+q\operatorname{div}_\Sigma DF(\nu).
\]
Therefore,
\[
\operatorname{div}_\Sigma Z
=-\langle\mathcal{S}DF(\nu)^T,e_4^T\rangle
+\langle\mathcal{S}e_4^T,DF(\nu)^T\rangle
-q\operatorname{div}_\Sigma DF(\nu)=-q\operatorname{div}_\Sigma DF(\nu)=0,
\]
where the last equality follows from the anisotropic minimal
surface equation.

Now we establish a local \(L^1\) estimate for \(z\).
Choose \(M\) so that \(|t|\le M\) and set
\[
h(t)=t+M+1,
\qquad
z_L=\min(z,L).
\]
For fixed \(L\), the tangent vector field
$
X=\eta h z_L\zeta_j Z
$
has compact support away from the singular set.
The divergence theorem and
\(\operatorname{div}_\Sigma Z=0\) give
\[
\begin{split}
0
&=
\int\operatorname{div}_\Sigma X
=
\int
\bigl\langle
Z,\nabla_\Sigma(\eta h z_L\zeta_j)
\bigr\rangle
\\
&=
\int\zeta_j
\bigl\langle
Z,\nabla_\Sigma(\eta h z_L)
\bigr\rangle
+
\int\eta h z_L
\bigl\langle Z,\nabla_\Sigma\zeta_j\bigr\rangle.
\end{split}
\]
By the Cauchy-Schwarz inequality,
\[\left| \int_\Sigma \eta h z_L \langle Z, \nabla_\Sigma \zeta_j \rangle \right|
\le CL \left( \int_\Sigma |\nabla_\Sigma \zeta_j|^2 \right)^{1/2}.\]
Letting \(j\to\infty\) and using Lemma \ref{lemma:capacity}, we obtain $0 = \int_\Sigma \langle Z, \nabla_\Sigma(\eta h z_L) \rangle$, which implies that
\begin{equation}\begin{split} \label{equation:logL1}
\int\eta F(\nu)z_L
&=\int\eta q F_{\xi_4}(\nu)z_L
-\int h z_LZ\cdot\nabla\eta
-\int h\eta Z\cdot\nabla z_L,\end{split}
\end{equation}
where $\nabla\eta=(D\eta,0)\in\mathbb{R}^4$.
Since \(Z\cdot\nabla_\Sigma\eta=Z\cdot\nabla\eta=-qD_{\xi'}F(\nu)\cdot D\eta\), and   \(qz_L\le ze^{-z}\le e^{-1}\), we obtain
\[\int|\eta qF_{\xi_4}(\nu)z_L|\leq C\int\eta\leq C\quad\text{and}\quad  \int| h z_LZ\cdot\nabla\eta|\leq C\int|D\eta|\leq C,\]
with constants independent of \(L\).
 Furthermore,
\[
\left|\int h\eta Z\cdot\nabla z_L\right|
\le (2M+1)\left(\int\eta^2|\nabla z|^2\right)^{1/2}
\left(\int_{\operatorname{spt}\eta}|Z|^2\right)^{1/2}\leq C,
\]
where we used  \eqref{equation:gradlog} in the last inequality. Hence \(\int\eta F(\nu)z_L\le C\) uniformly in \(L\), which implies that $z\in L^1_{\mathrm{loc}}(\Sigma).$

We next establish a stable Sobolev inequality on $\Sigma$.
For $\varphi\in C_c^1(\Sigma)$, the Michael-Simon inequality~\cite{MS} gives in dimension three
\[
\|\varphi\|_{L^6(\Sigma)}^2
\le
C\int
\left(
|\nabla_\Sigma\varphi|^2+H^2\varphi^2
\right),
\]
where $H$ is the classical mean curvature.
Since \(E\) is a local minimizer of \(P_F\) in
\(\Omega\times\mathbb R\), the second-variation formula yields
\begin{equation}\label{equation:stability}
 \int_\Sigma \operatorname{tr}_{\Sigma}(D^2F(\nu)\mathcal{S}^2)\varphi^2
\le\int_\Sigma\langle D^2F(\nu)\nabla_\Sigma\varphi,\nabla_\Sigma\varphi\rangle.
\end{equation}
Moreover, since $H^2\leq 3|\mathcal{S}|^2\leq C\operatorname{tr}_{\Sigma}(D^2F(\nu)\mathcal{S}^2)$, we obtain
\begin{equation}\label{equation:stableSobolev}
\left(\int|\varphi|^6\right)^{1/3}
\le C\int|\nabla_\Sigma\varphi|^2.
\end{equation}
By density, the same inequality holds for
\(\varphi\in W^{1,2}_0(\Sigma)\).

Next, we establish the local $L^2$ integrability of $z$.
Apply the stable Sobolev inequality \eqref{equation:stableSobolev} to \(\eta^3 z_L\zeta_j\). Letting \(j\to\infty\) and using \(\zeta_j\to1\), \(\int|\nabla\zeta_j|^2\to0\), we obtain
\[
\left(\int \eta^{18} z_L^6\right)^{1/3}
\le C\int |\nabla(\eta^3 z_L)|^2\leq C\int_{\operatorname{spt}\eta} |\nabla z|^2 + C\int \eta^4 z_L^2.
\]
Since \(0\le\eta\le1\), we have \(\eta^4\le\eta^{18/5}\).  H\"older's inequality therefore yields
\[
\int\eta^4z_L^2
\le
\left(\int\eta^{18}z_L^6\right)^{1/5}
\left(\int_{\operatorname{spt}\eta}z_L\right)^{4/5}.
\]
Then, applying Young's inequality, we get
\[
\left(\int\eta^{18}z_L^6\right)^{1/3}
\le
C\int_{\operatorname{spt}\eta}|\nabla_\Sigma z|^2
+
\frac12
\left(\int\eta^{18}z_L^6\right)^{1/3}
+
C\left(\int_{\operatorname{spt}\eta}z_L\right)^2.
\]
Consequently,
\[
\left(\int\eta^{18}z_L^6\right)^{1/3}
\le
C\int_{\operatorname{spt}\eta}|\nabla_\Sigma z|^2
+
C\left(\int_{\operatorname{spt}\eta}z_L\right)^2.
\]
The right-hand side is bounded independently of \(L\).
Letting \(L\to\infty\) and choosing \(\eta\equiv1\) near
\(K'\), we obtain
\[
z\in L^6\bigl(\Sigma\cap(K'\times\mathbb R)\bigr).
\]
The finite area of this portion then gives
\[
\int_{\Sigma\cap(K'\times\mathbb R)}
\left(z^2+|\nabla_\Sigma z|^2\right)<\infty
\qquad
\text{for every compact }K'\subset\subset\Omega.
\]

Since \(z_L\le z\) and \(z\in L^1_{\mathrm{loc}}(\Sigma)\), the right-hand side is bounded uniformly in \(L\). Letting \(L\to\infty\) and using monotone convergence, we obtain \(z\in L^6_{\mathrm{loc}}(\Sigma)\), hence \(z\in L^2_{\mathrm{loc}}(\Sigma)\).  Together with \(\nabla z\in L^2_{\mathrm{loc}}(\Sigma)\), this gives \(z\in W^{1,2}_{\mathrm{loc}}(\Sigma)\).  By Lemma~\ref{lemma:capacity}, for every \(\eta\in C_c^\infty(\partial E)\),
$
\eta z\in W^{1,2}_0(\Sigma).
$

It follows from \eqref{equation:logJacobi} that
\[
\operatorname{div}_\Sigma
\bigl(D^2F(\nu)\nabla_\Sigma z\bigr)\ge0.
\]
Let \(\varphi\) be the restriction to \(\Sigma\) of a
nonnegative function in \(C_c^\infty(\Omega\times\mathbb R)\).
Testing with \(\zeta_j^2\varphi\), we obtain
\[
\int_\Sigma \zeta_j^2\langle D^2F(\nu)\nabla_\Sigma z,\nabla_\Sigma\varphi\rangle
+2\int_\Sigma \zeta_j\varphi\langle D^2F(\nu)\nabla_\Sigma z,\nabla_\Sigma\zeta_j\rangle\le0.
\]
By the Cauchy-Schwarz inequality,
\[\left|
2\int\zeta_j\varphi
\bigl\langle
D^2F(\nu)\nabla_\Sigma z,\nabla_\Sigma\zeta_j
\bigr\rangle
\right|\le
C\|\varphi\|_{L^\infty}
\left(\int_{\operatorname{spt}\varphi}
|\nabla_\Sigma z|^2\right)^{1/2}
\left(\int|\nabla_\Sigma\zeta_j|^2\right)^{1/2}.
\]
Letting \(j\to\infty\), we obtain
\[
\int
\bigl\langle
D^2F(\nu)\nabla_\Sigma z,\nabla_\Sigma\varphi
\bigr\rangle\le0.
\]
Now let \(0\le\varphi\in W^{1,2}(\Sigma)\), and suppose that
\(\varphi=0\) almost everywhere outside
\(\Sigma\cap(K'\times\mathbb R)\) for some compact
\(K'\subset\subset\Omega\).
For fixed \(m>0\), set \(\varphi_m=\min(\varphi,m)\).
By Lemma~\ref{lemma:capacity},
\[
\begin{aligned}
\|\zeta_j\varphi_m-\varphi_m\|_{W^{1,2}(\Sigma)}
&\le
\|(\zeta_j-1)\varphi_m\|_{L^2(\Sigma)}
+
\|(\zeta_j-1)\nabla_\Sigma\varphi_m\|_{L^2(\Sigma)}
+
m\|\nabla_\Sigma\zeta_j\|_{L^2(\Sigma)}\to 0
\end{aligned}
\]
as $j\to\infty$.
Since \(\varphi_m\to\varphi\) in \(W^{1,2}(\Sigma)\),  a standard diagonal argument followed by mollification on \(\Sigma\) yields nonnegative \(\varphi_k\in C_c^\infty(\Sigma)\) with
\[
\|\varphi_k-\varphi\|_{W^{1,2}(\Sigma)}\to0.
\]
In particular, \(\varphi\in W^{1,2}_0(\Sigma)\).
Choose an open set \(U\) with \(K'\subset U\subset\subset\Omega\).
Multiplying by a fixed nonnegative cutoff function \(\eta\in C_c^\infty(\Omega)\) with \(\eta\equiv1\) near \(K'\) and \(\operatorname{spt}\eta\subset U\), we may assume
\[
\operatorname{spt}\varphi_k\subset\Sigma\cap(U\times\mathbb R)
\qquad\text{for every }k.
\]
Testing \eqref{equation:logJacobi} with \(\varphi_k\) gives
\[
\int_\Sigma \langle D^2F(\nu)\nabla_\Sigma z,\nabla_\Sigma\varphi_k\rangle\le0.
\]
Moreover, by the uniform ellipticity of \(F\), we get
\[\left|
\int
\bigl\langle
D^2F(\nu)\nabla_\Sigma z,
\nabla_\Sigma(\varphi_k-\varphi)
\bigr\rangle
\right|\le
C\|\nabla_\Sigma z\|_{L^2(\Sigma\cap(U\times\mathbb R))}
\|\nabla_\Sigma(\varphi_k-\varphi)\|_{L^2(\Sigma)}
\longrightarrow0.
\]
Therefore,
\begin{equation}\label{e0.10}
\int
\bigl\langle
D^2F(\nu)\nabla_\Sigma z,\nabla_\Sigma\varphi
\bigr\rangle\le0
\end{equation}
for every nonnegative \(\varphi \in W^{1,2}_0(\Sigma)\).

Set \(w=z+1\ge1\) and \(w_L=\min(w,L)\), where \(L>1\).
For \(p\ge2\), define
\[
\beta_L(w)=w w_L^{p-2},
\qquad
h_L(w)=w w_L^{(p-2)/2}.
\]
A direct computation gives, almost everywhere,
\[
\beta_L'\ge w_L^{p-2},
\qquad
\frac{\beta_L^2}{\beta_L'}\le h_L^2,
\qquad
|h_L'|^2\le Cp^2\beta_L'.
\]
Fix the compact set \(K\subset\subset\Omega\).
For \(s>0\), set
\[
\Sigma_s:=\Sigma\cap(U_s\times\mathbb R),
\]
where $U_s:=\{y\in\mathbb R^3:\operatorname{dist}(y,K)<s\}$.
Choose \(0<r<R\) such that
\(\overline{U_R}\subset\subset\Omega\).
For \(r\le s<t\le R\), choose
\(\eta\in C_c^\infty(U_t)\) satisfying
\[
\eta\equiv1\text{ on }U_s,
\qquad
0\le\eta\le1,
\qquad
|D\eta|\le\frac{C}{t-s}.
\]
Testing \eqref{e0.10} with \(\eta^2\beta_L(w)\) gives
\[
\begin{split}
\int\eta^2\beta_L'
\bigl\langle
D^2F(\nu)\nabla_\Sigma w,\nabla_\Sigma w
\bigr\rangle
\le&
2\left|
\int\eta\beta_L
\bigl\langle
D^2F(\nu)\nabla_\Sigma w,\nabla_\Sigma\eta
\bigr\rangle
\right|
\\
\le&
\frac12
\int\eta^2\beta_L'
\bigl\langle
D^2F(\nu)\nabla_\Sigma w,\nabla_\Sigma w
\bigr\rangle
+
C\int h_L^2|\nabla_\Sigma\eta|^2,
\end{split}
\]
which implies that
\[
\int\eta^2\beta_L'|\nabla_\Sigma w|^2
\le C\int h_L^2|\nabla_\Sigma\eta|^2.
\]
Applying \eqref{equation:stableSobolev} to
\(\eta h_L(w)\), we obtain
\[
\|\eta h_L(w)\|_{L^6(\Sigma)}^2
\le
C\int|\nabla_\Sigma(\eta h_L(w))|^2
\le
Cp^2\int h_L(w)^2|\nabla_\Sigma\eta|^2.
\]
Taking the limit as
\(L\to\infty\), we obtain
\[
\|w\|_{L^{3p}(\Sigma_s)}
\le
\left(\frac{Cp^2}{(t-s)^2}\right)^{1/p}
\|w\|_{L^p(\Sigma_t)}.
\]
Set
$
p_j=2\cdot3^j,$ and
$r_j=r+(R-r)2^{-j}.$
Applying the above
inequality with $p=p_j$, we get
\[
\|w\|_{L^{p_{j+1}}(\Sigma_{r_{j+1}})}
\le
\left(
\frac{Cp_j^2}{(r_j-r_{j+1})^2}
\right)^{1/p_j}
\|w\|_{L^{p_j}(\Sigma_{r_j})}.
\]
Since $r_j-r_{j+1}=(R-r)2^{-j-1}$
and
$
\prod_{j=0}^{\infty}
\left(
\frac{Cp_j^2}{(r_j-r_{j+1})^2}
\right)^{1/p_j}<\infty,
$
we deduce that
\[
\|w\|_{L^{p_N}(\Sigma_{r_N})}
\le
C_{r,R}\|w\|_{L^2(\Sigma_R)}
\]
for every $N$.
Since \(\Sigma_r\subset\Sigma_{r_N}\), letting
\(N\to\infty\) gives
\[
\operatorname*{ess\,sup}_{\Sigma_r}w
\le
C_{r,R}\|w\|_{L^2(\Sigma_R)}<\infty,
\]
which implies that \(w\in L^\infty_{\mathrm{loc}}(\Sigma)\).
By continuity of \(w\) on \(\Sigma\), this yields a pointwise upper
bound for \(z=w-1\) on every compact subset of the regular set.
This proves the proposition.

\end{proof}

\begin{corollary}\label{corollary:sharpgraph}
Let \(E=\{(y,t):t<f(y)\}\) be a minimizing subgraph
as in Lemma~\ref{lemma:heightBV}.
Then \(f\in W^{1,\infty}_{\mathrm{loc}}(\Omega)\) and
\begin{equation}\label{e0.64}
\operatorname{div} DF((-Df,1))=0
\qquad\text{ in }\Omega.
\end{equation}
\end{corollary}

\begin{proof}
By Lemma~\ref{lemma:verticalpositive},
\(f\in W^{1,1}_{\mathrm{loc}}(\Omega)\) and
\[
q(y,f(y))=\frac1{\sqrt{1+|Df(y)|^2}}
\qquad\text{for a.e. }y\in\Omega.
\]
For every ball \(B'\subset\subset\Omega\),
Proposition~\ref{proposition:logq} gives
\(q(y,f(y))\ge c_{B'}>0\) for a.e. \(y\in B'\).
Hence
\[
|Df(y)|^2=q(y,f(y))^{-2}-1
\le c_{B'}^{-2}-1
\qquad\text{for a.e. }y\in B'.
\]
Since \(|Df|\in L^\infty_{\mathrm{loc}}(\Omega)\), we have \(f\in W^{1,\infty}_{\mathrm{loc}}(\Omega)\) and satisfies \eqref{e0.64}.
\end{proof}


We now discuss the classification of entire minimizing subgraphs.
\begin{lemma}\label{lemma:halfspace}
Let \(E \subset \mathbb{R}^4\) be a nontrivial global \(F\)-perimeter minimizer. Suppose there exist a unit vector \(\nu \in \mathbb{S}^3\) and \(a \in \mathbb{R}\) such that the halfspace
\[
H := \{x \in \mathbb{R}^4 : x \cdot \nu < a\}\subset E.
\]
 Then \(E\) is  a halfspace parallel to \(H\).
\end{lemma}

\begin{proof}
Translate the halfspace to \(H=\{x\cdot e<0\}\) and fix \(p\in\partial E\). Let \(R_j\to\infty\) and set
\[
E_j:=\frac{E-p}{R_j}.
\]
Then, up to a subsequence, \(E_j\) converges locally to a minimizing set \(E_\infty\) with $0\in\partial E_\infty$ and $H\subset E_\infty.$
Let \(V\) be the unit-multiplicity integral varifold induced by the reduced boundary  \(\partial^*E_\infty\), defined by
\[
V(\varphi)=\int_{\partial^*E_\infty}\varphi\bigl(x,T_x\partial^*E_\infty\bigr)\,d\mathcal{H}^3(x),
\qquad
\varphi\in C_c^\infty(\mathbb{R}^4\times G(3)).
\]
Then \(V\) is stationary, supported in \(\{x\cdot e\ge0\}\), and touches \(\partial H\) at \(0\). By the one-sided maximum principle of Solomon-White \cite{SW}, we have
\begin{equation}\label{equation:planeinsupport}
\partial H\subset\operatorname{spt}V=\operatorname{spt}(\partial^*E_\infty).
\end{equation}


Since $\mathcal{H}^1(\mathcal S)=0$, $\mathcal{H}^3$-almost everywhere point of $P$ is a regular point of $\partial E_\infty$. At such a point, \(\partial E_\infty\) is locally a \(C^1\) graph containing \(P\), hence it coincides with \(P\) locally, with multiplicity one. Let \(V_{\partial H}\) be the unit-multiplicity integral varifold induced by the plane \(\partial H\), defined by
\[
V_{\partial H}(\varphi)
=
\int_{\partial H}
\varphi(x,T_x (\partial H ))\,d\mathcal{H}^3(x),
\qquad
\varphi\in C_c^\infty(\mathbb{R}^4\times G(3)).
\]
Since \(V\ge V_{\partial H}\) as varifold measures, the difference
\[
V':=V-V_{\partial H}
\]
is a nonnegative integral varifold  supported in \(\{x\cdot e\ge0\}\). By linearity of the first variation,
\[
\delta_FV'=\delta_FV-\delta_FV_{\partial H}=0.
\]


If $0\in\operatorname{spt} V'$,  the one-sided maximum principle gives $P\subset\operatorname{spt} V'$, which is a contradiction with the fact that $V'=0$ near any of the regular points of $P$. Thus $0\notin\operatorname{spt} V'$, and for some $\rho>0$,
\begin{equation}\label{equation:localsplit}
 V=|P|\quad\text{in }B_\rho.
\end{equation}
By the regularity theory established in \cite{SSA,SS}, there exist \(B_\rho'\subset B_\rho\) such that for large \(j\),
\[\partial E_j\cap B_\rho=\{y+u_j(y)e:y\in\partial H\cap B_\rho\},\]
where $u_j\in C^1(\partial H\cap B_\rho)$ satisfies $u_j(0)=0$ and $\|Du_j\|_{L^\infty}\to0.$
For \(x\in\partial E\), set \(z_j:=(x-p)/R_j\). Then \(z_j\in\partial E_j\) and \(z_j\to0\). Hence for large \(j\),
\[
z_j=y_j+u_j(y_j)e
\]
for some \(y_j\in\partial H\). Therefore
\[
|(x-p)\cdot e|
=R_j|z_j\cdot e|
=R_j|u_j(y_j)|
\le \|Du_j\|_{L^\infty}R_j|y_j|
=\|Du_j\|_{L^\infty}|\Pi_P(x-p)|,
\]
where $\Pi_P(x-p)$  is the projection of \(x-p\) onto \(P\).
Letting \(j\to\infty\),  we obtain
$
(x-p)\cdot e=0.
$
Thus \(x\in p+\partial H\). The arbitrariness of \(x\in\partial E\) implies that
$
\partial E\subset p+\partial H.
$
Moreover, since \(H\subset E\), we conclude that  \(\partial E=p+\partial H\).

\end{proof}

\begin{proposition}\label{theorem:sharpBernstein}
Let \(E=\{(y,t):t<f(y)\}\) be a global \(F\)-perimeter minimizer in \(\mathbb{R}^4\), where \(f:\mathbb{R}^3\to\mathbb{R}\cup\{\pm\infty\}\). Then \(E\) is empty, $\mathbb{R}^4$ or a halfspace.
\end{proposition}

\begin{proof}
Set $A_+=\{y:f(y)=+\infty\}$ and $A_-=\{y:f(y)=-\infty\}$. For each fixed \(s>0\), set
\[
E-se_4=\{(y,t):t<f(y)-s\},\qquad
E+se_4=\{(y,t):t<f(y)+s\}.
\]
Passing to the limit \(s\to+\infty\), by dominated convergence, we have
\begin{equation}\label{equation:infinitefibers}
 E-se_4\to A_+\times\mathbb{R},\qquad
 E+se_4\to(\mathbb{R}^3\setminus A_-)\times\mathbb{R} \quad\text{as}\quad s\to+\infty
\end{equation}
locally in $L^1_{loc}(\mathbb{R}^4)$, and the minimality is preserved in the limit.

 Let $\bar F(\xi')=F(\xi',0)$. We first show that \(A_+\) is a global
\(\bar F\)-perimeter minimizer in \(\mathbb{R}^3\).
Let \(\Omega\subset\subset\mathbb{R}^3\) be a bounded smooth open set,
and let \(B\) be a set of locally finite perimeter such that
$
B\triangle A_+\subset\subset\Omega.
$
For \(L>0\), define
\[
C_L
:=
\bigl(A_+\times(\mathbb{R}\setminus(-L,L))\bigr)
\cup
\bigl(B\times(-L,L)\bigr)
\quad\text{and}\quad
U_L:=\Omega\times(-L-1,L+1).
\]
Then
$
C_L\triangle(A_+\times\mathbb{R})
\subset
(A_+\triangle B)\times[-L,L]
\subset\subset U_L,
$
and
\begin{align*}
&P_F(C_L;U_L)-P_F(A_+\times\mathbb{R};U_L)
\\
&\qquad=
2L\bigl(
P_{\bar F}(B;\Omega)-P_{\bar F}(A_+;\Omega)
\bigr)
+2F(e_4)|A_+\triangle B|.
\end{align*}
 By the minimality of \(A_+\times\mathbb{R}\), it follows that
\[
2L\bigl(
P_{\bar F}(A_+;\Omega)-P_{\bar F}(B;\Omega)
\bigr)
\le
2F(e_4)|A_+\triangle B|.
\]
Dividing by \(2L\) and letting \(L\to+\infty\), we obtain
\[
P_{\bar F}(A_+;\Omega)\le P_{\bar F}(B;\Omega).
\]
Thus \(A_+\) is a global \(\bar F\)-perimeter minimizer. Applying the same argument to
\((\mathbb{R}^3\setminus A_-)\times\mathbb{R}\), we see that
\(\mathbb{R}^3\setminus A_-\) is also a global
\(\bar F\)-perimeter minimizer. Since \(\bar F\) is even,
\(A_-\) is a global \(\bar F\)-perimeter minimizer as well. Consequently, by \cite[Corollary~1]{S}, each of \(A_+\) and \(\mathbb{R}^3\setminus A_-\) is empty, full, or a halfspace.


If \(A_+=\mathbb{R}^3\), then \(f=+\infty\) almost everywhere,
and hence \(E=\mathbb{R}^4\).
If \(A_-=\mathbb{R}^3\), then \(f=-\infty\) almost everywhere,
and hence \(E=\varnothing\).
We may therefore assume that \(E\) is nontrivial.

If \(A_+\) is a halfspace, then $A_+\times\mathbb{R}\subset E.$
Since \(A_+\times\mathbb{R}\) is a halfspace in
\(\mathbb{R}^4\), Lemma~\ref{lemma:halfspace} implies that
\(E\) is a halfspace.

If \(A_-\) is a halfspace, then
$
A_-\times\mathbb{R}\subset E^c.
$
By the evenness of \(F\), \(E^c\) is also a nontrivial global
\(F\)-perimeter minimizer. Applying
Lemma~\ref{lemma:halfspace} to \(E^c\), we conclude that
\(E^c\), and hence \(E\), is a halfspace.

If \(f\) is finite almost everywhere in \(\mathbb{R}^3\),
applying Lemma~\ref{lemma:heightBV} and
Corollary~\ref{corollary:sharpgraph} on bounded base domains gives
$
f\in W^{1,\infty}_{\mathrm{loc}}(\mathbb{R}^3)
$
and
\[
\operatorname{div}_y
D_{\xi'}F((-Df,1))=0
\quad\text{weakly in }\mathbb{R}^3.
\]
The standard interior
elliptic regularity yields that
\(f\in C^2(\mathbb{R}^3)\).
We can therefore apply Simon's anisotropic Bernstein theorem
for entire graphs over \(\mathbb{R}^3\)
\cite[Corollary~1]{S} that
$
f(y)=a\cdot y+b
$
for some \(a\in\mathbb{R}^3\) and \(b\in\mathbb{R}\).
Consequently,
\[
E=\{(y,t):t<a\cdot y+b\}
\]
is a halfspace.  Thus \(E\) is empty, \(\mathbb{R}^4\), or a halfspace.

 \end{proof}

\begin{corollary}\label{corollary:planarinit}
After a rotation and translation, assume that
\(\partial_{x_4}u>0\) and \(u(0)=0\).
For every sequence \(R_j\to+\infty\), set
$
\varepsilon_j:=R_j^{-1},$ $ u_j(X):=u(R_jX).$
Then, after passing to a subsequence, there exists a halfspace
\(E\subset\mathbb R^4\), with \(0\in\partial E\), such that
\[
u_j\to\mathbf 1_{E^c}-\mathbf 1_E
\qquad\text{in }L^1_{\mathrm{loc}}(\mathbb R^4).
\]
Moreover, for every Lipschitz domain
\(\Omega\subset\subset\mathbb R^4\) with \(P_F(E;\partial\Omega)=0\),
\[
\mathcal E_{\varepsilon_j,F}(u_j;\Omega)
\longrightarrow c_WP_F(E;\Omega),
\]
and, for every fixed \(s\in(-1,1)\),
$
\{u_j=s\}\longrightarrow\partial E$ locally in Hausdorff distance.
\end{corollary}

\begin{proof}
By Proposition~\ref{proposition:calibration},
Lemma~\ref{lemma:growth} and
Proposition~\ref{proposition:phasecompact},
after passing to a subsequence,
\[
u_j\to\mathbf 1_{E^c}-\mathbf 1_E
\qquad\text{in }L^1_{\mathrm{loc}}(\mathbb R^4),
\]
where \(E\) is a global \(F\)-perimeter minimizer.
Monotonicity and \(L^1_{\mathrm{loc}}\) convergence imply that,
up to a null set,
\[
E=\{(y,t):t<f(y)\},
\qquad
f:\mathbb R^3\to\mathbb R\cup\{\pm\infty\}
\quad\text{measurable}.
\]
Since \(u_j(0)=0\), the level-set convergence in
Proposition~\ref{proposition:phasecompact} gives
\(0\in\partial E\).
Thus \(E\) is nontrivial, and
Proposition~\ref{theorem:sharpBernstein} implies that
\(E\) is a halfspace.
The energy and level-set convergence follow from
Proposition~\ref{proposition:phasecompact}.
\end{proof}

\section{Approximation of Level Sets}\label{section:normalization}
Let \(\nu\in\mathbb S^3\), and choose an orthonormal basis \(\tau_1,\tau_2,\tau_3\) of \(\nu^\perp\). Set \(Y_\nu=(\tau_1,\tau_2,\tau_3)\). Then $\ Y_\nu^TY_\nu=I_3,\ \nu^TY_\nu=0$ and $Y_\nu Y_\nu^T=I-\nu\otimes\nu$
is the orthogonal projection onto \(\nu^\perp\).   Define
\begin{equation}\label{equation:Lnu}
L_\nu(y,t)=Y_\nu y+tDF(\nu),\qquad
J_\nu=|\det L_\nu|,\qquad
Q_s^\nu=L_\nu(Q_s),
\end{equation}
where \(Q_s=B_s'(0)\times(-s,s)\subset\mathbb R^3\times\mathbb R\).

Let \(P_\nu=(Y_\nu,\nu)\). Since \(P_\nu\) is orthogonal, by Euler's identity,
\begin{equation}
J_\nu=|\det L_\nu|
=|\det P_\nu|\,|DF(\nu)\cdot\nu|
=F(\nu)\ge m,
\end{equation}
where \(m:=\min_{\mathbb S^3}F>0\).
For \((y,t)\in\mathbb R^3\times\mathbb R\),
\begin{equation}\label{e0.15}
|L_\nu(y,t)|
\le |y|+|t|\,|DF(\nu)|
\le C|(y,t)|.
\end{equation}
Conversely, let \(z=L_\nu(y,t)\). Since
$
\nu\cdot z=t\,\nu\cdot DF(\nu)=tF(\nu),
$
we have
\begin{equation}\label{e0.16}
t=\frac{\nu\cdot z}{F(\nu)},\qquad
|t|\le \frac{|z|}{m}\le C|z|.
\end{equation}
Moreover, it follows from $Y_\nu^Tz=Y_\nu^TY_\nu y+tY_\nu^TDF(\nu)
=y+tY_\nu^TDF(\nu)$ that
\begin{equation}\label{e0.17}
y=Y_\nu^Tz-tY_\nu^TDF(\nu),
\qquad
|y|\le |z|+C|t|\le C|z|.
\end{equation}
Therefore, from \eqref{e0.15}-\eqref{e0.17},
\[
|(y,t)|\le C|L_\nu(y,t)|.
\]
Hence \(L_\nu\) and \(L_\nu^{-1}\) are uniformly bounded uniformly in \(\nu\). Consequently, there exist constants \(c,C>0\), independent of \(\nu\), such that
\begin{equation}\label{e0.18}
B_{cs}(0)\subset Q_s^\nu\subset B_{Cs}(0).
\end{equation}


Set $\widehat F_\nu(p)=F(L_\nu^{-T}p)$. Since
$L_\nu^T\nu=F(\nu)e_4,$ and
$ L_\nu^{-T}e_4=\nu/F(\nu),$
we have
\begin{equation}\label{equation:normalizedF}
 \widehat F_\nu(e_4)=1,
 \qquad D\widehat F_\nu(e_4)=L_\nu^{-1}DF(\nu)=e_4.
\end{equation}
For any $p\in\mathbb{R}^4$, the convexity of $\widehat F_\nu$ gives
\[
 \widehat F_\nu(p)\geq \widehat F_\nu(e_4)+D\widehat F_\nu(e_4)\cdot(p-e_4)=1+e_4(p-e_4)=p_4.\]
 On the other hand, since $F$ is even, we get $ \widehat F_\nu(p)= \widehat F_\nu(-p)\geq -p_4$. Therefore,
 \begin{equation}\label{equation:verticalsupport}
  \widehat F_\nu(p)\geq |p_4|\qquad\text{for any}\quad p=(p',p_4)\in\mathbb R^3\times\mathbb R.
\end{equation}
Moreover, since \(L_\nu\) and \(L_\nu^{-1}\) are uniformly bounded for \(\nu\in\mathbb S^3\), and \(F\) has \(C^4\) bounds away from the origin, the Hessian
\[
D^2\widehat F_\nu(p)=L_\nu^{-1}D^2F(L_\nu^{-T}p)L_\nu^{-T}
\]
is uniformly bounded in \(\nu\).


For \(r>0\), define \(v(z):=u(rL_\nu z)\), where
\(\varepsilon:=r^{-1}\) and \(z\in Q_4\).
Then
\[
D_zv=rL_\nu^TD_xu,\qquad
D_Tv=rY_\nu^TD_xu,\qquad
dx=r^4J_\nu\,dz.
\]
Define the rescaled tilt excess and energy by
\begin{equation}\label{equation:TMphysical}
\mathcal T_r(\nu)
:=\frac{1}{r^3J_\nu}
\int_{Q^\nu_{4r}}|\Pi_{\nu^\perp}Du|^2\,dx,
\qquad
\mathcal M_r(\nu)
:=\frac{1}{r^3J_\nu}
\int_{Q^\nu_{4r}}[\Phi(Du)+W(u)]\,dx,
\end{equation}
where \(Q^\nu_{4r}:=L_\nu(Q_{4r})=rL_\nu(Q_4)\).
Since \(Y_\nu\) is an orthonormal basis of \(\nu^\perp\),
we have
\begin{equation}\label{e0.19}
\mathcal T_r(\nu)
=\frac{1}{r^3J_\nu}
\int_{Q^\nu_{4r}}|Y_\nu^TDu|^2\,dx
=\int_{Q_4}\varepsilon|D_Tv|^2\,dz.
\end{equation}
Similarly, we obtain
\begin{equation}\label{equation:TMnormalized1}
\mathcal M_r(\nu)
=r\int_{Q_4}
\bigl[r^{-2}\widehat\Phi_\nu(Dv)+W(v)\bigr]\,dz
=\mathcal E_{\varepsilon,\widehat F_\nu}(v;Q_4),
\end{equation}
where \(\widehat\Phi_\nu(p):=\frac12\widehat F_\nu(p)^2\).
Moreover, since the Jacobian \(r^4J_\nu>0\) is constant,
it preserves minimality, i.e. \(v\) is a minimizer of
\(\mathcal E_{\varepsilon,\widehat F_\nu}\).

Let \(p=(p_T,p_4)\in\mathbb R^3\times\mathbb R\), where \(p_T=(p_1,p_2,p_3)\). Under the normalization \eqref{equation:normalizedF} for \(\widehat F_\nu\), we have
\[
\widehat F_\nu(e_4)=1,\qquad D\widehat F_\nu(e_4)=e_4.
\]
Set \(\widehat A_\nu=D\widehat\Phi_\nu\). Since \(\widehat F_\nu\) is even and positively homogeneous of degree one,
\[
\widehat A_{\nu,T}(0,p_4)=0\in\mathbb{R}^3,\quad\text{and}\qquad \widehat\Phi_\nu(0,p_4)=\frac12 p_4^2.
\]
By convexity and the lower Hessian bound \(D^2\widehat\Phi_\nu\ge \lambda I\), integration along \(s\mapsto(sp_T,p_4)\) gives
\[\widehat\Phi_\nu(p_T,p_4)-\widehat\Phi_\nu(0,p_4)=\int_0^1 \widehat A_{\nu,T}(sp_T,p_4)\cdot p_T ds\geq \lambda \int_0^1 s|p_T|^2 ds=\frac{\lambda}{2}|p_T|^2,\]
which implies that
\begin{equation}\label{eq:horizontalPhi}
\widehat\Phi_\nu(p)-\frac12 p_4^2\ge \frac{\lambda}{2}|p_T|^2.
\end{equation}

Let \(v\) be  a minimizer of \(\mathcal E_{\varepsilon,\widehat F_\nu}\). Denote $Dv=(D_Tv,D_4v)\in\mathbb R^3\times\mathbb R,$
where \(D_Tv=(D_1v,D_2v,D_3v)\). Since \(W(v)=\frac12 g(v)^2\), the energy density
\[
e_{\varepsilon,\widehat F_\nu}(v)
=\varepsilon\widehat\Phi_\nu(Dv)+\varepsilon^{-1}W(v)
\]
satisfies
\begin{equation}\label{eq:verticaldefect}
e_{\varepsilon,\widehat F_\nu}(v)-g(v)|D_4v|
 =\varepsilon\left[\widehat\Phi_\nu(Dv)-\frac12(D_4v)^2\right]
  +\frac{(\varepsilon|D_4v|-g(v))^2}{2\varepsilon} \ge 0.
\end{equation}
If the minimizing phase limit is a horizontal plane of multiplicity one, then
\[
\int_{Q_4}\bigl[e_{\varepsilon,\widehat F_\nu}(v)-g(v)|D_4v|\bigr]\,dz\to 0.
\]
Since the integrand in \eqref{eq:verticaldefect} is nonnegative, both terms vanish in \(L^1\). In particular,
\[
\varepsilon\left[\widehat\Phi_\nu(Dv)-\frac12(D_4v)^2\right]\to 0
\quad\text{in }L^1(Q_4).
\]
Combining this with \eqref{eq:horizontalPhi} yields
\[
\int_{Q_4}\varepsilon|D_Tv|^2\,dz\to 0.
\]

The next lemma shows that small tilt excess and near-minimal energy force the limit to be a single horizontal interface, with energy and transition layer concentrated on it.


\begin{lemma}\label{lemma:onesheet}
Let \(F_k\) be normalized by
\eqref{equation:normalizedF} and let \(u_k(y,t)\) be   global
minimizers of \(\mathcal E_k:=\mathcal E_{\varepsilon_k,F_k}\), with $-1<u_k<1$ and $u_k(0)=0$. Let $\varepsilon_k\to0$ and suppose that
\begin{equation}\label{equation:oneinputs}
\delta_k^2:=\int_{Q_4}\varepsilon_k|D_Tu_k|^2dy\to0,
\qquad
\mathcal E_k(u_k;Q_4)
\le(1+\tau)c_W|B'_4|
\quad\text{for } 0<\tau<1.
\end{equation}
Then, after passing to a subsequence, there exists
\(\sigma\in\{-1,1\}\) such that
\[
u_k\longrightarrow \sigma\operatorname{sgn}(t)
\quad\text{in }L^1_{\mathrm{loc}}(Q_4),
\]
and
\[
\sup_{Q_r\cap\{|t|\ge a\}}
\bigl|u_k-\sigma\operatorname{sgn}(t)\bigr|
\longrightarrow0
\qquad\text{for every }0<a<r<4.
\]
For every Lipschitz domain \(U\subset\subset Q_4\) with
\(\mathcal H^3(\partial U\cap\{t=0\})=0\),
\[
\mathcal E_k(u_k;U)
\longrightarrow
c_W\mathcal H^3(U\cap\{t=0\}).
\]
Moreover, for every \(b\in(0,1)\),
\begin{equation}\label{equation:macroflat}
\sup_{\{|u_k|\le1-b\}\cap Q_3}|t|
\longrightarrow0.
\end{equation}
\end{lemma}


\begin{proof}
Recall that
\[
\omega(s)=\sqrt{2W(s)},\qquad
\Psi(s)=\int_{-1}^s\omega(a)\,da,
\qquad \Psi(1)=c_W.
\]
Since \(D_T\Psi(u_k)=\sqrt{2W(u_k)}D_Tu_k\), by Cauchy--Schwarz and \eqref{equation:oneinputs},
\[
\int_{Q_4}|D_T\Psi(u_k)|
\le\left(\int_{Q_4}\frac{2W(u_k)}{\varepsilon_k}\right)^{1/2}
\left(\int_{Q_4}\varepsilon_k|D_Tu_k|^2\right)^{1/2}
\to0.
\]
Hence the limit \(u_\infty\) depends only on \(t\). By Proposition~\ref{proposition:phasecompact}, \(u_\infty\) is a two-phase function, hence a one-dimensional step function. Let \(m\) denote the number of discontinuities of \(u_\infty\) in \((-4,4)\). Since \(F_\infty(e_4)=1\), we get
\[
m\,c_W|B'_4|\leq\mathcal E_\infty(u_\infty;Q_4)\le (1+\tau)c_W|B'_4|,\qquad \tau<1.
\]
Since $\tau<1$, we obtain $m\leq 1$.
On the other hand, \(u_k(0)=0\) forces the unique discontinuity to be located at \(t=0\), so \(m\ge 1\). Hence \(m=1\), and therefore
\[
u_\infty=\sigma\operatorname{sgn}(t)
\quad\text{for some }\sigma\in\{-1,1\}.
\]

Energy convergence follows from Proposition~\ref{proposition:phasecompact}, and hence for every fixed \(a>0\),
\begin{equation}\label{e0.20}
\mathcal E_k(u_k;Q_3\cap\{|t|\ge a\})\to0.
\end{equation}
Suppose there exist \(X_k\in Q_3\cap\{|t|\ge a\}\) with \(|u_k(X_k)|\le1-b\). By the uniform \(C^{1,\beta}\) estimate of Lemma~\ref{lemma:micro}, there is a ball \(B_{c\varepsilon_k}(X_k)\) on which \(|u_k|\le1-b/2\). Since \(W\ge c_b'>0\) on \(B_{c\varepsilon_k}(X_k)\),
\[
\int_{B_{c\varepsilon_k}(X_k)}e_{\varepsilon_k,F_k}(u_k)\,dx
\ge \varepsilon_k^{-1}c_b'|B_{c\varepsilon_k}|
\ge Cc_b'>0,
\]
contradicting \eqref{e0.20}. Hence
\[
\sup_{\{|u_k|\le1-b\}\cap Q_3}|t|\to0.
\]

Therefore, for fixed \(a>0\) and \(k\) large, we have \(|u_k|>1-b\) in \(Q_r\cap\{|t|\ge a\}\). By Lemma~\ref{lemma:modicaclear}, $e_{\varepsilon_k,F_k}(u_k)\to0$ uniformly in $Q_r\cap\{|t|\ge a\}$,
then we have
\[
W(u_k)\le\varepsilon_k e_{\varepsilon_k,F_k}(u_k)\to0
\quad\text{uniformly in }Q_r\cap\{|t|\ge a\}.
\]
Consequently,
\[
u_k\to\sigma\operatorname{sgn}(t)
\quad\text{uniformly in }Q_r\cap\{|t|\ge a\}
\]
for any \(0<a<r<4\).
\end{proof}


\begin{lemma}\label{lemma:masstransfer}
Let \(u\) be a global minimizer of \(\mathcal E_{1,F}\),
with \(-1<u<1\) and \(u(0)=0\).
Fix \(0<\tau<1\), \(0<\theta<1\) and \(C_N<\infty\).
Then there exist \(\varepsilon_0,\delta_0>0\) such that
the following holds.
If for some \(\nu\in\mathbb S^3\) and some
\(r\ge\varepsilon_0^{-1}\),
\[
\mathcal T_r(\nu)\le\delta_0^2,\qquad
\mathcal M_r(\nu)\le(1+\tau)c_W|B'_4|,
\]
then every \(\mu\in\mathbb S^3\) satisfying
\[
|\mu-\nu|\le C_N\sqrt{\mathcal T_r(\nu)}
\]
satisfies
\[
\mathcal M_{\theta r}(\mu)\le(1+\tau)c_W|B'_4|.
\]
\end{lemma}


\begin{proof}

Assume by contradiction that there exist sequences
\(r_k\to\infty\),
\(\delta_k:=\sqrt{\mathcal T_{r_k}(\nu_k)}\to0\),
and \(|\mu_k-\nu_k|\le C_N\delta_k\), such that
\begin{equation}\label{e0.21}
\mathcal M_{r_k}(\nu_k)\le(1+\tau)c_W|B'_4|,
\qquad
\mathcal M_{\theta r_k}(\mu_k)>(1+\tau)c_W|B'_4|.
\end{equation}
Set
\[
\varepsilon_k:=r_k^{-1},\qquad
v_k(z):=u(r_kL_{\nu_k}z),\qquad
\Omega_k:=L_{\nu_k}^{-1}L_{\mu_k}(Q_{4\theta}).
\]
Since \(|\mu_k-\nu_k|\to0\), the maps
\(L_{\nu_k}^{-1}L_{\mu_k}\) converge to an orthogonal
transformation, hence, up to a base rotation,
\[
\Omega_k\to Q_{4\theta}
\qquad\text{in the Hausdorff distance.}
\]
Then, for \(\theta\) fixed sufficiently small,
\(\Omega_k\subset\subset Q_4\) for all large \(k\).

 Recall that
\[
\mathcal M_{\theta r_k}(\mu_k)
=\frac{J_{\nu_k}}{J_{\mu_k}}\,\theta^{-3}
\int_{\Omega_k}
e_{\varepsilon_k,\widehat F_{\nu_k}}(v_k)\,dz.
\]
Since \(\Omega_k\to Q_{4\theta}\), for every \(\eta>0\)
and all large \(k\),
\[
Q_{4\theta-\eta}\subset\Omega_k\subset Q_{4\theta+\eta}.
\]
By Lemma~\ref{lemma:onesheet},
\[
\int_{Q_s}
e_{\varepsilon_k,\widehat F_{\nu_k}}(v_k)\,dz
\longrightarrow c_W|B'_s|
\]
for every fixed \(0<s<4\).
Taking \(s=4\theta\pm\eta\) and then letting \(\eta\to0\),
we obtain
\[
\int_{\Omega_k}
e_{\varepsilon_k,\widehat F_{\nu_k}}(v_k)\,dz
\longrightarrow c_W|B'_{4\theta}|.
\]
Moreover, since \(F\) is continuous and
\(\mu_k-\nu_k\to0\), the Jacobian ratio
$
 J_{\nu_k}/J_{\mu_k} \to1,
$
and we obtain
\[
\mathcal M_{\theta r_k}(\mu_k)
\longrightarrow
\theta^{-3}c_W|B'_{4\theta}|
=c_W|B'_4|,
\]
which is a contradiction with \eqref{e0.21}.

\end{proof}



In the rest of this section, we assume   \eqref{equation:oneinputs}  and, in addition, assume
\begin{equation}\label{equation:epsdelta}
 \varepsilon_k/\delta_k\longrightarrow0.
\end{equation}
Fix $r_0=1/32$. For $X\in Q_{5/2}$ define
\begin{equation}\label{equation:gamma}
 \gamma_k(X)^2=
 \max_{\varepsilon_k\le s\le r_0}s^{-3}
 \int_{B_s(X)}\varepsilon_k|D_Tu_k|^2\,dx.
\end{equation}
For large \(k\) the interval \([\varepsilon_k,r_0]\) is nonempty, and \(\gamma_k\) is continuous in \(X\). A point \(X\) is called \(\eta\)-good if \(\gamma_k(X)\le\eta\), and \(\eta\)-bad if \(\gamma_k(X)\ge\eta\). 

\begin{lemma}\label{lemma:goodpoint}
Fix \(b\in(0,1)\) and \(L,\rho,\ell>0\). There exist
\(\eta,\varepsilon_0,\delta_0>0\) such that if
\(\varepsilon\le\varepsilon_0\), \(\delta\le\delta_0\), and
\(X\in Q_{5/2}\) is \(\eta\)-good with \(u(X)\in I_b:=[-1+b,1-b]\),
then there exists \(\sigma_X\in\{-1,1\}\) such that
\begin{equation}\label{e0.22}
\bigl\|u(X+\varepsilon z)-g(\sigma_X z_4+a_X)\bigr\|_{C^1(B_L)}\le\rho,
\end{equation}
where \(a_X:=g^{-1}(u(X))\). Moreover,
\begin{equation}\label{equation:micronondeg}
\varepsilon|u_t|\ge c_{b,L}>0,\qquad
|D_Tu|\le\ell|u_t|
\quad\text{in }B_{L\varepsilon/2}(X),
\end{equation}
where \(c_{b,L}\) depends only on \(b,L,W\), and
\begin{equation}\label{equation:globalcone}
|Z_4-X_4|\le\ell|Z'-X'|
\qquad
\text{for every }Z\in Q_{5/2}\text{ with }u(Z)=u(X).
\end{equation}
\end{lemma}


\begin{proof}

Assume by contradiction  that  there are sequences $\eta_j\to0,$ $\varepsilon_j\to0,$   $\delta_j\to0,$
and \(\eta_j\)-good points \(X_j\in Q_{5/2}\) with \(u_j(X_j)\in I_b\), such that \eqref{e0.22} or \eqref{equation:micronondeg} fails. Define
\[
U_j(z):=u_j(X_j+\varepsilon_j z).
\]
For every fixed \(R\ge1\) with \(R\varepsilon_j\le r_0\),   the good-point condition yields
\begin{equation}\label{equation:microtilt}
\int_{B_R}|D_TU_j|^2\,dz
=\varepsilon_j^{-3}\int_{B_{R\varepsilon_j}(X_j)}\varepsilon_j|D_Tu_j|^2
\le\eta_j^2R^3\to0.
\end{equation}
The uniform local energy bound \eqref{equation:growth} gives
$
\mathcal E_{1,F_j}(U_j;B_R)\le CR^3.
$
Then, by Lemma~\ref{lemma:micro}, up to a subsequence, $U_j$ converges  to   \(U_\infty\) in \(C^1_{\mathrm{loc}}(\mathbb{R}^4)\).  Passing to the limit in \eqref{equation:microtilt} gives \(U_\infty\) depends only on \(z_4\)  and satisfies
\begin{equation}\label{e0.23}
U_\infty(z_4)''=W'(U_\infty(z_4)).
\end{equation}


By normalization and the \(C^1_{\mathrm{loc}}\) convergence,
for every \(R\ge2\),
\[
cR^3\int_{-R/2}^{R/2}
\left(\frac12|U_\infty'|^2+W(U_\infty)\right)\,dz_4
\le
\liminf_{j\to\infty}\mathcal E_{1,F_j}(U_j;B_R)
\le CR^3.
\]
Letting \(R\to\infty\), we obtain
\[
\int_{\mathbb R}
\left(\frac12|U_\infty'|^2+W(U_\infty)\right)\,dz_4
<\infty.
\]
Hence there exists \(t_m\to+\infty\) such that
\[
|U_\infty'(t_m)|^2+W(U_\infty(t_m))\to0.
\]
Since
$
\left(\frac12|U_\infty'|^2-W(U_\infty)\right)'
=
U_\infty'\bigl(U_\infty''-W'(U_\infty)\bigr)=0,
$
we have
\[
\frac12|U_\infty'|^2-W(U_\infty)
=
\lim_{m\to\infty}
\left(\frac12|U_\infty'(t_m)|^2
-W(U_\infty(t_m))\right)=0.
\]

Since \(U_\infty(0)\in I_b\subset(-1,1)\), we have \(W(U_\infty(0))>0\) and $U_\infty'(0)=2W(U_\infty(0))>0$.  If \(U_\infty'(t_0)=0\) for some \(t_0\), then the same identity yields \(W(U_\infty(t_0))=0\), hence \(U_\infty(t_0)=\pm1\). The uniqueness of the ODE \eqref{e0.23} implies that \(U_\infty\equiv\pm1\), which is a contradiction with \(U_\infty(0)\in(-1,1)\). Therefore \[
-1<U_\infty<1,\qquad
|U_\infty'|=\sqrt{2W(U_\infty)}>0.
\]
Consequently,
\[
U_\infty(z_4)=g(\sigma z_4+a),\]
where $\sigma:=\operatorname{sgn}U_\infty'(0)\in\{-1,1\},$ and $a:=g^{-1}(U_\infty(0))\in g^{-1}(I_b).$
Set
\[
a_j:=g^{-1}(u_j(X_j))
=g^{-1}(U_j(0))\longrightarrow a.
\]
Then
\begin{equation}\label{e0.24}
\begin{aligned}
\|U_j-g(\sigma z_4+a_j)\|_{C^1(B_L)}
&\le
\|U_j-U_\infty\|_{C^1(B_L)}\\
&\quad+
\|g(\sigma z_4+a)-g(\sigma z_4+a_j)\|_{C^1(B_L)}
\longrightarrow0.
\end{aligned}
\end{equation}
Moreover, since \(g^{-1}(I_b)\) is compact and \(g'>0\), we have
\[
\inf_{|z_4|\le L}|U_\infty'(z_4)|\ge2c_{b,L},
\]
where $c_{b,L}:=\frac12\min_{\substack{\xi\in g^{-1}(I_b)\\ |s|\le L}}g'(\xi+s)>0.$
The \(C^1\) convergence and \(D_TU_\infty=0\) imply that,
for all sufficiently large \(j\),
\begin{equation}\label{e0.25}
\varepsilon_j|u_{j,t}|\ge c_{b,L},
\qquad
|D_Tu_j|\le\ell|u_{j,t}|
\quad\text{in }B_{L\varepsilon_j/2}(X_j).
\end{equation}
Thus both \eqref{e0.22} and \eqref{equation:micronondeg}
hold for all sufficiently large \(j\), this is a contradiction with our assumption.


In the following, we prove \eqref{equation:globalcone}. Assume by contradiction that there exist
\(\eta_j\)-good points \(X_j\in Q_{5/2}\), and
\(Z_j\in Q_{5/2}\) such that
\[
u_j(X_j)=u_j(Z_j)\in I_b,\qquad
|Z_{j,4}-X_{j,4}|>\ell|Z_j'-X_j'|.
\]
By \eqref{equation:macroflat}, we have
\[
|X_{j,4}|+|Z_{j,4}|\to0 \qquad\text{and}\qquad
d_j^2:=|Z_j-X_j|^2
\le(1+\ell^{-2})|Z_{j,4}-X_{j,4}|^2\to0.
\]
Moreover,
\begin{equation}\label{equation:verticalfraction}
\frac{|Z_{j,4}-X_{j,4}|}{d_j}
\ge\frac{\ell}{\sqrt{1+\ell^2}}.
\end{equation}

Suppose first that \(d_j/\varepsilon_j\le M<\infty\).
By a similar argument as above, after passing to a subsequence,
\[
U_j(z):=u_j(X_j+\varepsilon_j z)
\longrightarrow g(\sigma z_4+a)
\quad\text{in }C^1(B_{M+1}),
\]
where $\sigma\in\{-1,1\} $ and $a\in g^{-1}(I_b).$
Since \(g'>0\), for all sufficiently large \(j\), $u_{j,t}$ has fixed sign and
\[
|D_Tu_j(P_j(s))|
\le\frac{\ell}{2}|u_{j,t}(P_j(s))|,
\quad 0\le s\le1,
\]
where \(P_j(s):=X_j+s(Z_j-X_j)\). Since \(u_j(Z_j)=u_j(X_j)\), we have
\[0=\int_0^1\frac{d}{ds}u(j)(P_j(s))ds=\int_0^1 Du_j(P_j(s))\cdot (Z_j-X_j)ds.\]
Using the decomposition $Du_j(P_j(s))\cdot (Z_j-X_j)=u_{j,t}(Z_{j,4}-X_{j,4})+D_T u_j\cdot(Z_j'-X_j')$ and \eqref{equation:micronondeg}, we obtain
\[\left|(Z_{j,4}-X_{j,4})\right| \int_0^1 |u_{j,t}|ds=\left|\int_0^1 D_Tu_j\cdot(Z_j'-X_j')ds\right|\leq \frac{\ell}{2}|Z_j'-X_j'|  \int_0^1 |u_{j,t}|ds,\]
which implies that $\left|(Z_{j,4}-X_{j,4})\right| \leq \frac{\ell}{2}|Z_j'-X_j'| .$ This is a contradiction with the assumption \(|Z_{j,4}-X_{j,4}|>\ell|Z_j'-X_j'|\).


In the remaining case \(d_j/\varepsilon_j\to\infty\), set
\[
V_j(z):=u_j(X_j+d_jz),\qquad
\widetilde\varepsilon_j:=\frac{\varepsilon_j}{d_j}\to0.
\]
Then \(V_j\) minimizes \(\mathcal E_{\widetilde\varepsilon_j,F_j}\).
For every fixed \(R>0\), we have
\(\varepsilon_j\le Rd_j\le r_0\) for all sufficiently large \(j\), and
\begin{equation}\label{equation:mesotilt}
\int_{B_R}\widetilde\varepsilon_j|D_TV_j|^2\,dz
=d_j^{-3}\int_{B_{Rd_j}(X_j)}
\varepsilon_j|D_Tu_j|^2\,dx
\le\eta_j^2R^3\to0.
\end{equation}
Moreover,
\[
\mathcal E_{\widetilde\varepsilon_j,F_j}(V_j;B_R)
=d_j^{-3}\mathcal E_{\varepsilon_j,F_j}
(u_j;B_{Rd_j}(X_j))
\le CR^3.
\]
By Proposition~\ref{proposition:phasecompact}, after passing
to a subsequence,
\[
F_j\to F_\infty,\qquad
V_j\to2\chi_E-1
\quad\text{in }L^1_{\mathrm{loc}}(\mathbb R^4),
\]
where \(E\) is a global \(F_\infty\)-perimeter minimizer.

Set \(\Psi(s):=\int_{-1}^s\sqrt{2W(\xi)}\,d\xi\). Then
\[
\Psi(V_j)\to c_W\chi_E
\quad\text{in }L^1_{\mathrm{loc}}(\mathbb R^4),
\]
and by Cauchy-Schwarz inequality,
\[
\int_{B_R}|D_T\Psi(V_j)|\,dz
\le
\left(\int_{B_R}
\frac{2W(V_j)}{\widetilde\varepsilon_j}\,dz\right)^{1/2}
\left(\int_{B_R}
\widetilde\varepsilon_j|D_TV_j|^2\,dz\right)^{1/2}
\le C\eta_jR^3\to0.
\]
Consequently, \(D_T\chi_E=0\), which means that \(E\) is a vertical cylinder. More precisely, there exists a set \(J\subset\mathbb{R}\) with \(\chi_J\in BV_{\mathrm{loc}}(\mathbb{R})\) such that \(E=\mathbb{R}^3\times J\) up to a set of measure zero.

Passing to a further subsequence, let
\[
\zeta_j:=\frac{Z_j-X_j}{d_j}\to\zeta\in\mathbb S^3.
\]
It follows from \eqref{equation:verticalfraction} that
$
|\zeta_4|\ge\frac{\ell}{\sqrt{1+\ell^2}}>0.
$
Since \(V_j(0)=V_j(\zeta_j)\in I_b\),  we get
\[
0,\zeta\in\partial E.
\]
Therefore, $0,\zeta_4\in \partial J$ and since \(\chi_J\in BV_{\mathrm{loc}}(\mathbb{R})\),  there exist \(a_-<a_+\) such that
\[
a_-,a_+\in \partial J,\qquad \partial J\cap(a_-,a_+)=\varnothing,
\]
and
\[
\partial E\cap
\bigl(\mathbb R^3\times[a_-,a_+]\bigr)
=
\bigl(\mathbb R^3\times\{a_-\}\bigr)
\cup
\bigl(\mathbb R^3\times\{a_+\}\bigr).
\]


Set \(d:=a_+-a_->0\) and
\[
C_R:=B_R'\times(a_-,a_+),\qquad
E_R:=E\triangle C_R,\qquad
\Omega_R:=B_{R+1}'\times(a_--1,a_++1).
\]
Then \(E_R\triangle E\subset\subset\Omega_R\).  Since $\partial J\cap(a_-,a_+)=\varnothing,$
by minimality and \(F_\infty(\pm e_4)=1\), we have
\[
0\le P_{F_\infty}(E_R;\Omega_R)
     -P_{F_\infty}(E;\Omega_R)
\le -2|B_R'|
     +Cd\,\mathcal H^2(\partial B_R')
\le -cR^3+CdR^2<0
\]
for all sufficiently large \(R\), a contradiction.
Thus \eqref{equation:globalcone} holds.

\end{proof}

\begin{lemma}\label{lemma:transfer}
If \(X,Z \in \overline{B'_2}\times[-1,1]\) and \(|Z-X|\le L\varepsilon_k\), then
\begin{equation}\label{equation:transfer}
 \gamma_k(Z)^2\le C_L\max\{\gamma_k(X)^2,\delta_k^2\},
\end{equation}
where $C_L=\max\{(1+L)^3,((1+L)/r_0)^3\}.$
\end{lemma}

\begin{proof}
Fix \(s\in[\varepsilon_k,r_0]\). Since \(|Z-X|\le L\varepsilon_k\le Ls\), we have \(B_s(Z)\subset B_{(1+L)s}(X)\).
If \((1+L)s\le r_0\), then
\[
s^{-3}\int_{B_s(Z)}\varepsilon_k|D_Tu_k|^2
\le (1+L)^3\bigl((1+L)s\bigr)^{-3}\int_{B_{(1+L)s}(X)}\varepsilon_k|D_Tu_k|^2
\le (1+L)^3\gamma_k(X)^2.
\]
On the other hand, if \((1+L)s>r_0\), then \(B_s(Z)\subset Q_4\) and \(s^{-3}<((1+L)/r_0)^3\), hence
\[
s^{-3}\int_{B_s(Z)}\varepsilon_k|D_Tu_k|^2
\le \left(\frac{1+L}{r_0}\right)^3\delta_k^2.
\]
Combining the two cases and taking the supremum over \(s\in[\varepsilon_k,r_0]\), the claim holds.
\end{proof}


\begin{lemma}\label{lemma:badballs}
Fix \(b\in(0,1)\) and \(\eta>0\),  set the bad set
\[
\mathcal B_k:=
\{X\in Q_{9/4}:u_k(X)\in I_b,\ \gamma_k(X)>\eta\}.
\]
Then, there exist \(X_{ki}\in\mathcal B_k\) and
\(\rho_{ki}\ge\varepsilon_k\) such that
$
\mathcal B_k\subset\bigcup_i B_{\rho_{ki}}(X_{ki}),
$
and
\begin{equation}\label{equation:badpacking}
\sum_i\rho_{ki}^3\le C\eta^{-2}\delta_k^2,
\qquad
\int_{\bigcup_iB_{2\rho_{ki}}(X_{ki})}e_k
\le C\eta^{-2}\delta_k^2.
\end{equation}
Moreover, let \(\pi(y,t):=y\) be the projection, then
$
|\pi(\mathcal B_k)|\le C\eta^{-2}\delta_k^2.$  
\end{lemma}

\begin{proof}
For each \(X\in\mathcal B_k\), choose
\(s_X\in[\varepsilon_k,r_0]\) such that
\[
\frac12\eta^2s_X^3
\le\int_{B_{s_X}(X)}\varepsilon_k|D_Tu_k|^2.
\]
By the Vitali \(5r\)-covering lemma, there exist pairwise disjoint
balls \(B_{s_{ki}}(X_{ki})\subset Q_4\) satisfying
\[
\mathcal B_k\subset\bigcup_i B_{5s_{ki}}(X_{ki}).
\]
Set \(\rho_{ki}:=5s_{ki}\). Then
\[
\sum_i\rho_{ki}^3
\le C\eta^{-2}\sum_i
\int_{B_{s_{ki}}(X_{ki})}\varepsilon_k|D_Tu_k|^2
\le C\eta^{-2}\delta_k^2.
\]
Applying the local growth bound \eqref{equation:growth}, we obtain
\[
\int_{\bigcup_iB_{2\rho_{ki}}(X_{ki})}e_k
\le\sum_i\int_{B_{2\rho_{ki}}(X_{ki})}e_k
\le C\sum_i\rho_{ki}^3
\le C\eta^{-2}\delta_k^2.
\]
Consequently, the projection of the bad set
satisfy
\[
|\pi(\mathcal B_k)|
\le\sum_i|B'_{\rho_{ki}}|
\le C\eta^{-2}\delta_k^2,
\]
The same construction applies to any subset of
\(\mathcal B_k\).
\end{proof}


Fix \(\ell\in(0,1/4)\) and let \(\eta_0>0\) be the good-point threshold satisfies Lemma~\ref{lemma:goodpoint}. For large \(k\), define
\begin{equation}\label{equation:Dzero}
D_{k,0}=\left\{y\in\overline{B'_2}:\exists t\in[-1,1],\ u_k(y,t)=0,\ \gamma_k(y,t)\le\eta_0\right\}.
\end{equation}
Since \(u_k\) and \(\gamma_k\) are continuous, 
the set \(D_{k,0}\) is compact. For each \(y\in D_{k,0}\),  there exists a unique \(t\in[-1,1]\) such that \(u_k(y,t)=0\). Indeed, if \((y,t_1)\) and \((y,t_2)\) are two good zeros, the global cone bound \eqref{equation:globalcone} gives
\[
|t_1-t_2|\le\ell|y-y|=0,
\]
then $t_1=t_2$. We then define \(f_k(y):=t\). Applying \eqref{equation:globalcone} to \((y,f_k(y))\) and \((z,f_k(z))\) yields
\[
|f_k(y)-f_k(z)|\le\ell|y-z|,
\]
where $y,z\in D_{k,0}$.
If \(y \notin D_{k,0}\), then every zero \(t\) of \(u_k(y,\cdot)\) satisfies \(\gamma_k(y,t)>\eta_0\), i.e., it is a bad zero.
Thus,  Lemma~\ref{lemma:badballs} yields
\[
|B'_2\setminus D_{k,0}|\le C\delta_k^2\to0,
\]
which implies that \(D_{k,0}\neq\emptyset\) for all sufficiently large \(k\).
For \(y\in B'_2\), define
\[
\widehat f_k(y):=
\inf_{z\in D_k}\bigl\{f_k(z)+\ell|y-z|\bigr\},
\qquad
\widetilde f_k(y):=
\max\{-m_k,\min\{m_k,\widehat f_k(y)\}\}.
\]
where $m_k:=\sup_{\{u_k=0\}\cap Q_3}|t|\to0$ by \eqref{equation:macroflat}.
Then
\[
\widetilde f_k=f_k\quad\text{on }D_k\cap B'_2,
\qquad
\operatorname{Lip}(\widetilde f_k;B'_2)\le\ell.
\]
Fix this extension and still denote it by \(f_k\). We have
\begin{equation}\label{equation:fsmall}
\|f_k\|_{L^\infty(B'_2)}\le m_k\longrightarrow0.
\end{equation}
We now use \(f_k\) to define the rescaled heights.


\begin{proposition}\label{proposition:common}
Fix \(b\in(0,1)\). There exist
\(\alpha_b,\beta_b,\varepsilon_b,\delta_b,L_b,c_b,C_b>0\)
such that, if
$
0<\varepsilon_k\le\varepsilon_b,$
$\delta_k\le\delta_b,$
the following hold. Define
\[
D_{k,b}:=
\left\{
y\in\overline{B'_2}:
\exists t\in[-1,1],\
u_k(y,t)=0,\ \gamma_k(y,t)\le\alpha_b
\right\}
\subset D_{k,0}.
\]
For every \((y,s)\in D_{k,b}\times I_b\), there exists a unique $f_k^s(y)$ such that $u_k(y,f_k^s(y))=s$,
and
\begin{equation}\label{equation:levelsnearzero}
|f_k^s(y)-f_k(y)|\le L_b\varepsilon_k,\qquad
\gamma_k(y,f_k^s(y))\le\beta_b,\qquad
\varepsilon_k|u_{k,t}(y,f_k^s(y))|\ge c_b.
\end{equation}
Conversely,
for every \((y,t)\in Q_{3/2}\cap\{u_k\in I_b\}\)
with \(\gamma_k(y,t)\le\beta_b\),
\begin{equation}\label{equation:centraltozero}
y\in D_{k,0},\qquad\text{and}\qquad
|t-f_k(y)|\le L_b\varepsilon_k.
\end{equation}
Moreover,
\begin{equation}\label{equation:badbaseb}
|B'_2\setminus D_{k,b}|\le C_b\delta_k^2,
\end{equation}
and there exist balls centered at \(\beta_b\)-bad central
points such that
\begin{equation}\label{equation:completecover}
\{u_k\in I_b\}\cap Q_{3/2}
\subset
N_{L_b\varepsilon_k}(\Gamma_k)
\cup\bigcup_i B_{\rho_{ki}}(X_{ki}),
\end{equation}
where   $\Gamma_k:=\{(y,f_k(y)): y\in B_1'\}$ is the graph of $f_k$, $\sum_i\rho_{ki}^3\le C_b\delta_k^2$ and $\sup_i |(X_{ki})_4|\longrightarrow0.$
\end{proposition}


\begin{proof}
Since \(g^{-1}(I_{b/2})\) is compact, choose \(L_b>0\) such that
\begin{equation}\label{e0.26}
1+\max_{s\in I_{b/2}}|g^{-1}(s)|<\frac{L_b}{2}.
\end{equation}
Since \(g\) is strictly increasing with \(g(0)=0\), for \(s\in I_b\) and \(\sigma\in\{-1,1\}\), we deduce that \(g(\sigma t+g^{-1}(s))=0\) if and only if \(t=-\sigma g^{-1}(s)\), and similarly, \(g(\sigma t)=s\) if and only if \(t=\sigma g^{-1}(s)\).
Now we apply Lemma~\ref{lemma:goodpoint} with \(b/2\) and \(L=L_b\).


First choose  \(0<\beta_b\le\eta_0\)
  where \(\eta_0\) is defined in \eqref{equation:Dzero}
such that
\[
C_{L_b}\beta_b^2\le\frac{\eta_0^2}{4},
\qquad
C_{L_b}\delta_k^2\le\frac{\eta_0^2}{4}.
\]
By \eqref{equation:macroflat}, we may also assume
\[
\sup_{\{u_k\in I_b\}\cap Q_{3/2}}|t|\le\frac12,
\qquad
L_b\varepsilon_k\le\frac12.
\]
Let \(X=(y,t)\in\{u_k\in I_b\}\cap Q_{3/2}\) satisfy
\(\gamma_k(X)\le\beta_b\). It follows from  Lemma~\ref{lemma:goodpoint} that
\[
X_0=(y,t_0),\qquad
u_k(X_0)=0,\qquad
|t_0-t|\le L_b\varepsilon_k,\qquad
|t_0|\leq|t|+|t_0-t|\le1.
\]
Moreover, by Lemma~\ref{lemma:transfer},
\[
\gamma_k(X_0)^2
\le C_{L_b}\bigl(\gamma_k(X)^2+\delta_k^2\bigr)
\le\frac{\eta_0^2}{2}.
\]
Therefore, we obtain
\[
y\in D_{k,0},\qquad
t_0=f_k(y),\qquad
|t-f_k(y)|\le L_b\varepsilon_k,
\]
which proves \eqref{equation:centraltozero}.


Next choose \(0<\alpha_b\le \eta_0\), with
\begin{equation}\label{e0.29}
C_{L_b}\alpha_b^2\le\frac{\beta_b^2}{4},
\qquad
C_{L_b}\delta_k^2\le\frac{\beta_b^2}{4}.
\end{equation}
By \eqref{equation:fsmall}, we may also assume
\begin{equation}\label{e0.27}
\|f_k\|_{L^\infty(B'_2)}+L_b\varepsilon_k<1.
\end{equation}
For \(y\in D_{k,b}\), set \(X_0:=(y,f_k(y))\). Then
$
u_k(X_0)=0$ and $\gamma_k(X_0)\le\alpha_b.$
Apply Lemma~\ref{lemma:goodpoint}   to \(X_0\) with \(L=L_b\), there exists \(\sigma\in\{-1,1\}\) such that
\begin{equation}\label{e0.28}
\|u_k(X_0+\varepsilon_k z)-g(\sigma z_4+a_{X_0})\|_{C^1(B_{L_b})}\le\rho,
\end{equation}
where \(a_{X_0}=g^{-1}(u_k(X_0))=0\). For any \(s\in I_b\), set \(z_4:=\sigma g^{-1}(s)\). Then, \eqref{e0.26} gives \(|z_4|<L_b/2\), so that \((0,z_4)\in B_{L_b}\).
Define \(X_s:=X_0+\varepsilon_k(0,z_4)=(y,f_k^s(y))\), the inequality \eqref{e0.28} yields \(u_k(X_s)=s\). Moreover,
\[
|f_k^s(y)-f_k(y)|=\varepsilon_k|g^{-1}(s)|\le L_b\varepsilon_k.
\]
Combining this with \eqref{e0.27} yields \(|f_k^s(y)|<1\).
Hence, by Lemma~\ref{lemma:transfer} and \eqref{e0.29}, we have
\[
\gamma_k(X_s)^2
\le C_{L_b}\bigl(\gamma_k(X_0)^2+\delta_k^2\bigr)
\le\frac{\beta_b^2}{2},
\]
which implies that $X_s$ is a $\beta_b$-good point.
Then, from \eqref{equation:globalcone}, we deduce that   for every $Z=(z',z_4)\in Q_{5/2}$ with $u_k(Z)=s$,
\[
|z_4-f_k^s(y)|\le\ell|z'-y|.
\]
Taking \(z'=y\) proves the uniqueness of $f_k^s(y)$, while taking
\(Z=(\widetilde y,f_k^s(\widetilde y))\) gives
\[
|f_k^s(y)-f_k^s(\widetilde y)|
\le\ell|y-\widetilde y|
\qquad
y,\widetilde y\in D_{k,b}.
\]
Finally, using \eqref{equation:micronondeg} at \(X_s\) yields
\[
\varepsilon_k|u_{k,t}(y,f_k^s(y))|\ge c_b>0.
\]
Then, we complete the proof of \eqref{equation:levelsnearzero}.


By Lemma~\ref{lemma:onesheet}, for all sufficiently large \(k\) and all \(y\in B'_2\), $u_k(y,-1)u_k(y,1)<0$.
Continuity and the definition of \(D_{k,b}\) give
\[
B'_2\setminus D_{k,b}
\subset
\pi\bigl(Q_{9/4}\cap\{u_k=0,\ \gamma_k>\alpha_b\}\bigr).
\]
Hence Lemma~\ref{lemma:badballs} yields
$
|B'_2\setminus D_{k,b}|
\le C\alpha_b^{-2}\delta_k^2
\le C_b\delta_k^2,
$
which proves \eqref{equation:badbaseb}.
Set
\[
\mathcal B_{k,b}:=
Q_{3/2}\cap\{u_k\in I_b,\ \gamma_k>\beta_b\}.
\]
Using Lemma~\ref{lemma:badballs} again, there exist
\(X_{ki}\in\mathcal B_{k,b}\) and
\(\rho_{ki}\ge\varepsilon_k\) such that
\[
\mathcal B_{k,b}\subset\bigcup_iB_{\rho_{ki}}(X_{ki}),
\qquad
\sum_i\rho_{ki}^3
\le C\beta_b^{-2}\delta_k^2
\le C_b\delta_k^2.
\]
Together with \eqref{equation:centraltozero}, we obtain
\[
\{u_k\in I_b\}\cap Q_{3/2}
\subset
N_{L_b\varepsilon_k}(\Gamma_k)
\cup\bigcup_iB_{\rho_{ki}}(X_{ki}),
\]
where $\Gamma_k:=\{(y,f_k(y)): y\in B_1'\}$, and then   \eqref{equation:completecover} holds.
Finally, \eqref{equation:macroflat} gives
$
\sup_i|(X_{ki})_4|
\le
\sup_{\{u_k\in I_b\}\cap Q_3}|t|
\to0$ as $k\to\infty.$


\end{proof}


By the McShane extension theorem, we extend \(f_k^s\) from the good set \(D_{k,b}\) to the whole ball \(B'_2\). We still denote the extension by \(f_k^s\), which are integrable over \(B'_2\times I_b\). Obviously, the extended functions \(f_k^s\) are uniformly Lipschitz continuous in $y$ and it follows from \eqref{equation:levelsnearzero} that $\sup_{s\in I_b}\|f_k^s\|_{L^\infty(B'_2)}\le s'_{k,b}\to0$.

\begin{lemma}\label{lemma:heightH1}
Assume \(\delta_k>0\) and \(\varepsilon_k/\delta_k\to0\). Then,
the fixed extension \(f_k\) satisfies
\begin{equation}\label{equation:zeroH1}
\int_{B'_2}|Df_k|^2\,dy\le C\delta_k^2.
\end{equation}
For every \(b\in(0,1)\), the extended functions \(f_k^s\) are uniformly Lipschitz continuous in $y$ with a constant \(C_b\), and satisfy
\begin{equation}\label{equation:levelH1}
\sup_{s\in I_b}\int_{B'_2}|Df_k^s|^2\,dy
\le C_b\delta_k^2.
\end{equation}
Moreover, up to a subsequence, there exists \(h\in H^1(B'_2)\)
such that
\begin{equation}\label{equation:heightlimit}
\frac{f_k-a_k}{\delta_k}\rightharpoonup h
\quad\text{in }H^1(B'_2),\qquad
\frac{f_k-a_k}{\delta_k}\longrightarrow h
\quad\text{in }L^2(B'_{3/2}),
\end{equation}
where $a_k:=\fint_{B'_2}f_k\,dy\to 0$ as $k\to\infty$.
Furthermore,
\[
\sup_{s\in I_b}
\left\|
\frac{f_k^s-a_k}{\delta_k}-h
\right\|_{L^2(B'_{3/2})}
\to0\]
uniformly in $s$.
\end{lemma}


\begin{proof}
For a.e. \(y\in D_{k,0}\cap B'_2\), set \(X:=(y,f_k(y))\).
The implicit function theorem and
\eqref{equation:micronondeg} give
\[
Df_k(y)=-\frac{D_Tu_k(X)}{u_{k,t}(X)},
\qquad
\varepsilon_k|u_{k,t}|\ge c_0>0
\quad\text{in }B_{c\varepsilon_k}(X).
\]
In particular, \(Du_k\ne0\) in this ball. For \(i=1,2,3\),
\(w_i:=\partial_i u_k\) satisfies
\[
\operatorname{div}\bigl(D^2\Phi_k(Du_k)Dw_i\bigr)
=\varepsilon_k^{-2}W''(u_k)w_i
\quad\text{in }B_{c\varepsilon_k}(X).
\]
We set $\widetilde w_i(z):=w_i(X+\varepsilon_kz),$ and $A_{k,X}(z):=D^2\Phi_k(Du_k(X+\varepsilon_kz)).$
Then
\begin{equation}\label{e0.30}
\operatorname{div}_z(A_{k,X}D_z\widetilde w_i)
=W''(u_k(X+\varepsilon_kz))\widetilde w_i
\quad\text{in }B_c,
\end{equation}
where  $\lambda I\le A_{k,X}\le\Lambda I,$ and $|W''(u_k(X+\varepsilon_kz))|\le\|W''\|_{L^\infty([-1,1])}.$

Let \(\varphi_k\) be the local implicit graph on a neighborhood
\(V\subset B'_2\) of \(y\). Then
\[
(f_k-\varphi_k)|_{D_{k,0}\cap V}=0,\]
which implies that  $Df_k=D\varphi_k$ a.e. on $D_{k,0}\cap V. $
By standard interior elliptic estimates,
\[
|D_Tu_k(X)|^2
\le C\varepsilon_k^{-4}
\int_{B_{c\varepsilon_k}(X)}|D_Tu_k|^2\,dx.
\]
Together with \eqref{equation:micronondeg} yields, for a.e.
\(y\in D_{k,0}\cap B'_2\),
\begin{equation}\label{equation:pointheightH1}
|Df_k(y)|^2
=\frac{|D_Tu_k(X)|^2}{|u_{k,t}(X)|^2}
\le C\varepsilon_k^2|D_Tu_k(X)|^2
\le C\varepsilon_k^{-3}
\int_{B_{c\varepsilon_k}(y,f_k(y))}
\varepsilon_k|D_Tu_k|^2\,dx.
\end{equation}
Therefore, by Fubini theorem,
\[
\begin{aligned}
\int_{D_{k,0}}|Df_k|^2\,dy
&\le C\varepsilon_k^{-3}
\int_{D_{k,0}}\int_{B_{c\varepsilon_k}(y,f_k(y))}
\varepsilon_k|D_Tu_k(x)|^2\,dx\,dy\\
&\le C\varepsilon_k^{-3}
\int_{Q_4}\varepsilon_k|D_Tu_k(x)|^2
|B'_{c\varepsilon_k}(x')|\,dx
\le C\delta_k^2.
\end{aligned}
\]
Since $f_k$ is Lipschitz in $B_2'$ and \(|B'_2 \setminus D_{k,0}| \le C\delta_k^2\) by Lemma~\ref{lemma:badballs}, we have
\[
\int_{B'_2\setminus D_{k,0}}|Df_k|^2\,dy
\le C\delta_k^2,
\]
where $C$ depend on $\eta_0$ and the Lipschitz constant of $f_k$. Combining with above two inequalities, we prove \eqref{equation:zeroH1}.

By \eqref{equation:levelsnearzero}, the same argument at
\((y,f_k^s(y))\) yields
\[
\sup_{s\in I_b}\int_{D_{k,b}}|Df_k^s|^2\,dy
\le C_b\delta_k^2.
\]
For the uniformly Lipschitz extensions, \eqref{equation:badbaseb}
then gives
\[
\begin{aligned}
\sup_{s\in I_b}\int_{B'_2}|Df_k^s|^2\,dy
&\le \sup_{s\in I_b}\int_{D_{k,b}}|Df_k^s|^2\,dy
+C_b|B'_2\setminus D_{k,b}|\le C_b\delta_k^2,
\end{aligned}
\]
which proves \eqref{equation:levelH1}.


On the compact set \(D_{k,b}\times I_b\), the map \((y,s)\mapsto f_k^s(y)\) is continuous. Indeed, given any sequence \((y_j,s_j)\to(y,s)\), the heights \(t_j:=f_k^{s_j}(y_j)\) satisfy \(u_k(y_j,t_j)=s_j\). By the boundedness of \(t_j\) and the continuity of \(u_k\), up to a  subsequence, \(t_{j}\to t\) satisfies \(u_k(y,t)=s\). Then, the uniqueness of the level height at \(y\in D_{k,b}\) forces \(t=f_k^s(y)\). Thus, $f_k^s(y)$ are continuous on $D_{k,b}\times I_b$

Define
\[
\widehat f_k^s(y)
:=\inf_{z\in D_{k,b}}
\bigl\{f_k^s(z)+\ell|y-z|\bigr\},\quad\text{and}\quad
\widetilde f_k^s(y)
:=\max\bigl\{-s'_{k,b},
\min\{s'_{k,b},\widehat f_k^s(y)\}\bigr\},
\]
where $s'_{k,b}:=\max_{(y,s)\in D_{k,b}\times I_b}|f_k^s(y)|
\to0$ by \eqref{equation:macroflat}.
These extensions are  continuous in \((y,s)\) and satisfy
\[
\widetilde f_k^s|_{D_{k,b}}=f_k^s,\qquad
\sup_{s\in I_b}\operatorname{Lip}(\widetilde f_k^s;B'_2)\le\ell,
\qquad
\sup_{s\in I_b}\|\widetilde f_k^s\|_{L^\infty(B'_2)}
\le s'_{k,b}.
\]
We still denote these extensions by \(f_k^s\) on \(B'_2\).


Set \(h_k:=(f_k-a_k)/\delta_k\). By Poincar\'e,
\eqref{equation:zeroH1}, and \eqref{equation:fsmall},
\[
\|h_k\|_{H^1(B'_2)}^2
\le C\delta_k^{-2}\int_{B'_2}|Df_k|^2\,dy\le C,
\qquad
|a_k|\le\|f_k\|_{L^\infty(B'_2)}\to0.
\]
After passing to a subsequence, by weak compactness and Rellich theorem, we obtain
\[
h_k\rightharpoonup h\quad\text{in }H^1(B'_2),
\qquad
h_k\to h\quad\text{in }L^2(B'_{3/2}),
\]
which proves \eqref{equation:heightlimit}. It follows from \eqref{equation:levelsnearzero} that
\eqref{equation:badbaseb},
\[\begin{split}
\sup_{s\in I_b}
\left\|\frac{f_k^s-f_k}{\delta_k}\right\|_{L^2(B'_{3/2})}^2
&\le C_b\frac{\varepsilon_k^2}{\delta_k^2}
+\frac{\bigl(s'_{k,b}+\|f_k\|_{L^\infty(B'_2)}\bigr)^2}
{\delta_k^2}|B'_2\setminus D_{k,b}|\notag\\
&\le C_b\left(\frac{\varepsilon_k}{\delta_k}\right)^2
+C_b\bigl(s'_{k,b}+\|f_k\|_{L^\infty(B'_2)}\bigr)^2
\longrightarrow0,
\end{split}\]
where \(\varepsilon_k/\delta_k\to0\). Consequently,
\begin{equation}\label{equation:commonheight}
\sup_{s\in I_b}
\left\|\frac{f_k^s-a_k}{\delta_k}-h\right\|_{L^2(B'_{3/2})}
\le
\sup_{s\in I_b}
\left\|\frac{f_k^s-f_k}{\delta_k}\right\|_{L^2(B'_{3/2})}
+\|h_k-h\|_{L^2(B'_{3/2})}
\longrightarrow0.
\end{equation}
Set \(v_k(y,s):=(f_k^s(y)-a_k)/\delta_k\).
Apply \eqref{equation:levelH1} and \eqref{equation:commonheight}  on \(B'_2\), we have
\[
\int_{I_b}\int_{B'_2}|D_yv_k|^2\,dy\,ds\le C_b,
\qquad
v_k\rightharpoonup h(y)
\quad\text{in }L^2(B'_2\times I_b).
\]
For every
\(\varphi\in C_c^\infty(B'_2\times(-1+b,1-b);\mathbb R^3)\),
\[
\lim_{k\to\infty}\iint D_yv_k\cdot\varphi
=-\lim_{k\to\infty}\iint v_k\operatorname{div}_y\varphi
=\iint Dh\cdot\varphi.
\]
Hence
$
 D_yf_k^s/\delta_k \rightharpoonup Dh(y)$
 in  $L^2(B'_2\times I_b;\mathbb R^3).$
By a diagonal argument for \(b_m=1/(m+1)\)  and passing to the limit, we obtain
\[
\frac{D_yf_k^s}{\delta_k }\rightharpoonup Dh(y)
\quad\text{in } L^2(B'_2\times (-1,1);\mathbb R^3).
\]
Finally, \eqref{equation:zeroH1} and Poincar\'e give
\[
\|h\|_{H^1(B'_2)}
\le\liminf_{k\to\infty}
\left\|\frac{f_k-a_k}{\delta_k}\right\|_{H^1(B'_2)}
\le C,
\]
with \(C\) independent of \(b\).

\end{proof}

\begin{lemma}\label{lemma:offsheet}
For \(b\in(0,1)\), set
$
A_{k,b}
:=\{(y,t)\in Q_1:\ u_k(y,t)\in I_b,\ y\notin D_{k,b}\}.
$
There exists \(C_b>0\) such that, for all sufficiently large \(k\),
\begin{equation}\label{e0.31}
\int_{A_{k,b}}e_k\,dx\le C_b\delta_k^2.
\end{equation}
\end{lemma}


\begin{proof}
Set the good set $G_{k,b}:=A_{k,b}\cap\{\gamma_k\le\beta_b\}.$
By Lemma~\ref{lemma:badballs},
\[
\int_{A_{k,b}\setminus G_{k,b}}e_k\,dx
\le C\beta_b^{-2}\delta_k^2
\le C_b\delta_k^2.
\]
By Lemma~\ref{lemma:goodpoint}, we get
\eqref{e0.22} and \eqref{equation:micronondeg},
\[
\varepsilon_k|Du_k|\le C,
\qquad
\varepsilon_k|u_{k,t}|\ge c_b
\quad\text{on }G_{k,b}.
\]
Since \(\Phi_k(p)\le C|p|^2\) and \(|u_k|<1\),
\[
e_k
=\varepsilon_k\Phi_k(Du_k)
+\frac{W(u_k)}{\varepsilon_k}
\le\frac{C}{\varepsilon_k},
\qquad
\frac{e_k}{|u_{k,t}|}\le C_b
\quad\text{on }G_{k,b}.
\]
Then,  applying  the area formula,
we obtain
\[
\int_{G_{k,b}}e_k\,dx
=
\int_{I_b}\int_{B'_1\setminus D_{k,b}}
\sum_{\substack{u_k(y,t)=s\\(y,t)\in G_{k,b}}}
\frac{e_k(y,t)}{|u_{k,t}(y,t)|}\,dy\,ds.
\]
By \eqref{equation:globalcone}, any two points
\((y,t_1),(y,t_2)\in G_{k,b}\) satisfy
\[
|t_1-t_2|\le\ell|y-y|=0,
\]
which implies \(t_1=t_2\).
The uniqueness of $t$ and \eqref{equation:badbaseb} give
\[
\int_{G_{k,b}}e_k\,dx
\le C_b|I_b|\,|B'_1\setminus D_{k,b}|
\le C_b\delta_k^2.
\]
Combining the two parts, \eqref{e0.31} holds.

\end{proof}

\section{Convergence of the Height Functions}\label{section:flux}

In this section, we retain the assumptions and notation of Section~5, with
  \(\Phi_k:=\frac12F_k^2\), $Du_k=(D_Tu_k,u_{k,4})$, $D\Phi_k=(D_T\Phi_k,\Phi_{k,4})$ and denote $A_k:=D\Phi_k,$ $A_T:=D_T\Phi_k$.

\begin{lemma}\label{lemma:algebra}
Let $F(e_4)=1$, $DF(e_4)=e_4$.
Then, for every $p=(p_T,p_4)$,
\begin{align}
 \lambda|p_T|^2\le A_T(p)\cdot p_T&\le\Lambda|p_T|^2,
 \qquad |A_T(p)| \le\Lambda|p_T|,\label{equation:coercivity}\\
 \frac{\lambda}{2}|p_T|^2\le\Phi(p)-\frac{1}{2}p_4^2
      &\le\frac{\Lambda}{2}|p_T|^2,\label{equation:PhiSplit}\\
 |p_4|\,|A_T(p)-D^2_{TT}\Phi(e_4)p_T|&\le C|p_T|^2.\label{equation:fluxremainder}
\end{align}
\end{lemma}
\begin{proof}
By normalization and evenness, $A_T(0,p_4)=0$ and $\Phi(0,p_4)=\frac12p_4^2.$
For \(p\ne0\), the global Lipschitz continuity of \(A\) gives
\[
A_T(p)
=\int_0^1D^2_{TT}\Phi(sp_T,p_4)p_T\,ds.
\]
Since $
\lambda I\le D^2_{TT}\Phi(sp_T,p_4)\le\Lambda I,$
we obtain \eqref{equation:coercivity}. Moreover,
\[
\begin{aligned}
\Phi(p)-\frac12p_4^2
&=\int_0^1 A_T(sp_T,p_4)\cdot p_T\,ds\\
&=\int_0^1(1-s)
\bigl\langle D^2_{TT}\Phi(sp_T,p_4)p_T,p_T\bigr\rangle\,ds,
\end{aligned}
\]
which yields \eqref{equation:PhiSplit}.
The case \(p=0\) is trivial.


For \(p_4=0\), \eqref{equation:fluxremainder} is immediate.
By oddness and homogeneity, $A_T(p)=p_4A_T(p_T/p_4,1)$ for any $p_4\ne0.$
Taylor expansion at \(e_4\) gives
\[
|A_T(a,1)-D^2_{TT}\Phi(e_4)a|\le C|a|^2
\qquad\text{for }|a|\le\frac12.
\]
Hence, if \(0<|p_4|\) and \(|p_T|\le|p_4|/2\),
\[
|p_4|\,|A_T(p)-D^2_{TT}\Phi(e_4)p_T|
=|p_4|^2
\left|A_T\!\left(\frac{p_T}{p_4},1\right)
-D^2_{TT}\Phi(e_4)\frac{p_T}{p_4}\right|
\le C|p_T|^2.
\]
If \(|p_T|>|p_4|/2\), together with \eqref{equation:coercivity} and
\(\|D^2_{TT}\Phi(e_4)\|\le\Lambda\), we obtain
\[
|p_4|\,|A_T(p)-Bp_T|
\le2\Lambda|p_4||p_T|
\le4\Lambda|p_T|^2.
\]
This completes the proof of \eqref{equation:fluxremainder} for all  \(p\in\mathbb R^4\).


\end{proof}

\begin{proposition}\label{proposition:flux}
Along the subsequence in Lemma~\ref{lemma:heightH1},
\begin{align}
 \varepsilon_k u_{k,4}^2\,dx
 &\rightharpoonup
 c_W\mathcal H^3\mathbin{\llcorner}\{t=0\},
 \label{equation:normaldensity}\\
 \frac{\varepsilon_k u_{k,4}D_Tu_k}{\delta_k}\,dx
 &\rightharpoonup
 -c_WDh\,\mathcal H^3\mathbin{\llcorner}\{t=0\},
 \label{equation:fluxlimit}
\end{align}
weakly as Radon measures in $Q_1$.
\end{proposition}

\begin{proof}

By Lemma~\ref{lemma:onesheet}, $u_k(\cdot,\pm2)$ converge to $\pm\sigma$ uniformly on $B'_2$, where  $\sigma\in\{-1,1\}.$
With $\Psi(s):=\int_{-1}^s\sqrt{2W(r)}\,dr$, we obtain
\begin{equation}\label{e0.32}
\liminf_{k\to\infty}
\int_{Q_2}\sqrt{2W(u_k)}\,|u_{k,4}|\,dx
\ge \lim_{k\to\infty}
\int_{B'_2}
\bigl|\Psi(u_k(y,2))-\Psi(u_k(y,-2))\bigr|\,dy=c_W|B'_2|.
\end{equation}
By \eqref{equation:PhiSplit}, \eqref{e0.32} and
Lemma~\ref{lemma:onesheet}, we obtain
\begin{equation}\label{equation:unscaleddefect}
0
\le \int_{Q_2}
\frac{\bigl(\varepsilon_k|u_{k,4}|-\sqrt{2W(u_k)}\bigr)^2}
     {2\varepsilon_k}\,dx
\leq \int_{Q_2}
\bigl(e_k-\sqrt{2W(u_k)}\,|u_{k,4}|\bigr)\,dx
\longrightarrow0.
\end{equation}
The energy bound and Cauchy-Schwarz give
\[
\int_{Q_2}
\left|\varepsilon_k u_{k,4}^2
-\sqrt{2W(u_k)}\,|u_{k,4}|\right|\,dx
 \le
\left(\int_{Q_2}
\frac{\bigl(\varepsilon_k|u_{k,4}|-\sqrt{2W(u_k)}\bigr)^2}
     {\varepsilon_k}\,dx\right)^{1/2}
\left(\int_{Q_2}\varepsilon_k u_{k,4}^2\,dx\right)^{1/2}
\longrightarrow0.
\]
Thus
$
 e_k-\varepsilon_k u_{k,4}^2 \to0$ in $L^1(Q_2)$ and Lemma~\ref{lemma:onesheet} yields
\eqref{equation:normaldensity}.


Fix $b\in(0,1)$ and define the good set
\[
\mathcal G_{k,b}
:=\{(y,f_k^s(y)):y\in D_{k,b}\cap B'_2,\ s\in I_b\}
\subset Q_2
\]
for all sufficiently large $k$. Set
$
b_k(y,s):=\varepsilon_k
|u_{k,4}(y,f_k^s(y))|.$
By Proposition~\ref{proposition:common} and
Lemma~\ref{lemma:goodpoint}, $c_b\le b_k\le C_b.$
Differentiating the identity $u_{k,4}(y,f_k^s(y))=s$ with respect to $y$ and applying the chain rule, we obtain $D_Tu_k=-u_{k,4}Df_k^s$ almost everywhere.
Thus coarea and \eqref{equation:unscaleddefect} yield
\begin{equation}\label{equation:weightdefect}
\begin{aligned}
\int_{I_b}\int_{D_{k,b}}
\frac{|b_k-\sqrt{2W(s)}|^2}{b_k}\,dy\,ds
&=
\int_{\mathcal G_{k,b}}
\frac{\bigl(\varepsilon_k|u_{k,4}|
-\sqrt{2W(u_k)}\bigr)^2}{\varepsilon_k}\,dx
\longrightarrow0,\\
\int_{I_b}\int_{D_{k,b}}
b_k|Df_k^s|^2\,dy\,ds
&=
\int_{\mathcal G_{k,b}}
\varepsilon_k|D_Tu_k|^2\,dx
\le\delta_k^2.
\end{aligned}
\end{equation}
For $\varphi\in C_c^\infty(Q_1)$, coarea formula gives
\begin{equation}\label{equation:coareaflux}
\int_{\mathcal G_{k,b}}
\frac{\varepsilon_k u_{k,4}D_Tu_k}{\delta_k}\,
\varphi\,dx
=
-\int_{I_b}\int_{D_{k,b}}
b_k\frac{Df_k^s}{\delta_k}\,
\varphi(y,f_k^s(y))\,dy\,ds.
\end{equation}
By Cauchy-Schwarz inequality and \eqref{equation:weightdefect}, we get
\begin{equation}\label{equation:replaceweight}
\begin{aligned}
&\left|
\int_{I_b}\int_{D_{k,b}}
\bigl(b_k-\sqrt{2W(s)}\bigr)
\frac{Df_k^s}{\delta_k}\,
\varphi(y,f_k^s(y))\,dy\,ds
\right|\\
&\quad\le \|\varphi\|_\infty
\left(
\int_{I_b}\int_{D_{k,b}}
\frac{|b_k-\sqrt{2W(s)}|^2}{b_k}\,dy\,ds
\right)^{1/2}
\longrightarrow0.
\end{aligned}
\end{equation}
By Lemma~\ref{lemma:heightH1}, we have
\[
 \sup_{s\in I_b}\|f_k^s\|_{L^\infty(B'_2)}
 \le s'_{k,b}\to0,
 \qquad
 \int_{I_b}\int_{B'_2}
 \frac{|Df_k^s|^2}{\delta_k^2}\,dy\,ds\le C_b.
\]
The mean value theorem gives
$
 |\varphi(y,f_k^s(y))-\varphi(y,0)|
 \le \|\partial_4\varphi\|_\infty |f_k^s(y)|
 \le \|\partial_4\varphi\|_\infty s'_{k,b}.
$
Hence, by Cauchy-Schwarz inequality,
\begin{equation}\label{e0.33}
\begin{aligned}
&\left|
 \int_{I_b}\int_{D_{k,b}}
 \sqrt{2W(s)}\,\frac{Df_k^s}{\delta_k}
 \bigl[\varphi(y,f_k^s(y))-\varphi(y,0)\bigr]\,dy\,ds
 \right|\\
&\quad\le C_{b,\varphi}s'_{k,b}
 \int_{I_b}\int_{D_{k,b}}
 \frac{|Df_k^s|}{\delta_k}\,dy\,ds\\
&\quad\le C_{b,\varphi}s'_{k,b}
 \bigl(|I_b|\,|D_{k,b}|\bigr)^{1/2}
 \left(
 \int_{I_b}\int_{D_{k,b}}
 \frac{|Df_k^s|^2}{\delta_k^2}\,dy\,ds
 \right)^{1/2}\le C_{b,\varphi}s'_{k,b}\longrightarrow0.
\end{aligned}
\end{equation}
Moreover, using  Cauchy-Schwarz inequality again, \eqref{equation:badbaseb} and
\eqref{equation:levelH1} give
\begin{equation}\label{e0.34}
\begin{aligned}
&\left|
\int_{I_b}\int_{B'_2\setminus D_{k,b}}
\sqrt{2W(s)}\,\frac{Df_k^s}{\delta_k}
\varphi(y,0)\,dy\,ds
\right|\\
&\quad\le C_{b,\varphi}|B'_2\setminus D_{k,b}|^{1/2}
\left(
\int_{I_b}\int_{B'_2}
\frac{|Df_k^s|^2}{\delta_k^2}\,dy\,ds
\right)^{1/2}
\le C_{b,\varphi}\delta_k\longrightarrow0.
\end{aligned}
\end{equation}
Therefore, by Lemma~\ref{lemma:heightH1} and \eqref{equation:coareaflux}-\eqref{e0.34}, we obtain
\begin{equation}\label{e0.35}
\lim_{k\to\infty}
\int_{\mathcal G_{k,b}}
\frac{\varepsilon_k u_{k,4}D_Tu_k}{\delta_k}\,
\varphi\,dx
=
-\left(\int_{I_b}\sqrt{2W(s)}\,ds\right)
\int_{B'_1}Dh(y)\varphi(y,0)\,dy.
\end{equation}

\noindent On the other hand, by Cauchy-Schwarz inequality and Lemma~\ref{lemma:offsheet},
\begin{equation}\label{e0.36}
\begin{aligned}
\frac1{\delta_k}\int_{A_{k,b}}
\varepsilon_k|u_{k,4}|\,|D_Tu_k|\,dx
&\le \frac1{\delta_k}
\left(\int_{A_{k,b}}\varepsilon_k|D_Tu_k|^2\,dx\right)^{1/2}
\left(\int_{A_{k,b}}\varepsilon_k u_{k,4}^2\,dx\right)^{1/2}\\
&\le C\left(\int_{A_{k,b}}e_k\,dx\right)^{1/2}
\le C_b\delta_k\longrightarrow0.
\end{aligned}
\end{equation}


For fixed $b\in(0,1)$, uniform convergence at $t=\pm2$ implies that
for all sufficiently large $k$,
\[
\int_{Q_2\cap\{u_k\in I_b\}}e_k\,dx
\ge
\int_{Q_2\cap\{u_k\in I_b\}}
\sqrt{2W(u_k)}\,|u_{k,4}|\,dx
\ge |B'_2|\int_{I_b}\sqrt{2W(s)}\,ds.
\]
By Lemma~\ref{lemma:onesheet}, we have
\begin{equation}
\limsup_{k\to\infty}
\int_{Q_2\cap\{|u_k|>1-b\}}e_k\,dx
\le
|B'_2|\left(c_W-\int_{I_b}\sqrt{2W(s)}\,ds\right).
\end{equation}
The same Cauchy-Schwarz estimate yields
\begin{equation}\label{e0.37}
\limsup_{k\to\infty}
\frac1{\delta_k}
\int_{Q_2\cap\{|u_k|>1-b\}}
\varepsilon_k|u_{k,4}|\,|D_Tu_k|\,dx
\le
C\left(c_W-\int_{I_b}\sqrt{2W(s)}\,ds\right)^{1/2}
\end{equation}
Since $\int_{I_b}\sqrt{2W(s)}\,ds\to c_W$ as $b\to0$, it follows from \eqref{e0.35}-\eqref{e0.37} that \eqref{equation:fluxlimit} holds.

\end{proof}

\begin{corollary}\label{corollary:varifoldplane}
The varifolds $V_k$ defined in \eqref{equation:diffuseV} converge locally to $c_W|\{t=0\}|$ in $Q_1.$
\end{corollary}

\begin{proof}
By the triangle inequality and the Lipschitz continuity of $F$, we obtain
\begin{align*}
\bigl||Du_k|F_k(Du_k)-u_{k,4}^2\bigr|
&\le |Du_k| \bigl|F_k(Du_k) - u_{k,4}\bigr|
+ u_{k,4} \bigl||Du_k| - u_{k,4}\bigr|\\
&\leq  C (|D_Tu_k| + u_{k,4}) |D_Tu_k|+ u_{k,4}\frac{|D_Tu_k|^2}{|Du_k| + u_{k,4}}\\
&\leq  C\bigl(|D_Tu_k|^2+|u_{k,4}|\,|D_Tu_k|\bigr).
\end{align*}
Hence, by Cauchy-Schwarz inequality, we get
\begin{equation}\label{e0.39}
\int_{Q_1}\varepsilon_k
\bigl||Du_k|F_k(Du_k)-u_{k,4}^2\bigr|\,dx
\le C\delta_k^2
+C\delta_k
\left(\int_{Q_1}\varepsilon_k u_{k,4}^2\,dx\right)^{1/2}
\longrightarrow0.
\end{equation}
Moreover, by the definition of the varifold  $V_k$, we have
\begin{equation}\label{e0.38}
\int_{Q_1\cap\{Du_k\ne0\}}
\left(1-\frac{u_{k,4}^2}{|Du_k|^2}\right)\,d\|V_k\|
=
\int_{Q_1\cap\{Du_k\ne0\}}
\varepsilon_k\frac{F_k(Du_k)}{|Du_k|}
|D_Tu_k|^2\,dx
\le C\delta_k^2\longrightarrow0.
\end{equation}
For \(\psi\in C_c^1(Q_1\times G(3))\), the preceding
tilt estimate \eqref{e0.38} and Cauchy-Schwarz give
\begin{equation}\label{e0.40}
\begin{aligned}
 &\left|
 \int_{Q_1\times G(3)}\psi(x,P)\,dV_k(x,P)
 -\int_{Q_1}\psi(x,e_4^\perp)\,d\|V_k\|(x)
 \right|\\
 &\quad\le C_\psi
 \int_{Q_1\cap\{Du_k\ne0\}}
 \left(1-\frac{u_{k,4}^2}{|Du_k|^2}\right)^{1/2}
 \,d\|V_k\|\\
 &\quad\le C_\psi\|V_k\|(Q_1)^{1/2}
 \left(
 \int_{Q_1\cap\{Du_k\ne0\}}
 \left(1-\frac{u_{k,4}^2}{|Du_k|^2}\right)
 \,d\|V_k\|
 \right)^{1/2}\le C_\psi\delta_k\longrightarrow0.
\end{aligned}
\end{equation}
The   estimate \eqref{e0.39} and \eqref{equation:normaldensity} therefore yield
\begin{equation}\label{e0.41}
 \|V_k\|\rightharpoonup
 c_W\mathcal H^3\mathbin{\llcorner}\{t=0\}
 \quad\text{in }Q_1.
\end{equation}
Then, it follows from \eqref{e0.40} and \eqref{e0.41} that
\begin{align*}
\lim_{k\to\infty}\frac{1}{c_W}\int_{Q_1\times G(3)}\psi(x,P)\,dV_k(x,P)
&= \lim_{k\to\infty}\frac{1}{c_W}\int_{Q_1}\psi(x,e_4^\perp)\,d\|V_k\|(x) \\
&= \int_{Q_1}\psi(x,e_4^\perp)\,d(\mathcal{H}^3\llcorner\{t=0\})
= \int_{B'_1}\psi((y,0),e_4^\perp)\,dy.
\end{align*}
For any \(\psi\in C_c(Q_1\times G(3))\), choose
\(\psi_j\in C_c^1(Q_1\times G(3))\) with
\(\|\psi_j-\psi\|_\infty\to0\).
Since \(\sup_k\|V_k\|(Q_1)<\infty\), we obtain
\[
 \limsup_{k\to\infty}
 \left|
 \frac1{c_W}\int_{Q_1\times G(3)}\psi\,dV_k
 -\int_{B'_1}\psi((y,0),e_4^\perp)\,dy
 \right|
 \le C\|\psi-\psi_j\|_\infty.
\]
Letting \(j\to\infty\), we conclude
\[
 c_W^{-1}V_k\rightharpoonup|\{t=0\}|
 \qquad\text{locally in }Q_1.
\]

\end{proof}


\begin{lemma}\label{lemma:harmonic}
Let \(h\) be as in Lemma~\ref{lemma:heightH1}. After a subsequence, $B_k=D^2_{TT}\Phi_k(e_4)\to B_\infty$, and  satisfies
\begin{equation}\label{equation:harmonic}
\operatorname{div}(B_\infty Dh)=0
\qquad\text{in }B'_1.
\end{equation}
Moreover,
$
|Dh(0)|+\|D^2h\|_{L^\infty(B'_{1/4})}\le C,
$
where $C$ is independent of \(b\in(0,1)\).
\end{lemma}


\begin{proof}
Let $\varphi\in C_c^\infty(B'_1)$ and choose
$\chi\in C_c^\infty((-1,1))$ with
$\chi=1$ on $[-1/2,1/2]$.
Taking $X=\varphi(y)\chi(t)e_4$ in \eqref{equation:stress} gives
\begin{equation}\label{e0.42}
\int_{Q_1}\chi\varepsilon_k u_{k,4}
D_T\Phi_k(Du_k)\cdot D\varphi\,dx
=
\int_{Q_1}\varphi\chi'
\bigl[e_k-\varepsilon_k u_{k,4}\Phi_{k,4}(Du_k)\bigr]\,dx.
\end{equation}
Since $\operatorname{spt}\chi'\subset\{|t|\ge1/2\}$,   from
Lemma~\ref{lemma:modicaclear} and the fact that $\varepsilon_k/\delta_k\to0$, we have
\begin{equation}\label{e0.43}
\frac1{\delta_k}\left|
\int_{Q_1}\varphi\chi'
\bigl[e_k-\varepsilon_k u_{k,4}\Phi_{k,4}(Du_k)\bigr]\,dx
\right|
\quad\le
\frac{C_{\varphi,\chi}}{\varepsilon_k\delta_k}
e^{-c/\varepsilon_k}
\le C_{\varphi,\chi}\varepsilon_k^{-2}
e^{-c/\varepsilon_k}
\longrightarrow0.
\end{equation}


\noindent By \eqref{equation:fluxremainder}, we have
\[
 \frac1{\delta_k}\int_{Q_1}
\varepsilon_k|u_{k,4}|
\bigl|D_T\Phi_k(Du_k)-B_kD_Tu_k\bigr|\,dx
 \le
\frac{C}{\delta_k}\int_{Q_1}
\varepsilon_k|D_Tu_k|^2\,dx
\le C\delta_k\longrightarrow0.
\]
After passing to a subsequence, $B_k\to B_\infty$ with $\lambda I\le B_\infty\le\Lambda I.$
It follows from  \eqref{e0.42}, \eqref{e0.43} and \eqref{equation:fluxlimit} that
\[
\int_{B'_1}B_\infty Dh\cdot D\varphi\,dy=0
\qquad\text{for every }\varphi\in C_c^\infty(B'_1).
\]
Interior estimates and Lemma~\ref{lemma:heightH1} therefore yield
\begin{equation}\label{equation:hestimate}
|Dh(0)|+\|D^2h\|_{L^\infty(B'_{1/4})}
\le C\|Dh\|_{L^2(B'_1)}
\le C.
\end{equation}

\end{proof}


\begin{lemma}\label{lemma:weightedheight}
Let $a_k$ and $h$ be as in Lemma~\ref{lemma:heightH1},
and assume
$
\delta_k>0,$ $\varepsilon_k/\delta_k\longrightarrow0.$
Along the same subsequence, for every compact $K\subset\subset Q_{1/2}$,
\begin{equation}\label{equation:weightedheight}
\lim_{k\to\infty}
\frac1{\delta_k^2}
\int_K e_k(x)\,
|t-a_k-\delta_kh(y)|^2\,dx=0.
\end{equation}
\end{lemma}


\begin{proof}

Let $b_w$ be as in Lemma~\ref{lemma:modicaclear} and set
\[
\Sigma_k:=\{|u_k|\le1-b_w\}\cap Q_{3/2}.
\]
Since $\operatorname{dist}(Q_{1/2},\partial Q_{3/2})=1$,
\eqref{equation:clearing} yields
\begin{equation}\label{e0.46}
e_k(x)\le\frac{C}{\varepsilon_k}
\left(
e^{-c\,\operatorname{dist}(x,\Sigma_k)/\varepsilon_k}
+e^{-c/\varepsilon_k}
\right)
\qquad \text{in}\quad Q_{1/2}.
\end{equation}
By Proposition~\ref{proposition:common},
\begin{equation}\label{e0.44}
\Sigma_k\subset
N_{L_{b_w}\varepsilon_k}(\Gamma_k)
\cup\bigcup_i B_{\rho_{ki}}(X_{ki}),
\end{equation}
where     $\Gamma_k:=\{(y,f_k(y)): y\in B_1'\}$ ,
$
\sum_i\rho_{ki}^3\le C\delta_k^2,$ and
$\rho_k^*:=\sup_i\rho_{ki}
\le C\delta_k^{2/3}\to0.$

By Lemma~\ref{lemma:harmonic} and standard interior estimates, \(\|h\|_{W^{1,\infty}(B'_{3/4})} \le C\).
We extend \(h\) by a standard cut-off argument to a smooth function with compact support in \(\mathbb{R}^3\),
still denoted by \(h\), such that \(\|h\|_{W^{1,\infty}(\mathbb{R}^3)} \le C\).
Set
\[
\ell_k(y):=a_k+\delta_kh(y).
\]
Choose $s_k\to0$ such that
\[
\max\left\{
\|f_k\|_{L^\infty(B'_2)},
\|\ell_k\|_{L^\infty(\mathbb R^3)},
|a_k|
\right\}
+\sup_{(y,t)\in\Sigma_k}|t|
+\sup_i|X_{ki,4}|
\le s_k.
\]


Let \(U_k:=\bigcup_iB_{2\rho_{ki}}(X_{ki})\). Since $|t - \ell_k(y)| \le |t - X_{ki,4}| + |X_{ki,4}| + |\ell_k(y)| \le   C(s_k+\rho_{ki}+\delta_k)$ on $B_{2\rho_{ki}}(X_{ki})$,
\eqref{equation:growth} gives
\begin{equation}\label{e0.50}
\begin{split}
 \int_{K\cap U_k}e_k|t-\ell_k(y)|^2\,dx
 &\le C\sum_i(s_k+\rho_{ki}+\delta_k)^2\rho_{ki}^3 \\
 &\le C(s_k+\rho_k^*+\delta_k)^2\delta_k^2
 =o(\delta_k^2).
\end{split}
\end{equation}


For $x\in K\setminus U_k$, the relation \eqref{e0.44} gives
\begin{equation}\label{e0.47}
e^{-c\,\operatorname{dist}(x,\Sigma_k)/\varepsilon_k}
\le
e^{cL_{b_w}}
e^{-c\,\operatorname{dist}(x,\Gamma_k)/\varepsilon_k}
+\sum_i
e^{-c\,\operatorname{dist}(x,B_{\rho_{ki}}(X_{ki}))/\varepsilon_k}.
\end{equation}
Extend $f_k$ to $\mathbb R^3$ such that $\|f_k\|_{W^{1,\infty}}(\mathbb{R}^3)\leq C$. Then
$
\operatorname{dist}((y,t),\Gamma_k)
\ge c|t-f_k(y)|.
$
With the change of variables \(z=(t-f_k(y))/\varepsilon_k\), we obtain, for every $c>0$,
\begin{equation}\label{equation:graphkernel}
\begin{aligned}
 &\int_{\mathbb R}\varepsilon_k^{-1}
 e^{-c|t-f_k(y)|/\varepsilon_k}|t-\ell_k(y)|^2\,dt\\
 &\quad=\int_{\mathbb R}e^{-c|z|}
 \bigl(f_k(y)-\ell_k(y)+\varepsilon_k z\bigr)^2\,dz\\
 &\quad=2|f_k(y)-\ell_k(y)|^2
 \int_0^\infty e^{-cz}\,dz
 +2\varepsilon_k^2\int_0^\infty z^2e^{-cz}\,dz\\
 &\quad=\frac2c|f_k(y)-\ell_k(y)|^2
 +\frac4{c^3}\varepsilon_k^2 \le C\bigl(|f_k(y)-\ell_k(y)|^2+\varepsilon_k^2\bigr).
\end{aligned}
\end{equation}
By \eqref{equation:heightlimit} and \eqref{equation:epsdelta},
\begin{equation}\label{e0.45}
\begin{aligned}
\int_{B'_{1/2}}
\bigl(|f_k-\ell_k|^2+\varepsilon_k^2\bigr)\,dy
&=
\delta_k^2\int_{B'_{1/2}}
\left|\frac{f_k-a_k}{\delta_k}-h\right|^2\,dy
+|B'_{1/2}|\varepsilon_k^2=o(\delta_k^2).
\end{aligned}
\end{equation}


For
$s=|x-X_{ki}|\ge2\rho_{ki}$, by triangle inequality, we have
\[
|t-\ell_k(y)|   \le |t - X_{ki, 4}| + |X_{ki, 4}| + |\ell_k(y)|\le s+2s_k.
\]
By the substitution \(w = (s-\rho_{ki})/\varepsilon_k\), we have
\begin{align}\label{equation:ballmoments3}
\int_{2\rho_{ki}}^{\infty}
\varepsilon_k^{-1}e^{-c(s-\rho_{ki})/\varepsilon_k}s^3\,ds
&=
\int_{\rho_{ki}/\varepsilon_k}^{\infty}
e^{-cw}(\rho_{ki}+\varepsilon_k w)^3\,dw
\notag\\
&\le
\rho_{ki}^3\int_1^\infty e^{-cw}(1+w)^3\,dw
\le C\rho_{ki}^3,
\end{align}
and similarly,
\begin{align}\label{equation:ballmoments5}
\int_{2\rho_{ki}}^{\infty}
\varepsilon_k^{-1}e^{-c(s-\rho_{ki})/\varepsilon_k}s^5\,ds
\le C\rho_{ki}^5.
\end{align}
Since $K\setminus U_k\subset\mathbb R^4\setminus
B_{2\rho_{ki}}(X_{ki})$ for every $i$, by \eqref{equation:ballmoments3}, \eqref{equation:ballmoments5} and the fact that $\sum_i\rho_{ki}^3\le C\delta_k^2$, we have
\begin{equation}\label{e0.48}
\begin{aligned}
&\sum_i\int_{K\setminus U_k}
\varepsilon_k^{-1}
e^{-c\,\operatorname{dist}(x,B_{\rho_{ki}}(X_{ki}))/\varepsilon_k}
|t-\ell_k(y)|^2\,dx\\
&\quad\le \sum_i\int_{K\setminus U_k}
\varepsilon_k^{-1}
e^{-c(s-\rho_{ki})/\varepsilon_k}
(2s^2+8s_k^2) s^3ds\\
&\quad\le C\sum_i
\bigl(\rho_{ki}^5+s_k^2\rho_{ki}^3\bigr)
\le C\bigl((\rho_k^*)^2+s_k^2\bigr)\delta_k^2
=o(\delta_k^2).
\end{aligned}
\end{equation}
Moreover, \eqref{equation:epsdelta} gives
\begin{equation}\label{e0.49}
\frac{e^{-c/\varepsilon_k}}{\varepsilon_k\delta_k^2}
\int_K|t-\ell_k(y)|^2\,dx
\le C\varepsilon_k^{-3}e^{-c/\varepsilon_k}
\longrightarrow0.
\end{equation}
Together with \eqref{e0.46} and \eqref{e0.47}--\eqref{e0.49}, we obtain
\begin{equation}\label{e0.51}
 \int_{K\setminus U_k}e_k|t-\ell_k(y)|^2\,dx
 =o(\delta_k^2).
\end{equation}
Combining this with \eqref{e0.50}, we complete the proof of \eqref{equation:weightedheight}.

\end{proof}

\begin{proposition}\label{proposition:strongtilt}
Let $h$ be as in Lemma~\ref{lemma:harmonic}, and suppose that
$\delta_k>0,$
$ \varepsilon_k/\delta_k \to0.$
Along the same subsequence, for every compact $K\subset\subset Q_{1/2}$,
\begin{align}
\int_K\varepsilon_k
\bigl|D_Tu_k+\delta_k u_{k,4}Dh(y)\bigr|^2\,dx
&=o(\delta_k^2),
\label{equation:strongtilt}\\
\int_K
\left[\frac{W(u_k)}{\varepsilon_k}
-\varepsilon_k\Phi_k(Du_k)\right]\,dx
&=o(\delta_k^2).
\label{equation:strongdiscrep}
\end{align}
\end{proposition}

\begin{proof}
Set
\[
\mathcal R_k:=
\frac{W(u_k)}{\varepsilon_k}
-\varepsilon_k\Phi_k(Du_k)\ge0,
\qquad
S_k:=e_k-\varepsilon_k u_{k,4}\Phi_{k,4}(Du_k).
\]
Since $D\Phi_k(Du_k)\cdot Du_k=2\Phi_k(Du_k)$,
\begin{equation}\label{equation:Sdecomp}
S_k=\mathcal R_k
+\varepsilon_kD_T\Phi_k(Du_k)\cdot D_Tu_k.
\end{equation}
Choose \(0<r<R<1/2\) such that \(K\subset Q_r\), and choose
\(0\le\eta,\chi\le1\) satisfying
\[
\begin{gathered}
 \eta\in C_c^\infty(B'_R),
 \quad \eta=1\quad\text{on }\overline{B'_r};\qquad
 \chi\in C_c^\infty((-R,R)),
 \quad \chi=1\quad\text{on }[-r,r].
\end{gathered}
\]
Then \(\eta(y)\chi(t)=1\) on \(Q_r\), and
$
 \operatorname{spt}\eta\times\operatorname{spt}\chi'
 \subset
 \overline{B'_R}\times\bigl([-R,-r]\cup[r,R]\bigr).
$
By Lemma~\ref{lemma:onesheet}, we have
\[
 1-|u_k|\longrightarrow0
 \quad\text{uniformly on }
 \operatorname{spt}\eta\times\operatorname{spt}\chi'.
\]
Apply \eqref{equation:stress} to
\[
X=\eta(y)^2\chi(t)^2\bigl(t-\ell_k(y)\bigr)e_4,
\qquad
\ell_k(y)=a_k+\delta_kh(y),
\]
we obtain
\begin{equation}\label{e0.52}
\begin{split}
0
&=\int_{Q_{1/2}}
\left[
e_k\partial_4X_4
-\varepsilon_k u_{k,4}
\left(
D_T\Phi_k(Du_k)\cdot D_TX_4
+\Phi_{k,4}(Du_k)\partial_4X_4
\right)
\right]\,dx\\
&=\int_{Q_{1/2}}
\left[
S_k\partial_4X_4
-\varepsilon_k u_{k,4}
D_T\Phi_k(Du_k)\cdot D_TX_4
\right]\,dx.
\end{split}
\end{equation}
Since
\begin{equation}\label{e0.53}
\partial_4X_4
=\eta^2\chi^2
+2\eta^2\chi\chi'(t-\ell_k),\qquad
D_TX_4
=2\eta\chi^2(t-\ell_k)D\eta
-\eta^2\chi^2D\ell_k,
\end{equation}
substituting \eqref{e0.53} into \eqref{e0.52} and using \eqref{equation:Sdecomp}, we obtain
\begin{equation}\label{equation:relativeidentity}
\begin{split}
&\int_{Q_{1/2}}\eta^2\chi^2
\left[
\mathcal{R}_k
+\varepsilon_kD_T\Phi_k(Du_k)\cdot
\bigl(D_Tu_k+u_{k,4}D\ell_k\bigr)
\right]\,dx\\
&\quad=
2\int_{Q_{1/2}}
\varepsilon_k\eta\chi^2(t-\ell_k)u_{k,4}
D_T\Phi_k(Du_k)\cdot D\eta\,dx
-2\int_{Q_{1/2}}
\eta^2\chi\chi'(t-\ell_k)S_k\,dx.
\end{split}\end{equation}


Since
$
 |S_k|
 \le e_k+\varepsilon_k|u_{k,4}|
       |\Phi_{k,4}(Du_k)|
 \le Ce_k,$
$ \|t-\ell_k(y)\|_{L^\infty(Q_{1/2})}\le C,$
and \(|t|\ge r\) on
\(\operatorname{spt}\eta\times\operatorname{spt}\chi'\),
by Lemma~\ref{lemma:modicaclear} and $\varepsilon_k/\delta_k\to0$, we have
\[
\left|
\int_{Q_{1/2}}\eta^2\chi\chi'(t-\ell_k)S_k\,dx
\right|
\le C\varepsilon_k^{-1}e^{-c/\varepsilon_k}
=o(\delta_k^2).
\]
Using \eqref{equation:coercivity} and
$\varepsilon_k u_{k,4}^2\le2e_k$, Cauchy-Schwarz inequality yields
\[
\begin{aligned}
 &\left|
 \int_{Q_{1/2}}
 \varepsilon_k\eta\chi^2(t-\ell_k)u_{k,4}
 D_T\Phi_k(Du_k)\cdot D\eta\,dx
 \right|\\
 &\quad\le C
 \left(
 \int_{Q_{1/2}}\eta^2\chi^2
 \varepsilon_k u_{k,4}^2|t-\ell_k|^2\,dx
 \right)^{1/2}
 \left(
 \int_{Q_{1/2}}\chi^2|D\eta|^2
 \varepsilon_k|D_Tu_k|^2\,dx
 \right)^{1/2}\\
 &\quad\le C\delta_k
 \left(
 \int_{Q_{1/2}}\eta^2\chi^2e_k
 |t-\ell_k|^2\,dx
 \right)^{1/2}
 =o(\delta_k^2),
\end{aligned}
\]
where we used Lemma~\ref{lemma:weightedheight} in the last equality.


Let
\[
r_k:=D_Tu_k+\delta_k u_{k,4}Dh,
\qquad
q_k:=(-\delta_k u_{k,4}Dh,u_{k,4}).
\]
Uniform convexity gives
\begin{equation}\label{equation:relativecoercivity}
\bigl[D_T\Phi_k(Du_k)-D_T\Phi_k(q_k)\bigr]\cdot r_k
\ge\lambda|r_k|^2.
\end{equation}
Since \(D\ell_k=\delta_kDh\), we have
\(D_Tu_k+u_{k,4}D\ell_k=r_k\).
By \eqref{equation:relativeidentity}, we get
\begin{equation}\label{e0.54}
\int_{Q_{1/2}}\eta^2\chi^2
\bigl[
\mathcal{R}_k
+\varepsilon_kD_T\Phi_k(Du_k)\cdot r_k
\bigr]\,dx
=o(\delta_k^2).
\end{equation}
By \eqref{equation:relativecoercivity},
\[
\begin{aligned}
D_T\Phi_k(Du_k)\cdot r_k
&=
\bigl[D_T\Phi_k(Du_k)-D_T\Phi_k(q_k)\bigr]\cdot r_k
+D_T\Phi_k(q_k)\cdot r_k\\
&\ge \lambda|r_k|^2+D_T\Phi_k(q_k)\cdot r_k.
\end{aligned}\]
Substituting this into \eqref{e0.54}, we obtain
\begin{equation}\label{e0.55}
\int_{Q_{1/2}}\eta^2\chi^2
\bigl(\mathcal{R}_k+\lambda\varepsilon_k|r_k|^2\bigr)\,dx
\le
-\int_{Q_{1/2}}\eta^2\chi^2\varepsilon_k
D_T\Phi_k(q_k)\cdot r_k\,dx
+o(\delta_k^2).
\end{equation}
By   homogeneity and Taylor expansion,
\begin{equation}\label{equation:qexpansion}
\begin{aligned}
D_T\Phi_k(q_k)
&=u_{k,4}D_T\Phi_k(-\delta_kDh,1)\\
&=-\delta_k u_{k,4}B_kDh
+O(\delta_k^2|u_{k,4}|),
\end{aligned}
\end{equation}
where $B_k=D^2_{TT}\Phi_k(e_4)$.
Moreover,
\begin{equation}\label{e0.56}
\int_{Q_{1/2}}\varepsilon_k|r_k|^2\,dx
\le
2\delta_k^2
+2\delta_k^2\|Dh\|_{L^\infty(B'_{1/2})}^2
\int_{Q_{1/2}}\varepsilon_k u_{k,4}^2\,dx
\le C\delta_k^2.
\end{equation}
Consequently, by \eqref{equation:qexpansion} and Cauchy-Cauchy-Schwarz inequality, we obtain
\begin{equation}\label{e0.58}
\begin{split}
&\left|
\int_{Q_{1/2}}\eta^2\chi^2\varepsilon_k
\bigl[D_T\Phi_k(q_k)+\delta_k u_{k,4}B_kDh\bigr]
\cdot r_k\,dx
\right|\\
&\quad\leq C \int_{Q_{1/2}}\eta^2\chi^2\varepsilon_k\delta_k^2|u_{k,4}|| r_k|\,dx\\
&\quad\le C\delta_k^2
\left(\int_{Q_{1/2}}\varepsilon_k u_{k,4}^2\,dx\right)^{1/2}
\left(\int_{Q_{1/2}}\varepsilon_k|r_k|^2\,dx\right)^{1/2}
\le C\delta_k^3.
\end{split}
\end{equation}


\noindent By \eqref{equation:fluxlimit} and \eqref{equation:normaldensity},
\begin{equation}\label{e0.57}
\begin{split}
&\int_{Q_{1/2}}\eta^2\chi^2
\left[
\frac{\varepsilon_k u_{k,4}D_Tu_k}{\delta_k}\cdot B_kDh
+\varepsilon_k u_{k,4}^2\langle B_kDh,Dh\rangle
\right]\,dx\\
&\quad\longrightarrow
-c_W\int_{B'_{1/2}}\eta^2\langle B_\infty Dh,Dh\rangle\,dy
+c_W\int_{B'_{1/2}}\eta^2\langle B_\infty Dh,Dh\rangle\,dy
=0.
\end{split}
\end{equation}
Since $(\delta_k u_{k,4}B_kDh)\cdot r_k=\delta_k u_{k,4} \langle B_kDh, D_Tu_k \rangle + \delta_k^2 u_{k,4}^2 \langle B_kDh, Dh \rangle$, by \eqref{e0.58} and \eqref{e0.57}, we obtain
\[
\begin{aligned}
0
&\le\frac1{\delta_k^2}
\int_K\bigl(\mathcal{R}_k+\lambda\varepsilon_k|r_k|^2\bigr)\,dx\\
&\le\frac1{\delta_k^2}
\int_{Q_{1/2}}\eta^2\chi^2
\bigl(\mathcal{R}_k+\lambda\varepsilon_k|r_k|^2\bigr)\,dx
\longrightarrow0.
\end{aligned}
\]
This proves \eqref{equation:strongtilt} and
\eqref{equation:strongdiscrep}.

\end{proof}

\section{Improvement of Flatness}\label{section:onestep}

We are now in a position to prove the following improvement
of flatness.

\begin{proposition}\label{theorem:onestep}
Let \(u\) be a global minimizer with \(-1<u<1\) and \(u(0)=0\).  There exist constants $0<\theta<1/128$, $C_N,K,\delta_*$, $\varepsilon_*>0$, \(r>0\) and \(\nu\in\mathbb S^3\) such that \begin{equation}\label{equation:stepinputs}
r^{-1}\le\varepsilon_*,
\qquad
\mathcal T_r(\nu)\le\delta_*^2,
\qquad
\mathcal M_r(\nu)\le\frac32c_W|B'_4|,
\end{equation}
where $\mathcal T_r(\nu)$ and $\mathcal M_r(\nu)$ are defined by \eqref{equation:TMphysical}.
If \( \mathcal T_r(\nu) \ge K^2 r^{-2}\), then there exists \(\mu\in\mathbb S^3\) such that
\begin{equation}\label{equation:highstep}
|\mu-\nu|\le C_N\sqrt{\mathcal T_r(\nu)},
\qquad
\mathcal T_{\theta r}(\mu)\le \frac{1}{4} {\mathcal T_r(\nu)},
\qquad
\mathcal M_{\theta r}(\mu)\le\frac32c_W|B'_4|.
\end{equation}
If \( \mathcal T_r(\nu) <K^2r^{-2}\), then \(\mu=\nu\) satisfies
\begin{equation}\label{equation:lowstep}
\mathcal T_{\theta r}(\mu)
\le K^2\theta^{-1}(\theta r)^{-2},
\qquad
\mathcal M_{\theta r}(\mu)\le\frac32c_W|B'_4|.
\end{equation}
In both cases, $Q_{4\theta r}^{\mu}\subset Q_{4r}^{\nu}.$
\end{proposition}


\begin{proof}

For \(r>0\), define $v(z):=u(rL_\nu z)$ with $r=\varepsilon^{-1}$, where $z\in Q_4$ as in Section \ref{section:normalization}. Then, from \eqref{e0.19}, \eqref{equation:TMnormalized1}
and the assumptions \eqref{equation:stepinputs}, we have
\[
\int_{Q_4}\varepsilon|D_Tv|^2\,dz=\mathcal T_r(\nu)\le\delta_*^2,
\qquad
\mathcal E_{\varepsilon,\widehat F_\nu}(v;Q_4)=\mathcal M_r(\nu)\le\frac32c_W|B'_4|,
\]
which satisfies the hypotheses of Lemma~\ref{lemma:onesheet} with $\tau=1/2$.
By Lemmas~\ref{lemma:heightH1} and Lemma~\ref{lemma:harmonic},
every height limit \(h\) satisfies
\begin{equation}\label{e0.60}
|Dh(0)|+\|D^2h\|_{L^\infty(B'_{1/4})}\le C_H,
\end{equation}
where \(C_H\)   is independent of the chosen sequence.
Let \(C_0\) be a uniform constant such that
$
L_\nu^{-1}L_\mu(Q_{4\theta})\subset Q_{C_0\theta}.
$
Choose \(C_N\ge C_0C_H\), and then choose \(\theta\) such that
\[
0<\theta<\frac1{128},
\qquad
C_0\theta<\frac14,
\qquad
C_1\theta^2<\frac18,
\]
where \(C_1>0\) is the constant in
\eqref{equation:stepcontradiction}, independent of \(\theta\)
and of the sequence.


Now we prove  \eqref{equation:highstep}.
Assume by contradiction that for the fixed \(C_N,\theta\) and for any $k>0$, $\varepsilon_k=r_k^{-1}$  and \(\delta_*=\varepsilon_*=1/k\), if
\[
\delta_k:=\sqrt{\mathcal T_{r_k}(\nu_k)}\le\frac1k,
\qquad
\varepsilon_k\le\frac1k,
\qquad
\delta_k\ge k\varepsilon_k,\qquad
\mathcal M_{r_k}(\nu_k)\le\frac32c_W|B'_4|,
\]
then
 \begin{equation}\label{e0.59}\mathcal T_{\theta r_k}(\mu_k)\ge \frac{1}{4} {\mathcal T_{r_k}(\nu_k)}\quad\text{ for any }  \mu_k\in\mathbb{S}^n\text{ with } |\mu_k-\nu_k|\leq C_N\delta_k.\end{equation}
The normalized functions \(v_k(z)=u_k(r_kL_{\nu_k}z)\)
satisfy \eqref{equation:oneinputs} with \(\tau=1/2\), and
$
 \varepsilon_k/\delta_k \le 1/k\to0.
$
By Lemma~\ref{lemma:heightH1} and
Proposition~\ref{proposition:strongtilt}, after passing to
a subsequence, there is a uniform  limit \(h\) such that
\[
\int_\Omega\varepsilon_k
\bigl|D_Tv_k+\delta_kv_{k,4}Dh\bigr|^2\,dz
=o(\delta_k^2)
\]
for every compact subcylinder $\Omega\subset\subset Q_{1/2}.$
Define
\begin{equation}\label{equation:newnormal}
b_k=\delta_kDh(0)\qquad\text{and}\qquad\mu_k=
\frac{L_{\nu_k}^{-T}(-b_k,1)}
     {|L_{\nu_k}^{-T}(-b_k,1)|}.
\end{equation}
Since \(L_{\nu_k}^{-T}e_4=\nu_k/F(\nu_k)\),
the uniform bounds for \(L_{\nu_k}^{\pm1}\) and \eqref{e0.60} imply
\begin{equation}\label{e0.61}
|\mu_k-\nu_k|
\le
\frac{2|L_{\nu_k}^{-T}(-b_k,0)|}
     {|L_{\nu_k}^{-T}e_4|}
\le C_0|b_k|
\le C_0C_H\delta_k
\le C_N\delta_k.
\end{equation}
Then, by Lemma~\ref{lemma:masstransfer}, we have
\begin{equation}\label{e0.62}
\mathcal M(\theta r_k,\mu_k)\le\frac32c_W|B'_4|
\end{equation}
for all sufficiently large \(k\).



For \(x=r_kL_{\nu_k}z\), the chain rule gives
\[
D_xu_k(x)=r_k^{-1}L_{\nu_k}^{-T}D_zv_k(z),\qquad\text{and}\qquad D_zv_k
=v_{k,4}(-b_k,1)
+\bigl(D_Tv_k+b_kv_{k,4},0\bigr).
\]
By \eqref{equation:newnormal},
$
\Pi_{\mu_k^\perp}L_{\nu_k}^{-T}(-b_k,1)=0.
$
Hence the uniform bound for \(L_{\nu_k}^{-1}\) yields
\begin{equation}\label{equation:physicalprojection}
\begin{aligned}
\bigl|\Pi_{\mu_k^\perp}D_xu_k(x)\bigr|
&=r_k^{-1}
\left|
\Pi_{\mu_k^\perp}L_{\nu_k}^{-T}
\bigl(D_Tv_k(z)+b_kv_{k,4}(z),0\bigr)
\right|\\
&\le Cr_k^{-1}
\bigl|D_Tv_k(z)+b_kv_{k,4}(z)\bigr|.
\end{aligned}
\end{equation}


\noindent Let
$
\Omega_k=L_{\nu_k}^{-1}L_{\mu_k}(Q_{4\theta})
\subset Q_{C_0\theta}\subset\subset Q_{1/4}.
$
By \eqref{equation:TMphysical} and
\eqref{equation:physicalprojection},
\[
\mathcal T_{\theta r_k}(\mu_k)
\le C\theta^{-3}
\int_{Q_{C_0\theta}}\varepsilon_k
|D_Tv_k+b_kv_{k,4}|^2\,dz.
\]
Since \(b_k=\delta_kDh(0)\), we have
\[
D_Tv_k+b_kv_{k,4}
 =D_Tv_k+\delta_kv_{k,4}Dh(y)
+\delta_kv_{k,4}\bigl(Dh(0)-Dh(y)\bigr).
\]
By mean value theorem and \eqref{e0.60}, $|Dh(y)-Dh(0)|\le C_H|y|\le C_0C_H\theta$
 on $B'_{C_0\theta}.$
Consequently,
\[
\frac{\mathcal T_{\theta r_k}(\mu_k)}{\delta_k^2}
 \le
\frac{C\theta^{-3}}{\delta_k^2}
\int_{Q_{C_0\theta}}\varepsilon_k
|D_Tv_k+\delta_kv_{k,4}Dh|^2\,dz
 +
C\theta^{-1}
\int_{Q_{C_0\theta}}\varepsilon_kv_{k,4}^2\,dz.
\]
For the fixed \(\theta\), the first term tends to zero by
Proposition~\ref{proposition:strongtilt}, while
\eqref{equation:normaldensity} gives
\[
\int_{Q_{C_0\theta}}\varepsilon_kv_{k,4}^2\,dz
\longrightarrow c_W|B'_{C_0\theta}|
\le C\theta^3.
\]
Hence
\begin{equation}\label{equation:stepcontradiction}
\limsup_{k\to\infty}
\frac{\mathcal T_{\theta r_k}( \mu_k)}{\delta_k^2}
\le C_1\theta^2<\frac18,
\end{equation}
which is a contradiction with \eqref{e0.59}. Therefore, \eqref{equation:highstep} holds.


Decrease \(\delta_*,\varepsilon_*\), if necessary, so that
Lemma~\ref{lemma:masstransfer} applies with the fixed
\(C_N,\theta\) and \(\tau=1/2\).
Then every \(\mu\in\mathbb S^3\) satisfying
\(|\mu-\nu|\le C_N\sqrt{\mathcal T_r(\nu)}\), we get
\[
\mathcal M_{\theta r}(\mu)\le\frac32c_W|B'_4|.
\]
Moreover, the choice of \(\theta\) gives
\[
Q_{4\theta r}^{\mu}
 =rL_\nu\bigl(L_\nu^{-1}L_\mu(Q_{4\theta})\bigr)
 \subset rL_\nu(Q_{C_0\theta})
\subset Q_{4r}^{\nu}.
\]


If \(\mathcal T_r(\nu) <K^2r^{-2}\), take \(\mu=\nu\).
Since \(Q_{4\theta r}^{\nu}\subset Q_{4r}^{\nu}\),
\eqref{equation:TMphysical} gives
\[
\mathcal T_{\theta r}( \nu)
 \le\theta^{-3}\mathcal T(r,\nu)
 <K^2\theta^{-3}r^{-2}
=K^2\theta^{-1}(\theta r)^{-2}.
\]
Apply Lemma~\ref{lemma:masstransfer}  with \(\mu=\nu\), we obtain
\[
\mathcal M_{\theta r}( \nu)\le\frac32c_W|B'_4|.
\]
This proves \eqref{equation:lowstep} and completes the proof.


\end{proof}

\section{Proof of Theorem \ref{theorem:main}}\label{section:iteration}


With the constants in Proposition~\ref{theorem:onestep} fixed, set
\begin{equation}\label{equation:floorconstants}
C_0=\frac{K^2}{\theta},
\qquad
C_f=\frac{C_0}{1-\theta^2/4},
\end{equation}
and choose
\begin{equation}\label{equation:Rstar}
R_*\ge
\max\left\{
\varepsilon_*^{-1},
\frac{\sqrt{2C_f}}{\delta_*},
1
\right\}.
\end{equation}

\begin{proposition}\label{proposition:iteration}
Let \(u\) be a  global minimizer with
\(-1<u<1\) and \(u(0)=0\).
Suppose \(R\ge R_*\) and \(\nu_0\in\mathbb S^3\) satisfy
\[
 \mathcal T_R(\nu_0)\le\frac{\delta_*^2}{2},
\qquad
\mathcal M_R(\nu_0)\le\frac32c_W|B'_4|.
\]
Set
$
r_j=\theta^jR,$ and
$N=\max\{j\in\mathbb N_0:r_j\ge R_*\}.$
Then, there exists  normals \(\nu_1,\ldots,\nu_N\)  such that
\begin{equation}\label{equation:recursionsolved}
 \mathcal T_{r_j}(\nu_j)\le 4^{-j}\mathcal T_R(\nu_0)+C_fr_j^{-2},\qquad \mathcal M_{r_j}(\nu_j)\le\frac32c_W|B'_4|
\qquad\text{for}\quad0\le j\le N,
\end{equation}
and
\begin{equation}\label{equation:normalsum}
|\nu_N-\nu_0|\le C\sqrt{\mathcal T_R(\nu_0)}+\frac{C}{r_N},
\end{equation}
where \(C\) depends only on the constants in
Proposition~\ref{theorem:onestep}.
\end{proposition}


\begin{proof}
We construct the normals by induction. The case \(j=0\) is trivial.
Suppose that \eqref{equation:recursionsolved} hold for some \(0\le j<N\).
Since \(r_j\ge R_*\), \eqref{equation:Rstar} gives
\[
\mathcal T_{r_j}(\nu_j)\le 4^{-j}\mathcal T_{R}(\nu_0)+C_fr_j^{-2}
\le\frac{\delta_*^2}{2}+C_fR_*^{-2}
\le\delta_*^2.
\]
Proposition~\ref{theorem:onestep} therefore yields \(\nu_{j+1}\) with
\begin{equation}\label{equation:recurrence}
\mathcal T_{r_{j+1}}(\nu_{j+1})\le \frac{1}{4}\mathcal T_{r_j}(\nu_j)+C_0r_{j+1}^{-2},
\quad\quad \mathcal M(r_{j+1},\nu_{j+1})
\le\frac32c_W|B'_4|,
\end{equation}
and
\begin{equation}\label{e0.63}|\nu_{j+1}-\nu_j|\le C_N\sqrt{\mathcal T_{r_j}(\nu_j)}.\end{equation}
Using \(r_{j+1}=\theta r_j\) and
\(\frac{1}{4}\theta^2C_f+C_0=C_f\), we obtain
\[
\begin{aligned}
\mathcal T_{r_{j+1}}(\nu_{j+1})
&\le \frac{1}{4}\bigl(4^{-j}E_0+C_fr_j^{-2}\bigr)
+C_0r_{j+1}^{-2}\\
&=4^{-j-1}\mathcal T_{R}(\nu_0)+(a\theta^2C_f+C_0)r_{j+1}^{-2}\\
&=4^{-j-1}\mathcal T_{R}(\nu_0)+C_fr_{j+1}^{-2}.
\end{aligned}
\]
This proves \eqref{equation:recursionsolved}
for \(0\le j\le N\). Then, by \eqref{e0.63} and
\eqref{equation:recursionsolved}, we obtain
\[
\begin{aligned}
|\nu_N-\nu_0|
&\le C_N\sum_{j=0}^{N-1}\sqrt{\mathcal T_{r_j}(\nu_j)}\\
&\le C_N\sqrt{\mathcal T_{R}(\nu_0)}\sum_{j=0}^{N-1}4^{-j/2}
+C_N\sqrt{C_f}\sum_{j=0}^{N-1}r_j^{-1}\\
&=C_N\sqrt{\mathcal T_{R}(\nu_0)}\sum_{j=0}^{N-1}4^{-j/2}
+\frac{C_N\sqrt{C_f}}{r_N}
\sum_{j=0}^{N-1}\theta^{N-j}\\
&\le2 C_N \sqrt{\mathcal T_{R}(\nu_0)}
+\frac{C_N\sqrt{C_f}\,\theta}{1-\theta}\,\frac1{r_N}.
\end{aligned}
\]
This proves \eqref{equation:normalsum}.


\end{proof}

\begin{proposition}\label{proposition:fixedplane}
Let \(u\in C^2(\mathbb{R}^4)\) be a solution of \eqref{equation:mainPDE} satisfying \eqref{equation:layer}.
Then there exist \(P\in G(3)\) and constants \(C, r_*>0\) such that
\begin{equation}\label{equation:transversegrowth}
\int_{B_r}|\Pi_PDu|^2\,dx\le Cr
\qquad\text{for all }r\ge r_*.
\end{equation}
\end{proposition}


\begin{proof}
 Choose \(R_k\to\infty\) and a planar limit with unit normal
\(\nu_0\) as in Corollary~\ref{corollary:planarinit}. For \(z \in Q_5\), set $v_k(z)=u(R_kL_{\nu_0}z)$   with $\varepsilon_k=R_k^{-1}.$
By Proposition~\ref{proposition:phasecompact},
\eqref{equation:TMnormalized1} and
\(\widehat F_{\nu_0}(e_4)=\widehat F_{\nu_0}(-e_4)=1\),
\begin{equation}\label{e0.65}
\mathcal M_{R_k}(\nu_0)
=\mathcal E_{R_k^{-1},\widehat F_{\nu_0}}(v_k;Q_4)
\longrightarrow
c_W\int_{B'_4}\widehat F_{\nu_0}(e_4)\,dy
=c_W|B'_4|.
\end{equation}
Recall
\[
\Psi(s):=\int_0^s\sqrt{2W(a)}\,da.
\]
Since \(-1<v_k<1\) and
\(v_k\to\sigma\operatorname{sgn}(z_4)\) in \(L^1(Q_4)\), we have $\Psi(v_k)\to\Psi(\sigma\operatorname{sgn}(z_4))$ in $L^1(Q_4).$
By lower semicontinuity, we obtain
\begin{equation}\label{e0.67}
\begin{split}
\liminf_{k\to\infty}
\int_{Q_4}\sqrt{2W(v_k)}\,|v_{k,4}|\,dz
&=\liminf_{k\to\infty}|D_4\Psi(v_k)|(Q_4)\\
&\ge
|D_4[\Psi(\sigma\operatorname{sgn}(z_4))]|(Q_4)\\
&=|\Psi(1)-\Psi(-1)|\,|B'_4|=c_W|B'_4|.
\end{split}
\end{equation}
It follows from \eqref{equation:PhiSplit} that
\begin{equation}\label{e0.66}
0\le\frac{\lambda}{2} \mathcal T_{R_k}(\nu_0)
\le
\mathcal M_{R_k}(\nu_0)
-\int_{Q_4}\sqrt{2W(v_k)}\,|v_{k,4}|\,dz
\longrightarrow0.
\end{equation}


By the uniform bounds for \(L_\nu^{-1}\), choose \(A_0\ge1\) such that
$
B_r\subset Q_{4A_0r}^{\nu}$ for every $r>0,\ \nu\in\mathbb S^3.$
Fix \(r\) with \(A_0r\ge R_*\).
For all sufficiently large \(k\), define
$
r_{k,j}=\theta^jR_k,$ and
$N_k=\max\{j\in\mathbb N_0:r_{k,j}\ge A_0r\}.$
Then
\begin{equation}\label{equation:terminalradius}
A_0r\le r_{k,N_k}<\frac{A_0r}{\theta},
\end{equation}
and $r_{k,j+1}\ge r_{k,N_k}\ge A_0r\ge R_*$ for every \(0\le j<N_k\).


Set \(P=\nu_0^\perp\) and \(\mu_k=\nu_{k,N_k}\). Since
$
\|\Pi_P-\Pi_{\mu_k^\perp}\|
\le C|\mu_k-\nu_0|,$
we have
\[
|\Pi_Pp|^2
\le 2|\Pi_{\mu_k^\perp}p|^2
   +C|\mu_k-\nu_0|^2|p|^2]\qquad \text{for any}\quad p\in\mathbb{R}^4.
\]
Therefore, by Proposition~\ref{proposition:iteration},   we obtain 
\[
\begin{aligned}
\int_{B_r}|\Pi_PDu|^2\,dx
&\le
2\int_{Q_{4r_{k,N_k}}^{\mu_k}}
|\Pi_{\mu_k^\perp}Du|^2\,dx
+C|\mu_k-\nu_0|^2\int_{B_r}|Du|^2\,dx\\
&\le Cr_{k,N_k}^3\mathcal{T}_{r_{k,N_k}}(\mu_k)
   +Cr^3|\mu_k-\nu_0|^2\\
&\le Cr_{k,N_k}^34^{-N_k}\mathcal{T}_{R_k}(\nu_0)
   +Cr_{k,N_k}
   +Cr^3\bigl(\mathcal{T}_{R_k}(\nu_0)+r_{k,N_k}^{-2}\bigr).
\end{aligned}
\]
By \eqref{equation:terminalradius} and \eqref{e0.66},
 for each fixed  large \(r\), we have
 \[
\int_{B_r}|\Pi_PDu|^2\,dx
\le Cr^3(4^{-N_k}+1)\mathcal{T}_{R_k}(\nu_0)+Cr
\longrightarrow Cr
\qquad\text{as}\quad k\to\infty.
\]
Hence \eqref{equation:transversegrowth} holds with the same
\(P=\nu_0^\perp\), where \(C\) is independent of \(r\) and \(k\).

\end{proof}

To conclude the proof, we  apply the quotient argument of Ambrosio--Cabr\'e~\cite{AC}
to derivatives tangent to \(P\).

\begin{lemma}\label{lemma:rigidity}
Let \(u\) solve \eqref{equation:mainPDE} in \(\mathbb R^4\),
with \(\partial_{x_4} u > 0\).
If \(P\in G(3)\) satisfies
\begin{equation}\label{equation:subquadratic}
\lim_{R\to\infty}\frac{1}{R^2}\int_{B_R}|\Pi_PDu|^2\,dx=0,
\end{equation}
then $\Pi_PDu\equiv0$ in $\mathbb R^4.$
\end{lemma}


\begin{proof}
Since \(\partial_{x_4} u >0\), we have \(Du\ne0\). Define \(v=\partial_{x_4} u\).
Differentiating \eqref{equation:mainPDE} in the directions \(x_4\) and \(\tau\in P\) gives
\begin{equation}\label{equation:Jacobipair}
\operatorname{div}\bigl(D^2\Phi(Du)\,Dv\bigr)=W''(u)v,
\qquad
\operatorname{div}\bigl(D^2\Phi(Du)\,D u_\tau\bigr)=W''(u)u_\tau.
\end{equation}

For \(s_\tau=u_\tau/v\), the symmetry of \(D^2\Phi(Du)\) and \eqref{equation:Jacobipair} give
\[
v^2D^2\Phi(Du)\,Ds_\tau
=vD^2\Phi(Du)\,Du_\tau-u_\tau D^2\Phi(Du)\,Dv,
\]
and hence
\begin{equation}\label{equation:quotient}
\operatorname{div}\bigl(v^2D^2\Phi(Du)\,Ds_\tau\bigr)
=v\operatorname{div}\bigl(D^2\Phi(Du)\,Du_\tau\bigr)
-u_\tau\operatorname{div}\bigl(D^2\Phi(Du)\,Dv\bigr)=0.
\end{equation}
Choose \(\eta_R\in C_c^\infty(B_{2R})\) satisfying
\[
0\le\eta_R\le1,\qquad
\eta_R=1\ \text{on }B_R,\qquad
|D\eta_R|\le C/R.
\]
Since \(\min_{\overline{B_{2R}}}v>0\), we   test
\eqref{equation:quotient} with \(\eta_R^2s_\tau\).
Uniform ellipticity and Young's inequality yield
\[
\begin{aligned}
\lambda\int_{\mathbb R^4}\eta_R^2v^2|Ds_\tau|^2\,dx
&\le
-2\int_{\mathbb R^4}
\eta_Rs_\tau v^2\bigl\langle D^2\Phi(Du)\,Ds_\tau,D\eta_R\bigr\rangle\,dx\\
&\le
\frac{\lambda}{2}
\int_{\mathbb R^4}\eta_R^2v^2|Ds_\tau|^2\,dx
+C\int_{\mathbb R^4}u_\tau^2|D\eta_R|^2\,dx.
\end{aligned}
\]
Thus, by \eqref{equation:subquadratic}, for every fixed \(\rho>0\) and \(R>\rho\),
\[
\int_{B_\rho}v^2|Ds_\tau|^2\,dx\le \frac{C}{R^2}\int_{B_{2R}}u_\tau^2\,dx
\le \frac{C|\tau|^2}{R^2}
\int_{B_{2R}}|\Pi_PDu|^2\,dx
\longrightarrow0
\]
  Since \(v>0\), we obtain
\(Ds_\tau\equiv0\), and therefore
\begin{equation}\label{equation:derivativeratios}
u_\tau=c(\tau)v
\qquad\text{for any}\quad\tau\in P,
\end{equation}
where \(c(\tau)=u_\tau(0)/v(0)\) is linear on \(P\).


It remains to show that \(c(\tau)=0\). Suppose that \(c(\tau)\not\equiv 0\)  and set
\[
K:=\ker c\subset P,
\]
then $\dim K=\dim P-1=2.$ By \eqref{equation:derivativeratios},
\begin{equation}\label{e0.69}
u_\tau=c(\tau)v=0
\qquad\text{for every }\tau\in K.
\end{equation}
Hence, for \(x\in\mathbb R^4\) and \(y\in K\),
\[
u(x+y)-u(x)=\int_0^1 Du(x+sy)\cdot y\,ds=0.
\]
Choose \(\tau_0\in P\) with \(c(\tau_0)\ne0\). Since
$
v= u_{\tau_0}/c(\tau_0),
$
we also have \(v(x+y)=v(x)\) for every \(y\in K\). Set $\bar v(z):=v(y+z)$ for any $y\in K,\ z\in K^\perp.$
It follows from \eqref{equation:subquadratic} that
\begin{equation}\label{e0.68}
\lim_{R\to\infty} \frac{1}{R^2} \int_{B_R} v^2\,dx
= \lim_{R\to\infty} \frac{1}{|c(\tau_0)|^2 R^2} \int_{B_R} |u_{\tau_0}|^2\,dx
= 0.
\end{equation}
On the other hand, let \(B_{R/2}^K := K \cap B_{R/2}\) and \(B_{1}^{K^\perp} := K^\perp \cap B_1\). For \(R\ge4\) and \((y,z)\in B_{R/2}^K\times B_1^{K^\perp}\subset B_R\)  since \(v>0\), we obtain
\[
\int_{B_R}v^2\,dx
\ge
\int_{B_{R/2}^K}\int_{B_1^{K^\perp}}
\bar v(z)^2\,dz\,dy
=|B_{R/2}^K|\int_{B_1^{K^\perp}}\bar v(z)^2\,dz
=\frac{\pi  }{4}R^2\int_{B_1^{K^\perp}}\bar v(z)^2\,dz>0.
\]
Consequently,
\[
0<\frac{\pi  }{4}R^2\int_{B_1^{K^\perp}}\bar v(z)^2\,dz
\le\liminf_{R\to\infty}
\frac1{R^2}\int_{B_R}v^2\,dx=0,
\]
a contradiction. Thus $c=0$,   and hence $u_\tau=0$ for every $\tau\in P.$
Since \(\partial_{x_4} u>0\), we conclude that \(e_4\notin P\).

\end{proof}

\begin{proof}[Proof of Theorem~\ref{theorem:main}]
By Proposition~\ref{proposition:fixedplane}, there exists
\(P\in G(3)\) such that
\[
\int_{B_R}|\Pi_PDu|^2\,dx\le CR,
\]
which implies that \[\lim_{R\to\infty}\frac{1}{R^2}\int_{B_R}|\Pi_PDu|^2\,dx=0.\]
Then, it follows from Lemma~\ref{lemma:rigidity} that \(\Pi_PDu\equiv0\).
Hence, for a unit normal \(\nu\) to \(P\),
\begin{equation}\label{e0.70}
u(x)=h(x\cdot\nu-a_0)
\qquad\text{for}\quad x\in\mathbb R^4 ,
\end{equation}
where \(h:\mathbb R\to(-1,1)\) and \(a_0\in\mathbb R\).

Orient \(\nu\) so that \(e_4\cdot\nu>0\). Since \(\partial_{x_4} u=(e_4\cdot\nu)h'>0\), it follows that \(h'>0\). By homogeneity,
\[
D\Phi(h'\nu)=h'D\Phi(\nu),
\qquad
D\Phi(\nu)\cdot\nu=F(\nu)^2.
\]
Equation~\eqref{equation:mainPDE} therefore becomes
\begin{equation}\label{equation:finalODE}
F(\nu)^2h''=W'(h).
\end{equation}

The condition \eqref{equation:layer} implies that \(h(\pm\infty)=\pm1\).
Since \(h'>0\) and \(\int_{\mathbb R}h'\,dt=2\), there exist
\(t_j\to\infty\) with \(h'(t_j)\to0\).
Multiplying \eqref{equation:finalODE} by \(h'\) gives
\[
\frac{F(\nu)^2}{2}h'^2-W(h)\equiv
\lim_{j\to\infty}
\left[\frac{F(\nu)^2}{2}h'(t_j)^2-W(h(t_j))\right]=0.
\]
Consequently,
\begin{equation}\label{e0.71}
h'=\frac{\sqrt{2W(h)}}{F(\nu)}.
\end{equation}
Let \(g\) be the increasing one-dimensional solution of
\eqref{e0.2} on \(\mathbb R\), with \(g(0)=0\).
Choose \(t_0\in\mathbb R\) such that \(h(t_0)=0\).
Since \(g'=\sqrt{2W(g)}\), the functions \(h(t)\) and
\(g((t-t_0)/F(\nu))\) satisfy \eqref{e0.71}
with the same initial value at \(t=t_0\).
 By uniqueness for the first-order ODE \eqref{e0.71}, we have
\[
h(t)=g\!\left(\frac{t-t_0}{F(\nu)}\right).
\]

Substituting this into \eqref{e0.70} and setting
\(a:=a_0+t_0\), we conclude
\[
u(x)=h(x\cdot\nu-a_0)
=g\!\left(\frac{x\cdot\nu-a}{F(\nu)}\right),
\]
which completes the proof.

\end{proof}



\end{document}